\documentclass[twoside,12pt]{article} 

\usepackage[a4paper]{geometry} 
\usepackage{indentfirst}
\usepackage{amsfonts}
\usepackage{amsmath,amssymb}
\usepackage{amsthm}
\usepackage{mathrsfs,upref} 
\usepackage{graphicx}
\usepackage{epsfig}
\usepackage[pdftex]{hyperref}
\hypersetup{
            colorlinks=true,
            linkcolor=blue,
            anchorcolor=green,
            citecolor=red,
            bookmarksnumbered,
            bookmarksopen,
            bookmarksopenlevel=2
            }
\numberwithin{equation}{section}

\newtheorem {definition}{Definition}[section]

\newtheorem {theorem}   [definition]{Theorem}
\newtheorem {corollary} [definition]{Corollary}
\newtheorem {lemma}     [definition]{Lemma}
\newtheorem {notation}  [definition]{Notation}

\newtheorem {proposition}{Proposition}  [section]
\newtheorem {remark}{Remark}            [section]
\newtheorem {illustration}{Illustration}[section]
\newtheorem {tables}{Table}             [section]

\newcommand{\sldef}[1]{\begin{definition}{\textsl{#1}}\end{definition}}
\newcommand{\slthm}[1]{\begin{theorem}{\textsl{#1}}\end{theorem}}

\newcommand{\slcor}[1]{\begin{corollary}{\textsl{#1}}\end{corollary}}
\newcommand{\sllem}[1]{\begin{lemma}{\textsl{#1}}\end{lemma}}
\newcommand{\slnot}[1]{\begin{notation}{\textsl{#1}}\end{notation}}
\newcommand{\slprop}[1]{\begin{proposition}{\textsl{#1}}\end{proposition}}
\newcommand{\slrem}[1]{\begin{remark}{\textsl{#1}}\end{remark}}

\newcommand{\slprf}[1]{\begin{proof}{\textsl{#1}}\end{proof}} 

\newcommand{\isomeanvalue}[3]{{\overline {{\displaystyle #1}_{\scriptstyle #2}}}\bigl|_{#3}}
\newcommand{\isomeanvalueII}[3]{{\overline {{\displaystyle #1}_{\scriptstyle #2}}}\bigl|_{#3}^{II}}  
\newcommand{\isomeanvalueIII}[3]{{\overline {{\displaystyle #1}_{\scriptstyle #2}}}\bigl|_{#3}^{III}}
\newcommand{\tabincell}[2]{\begin{tabular}{@{}#1@{}}#2\end{tabular}}

\hfuzz=\maxdimen
\title{\Large{\textbf{AN INTRODUCTION TO \\ISOMORPHIC MATHEMATICAL ANALYSIS SYSTEM}}}
\author {Tim Y. Liu
        \thanks{
                Finished Date: ~2014-Jan-28 \hspace*{2.8em}Revised Date: ~2024-Jan-29\newline 
                \hspace*{1.8em}Liu Yuan(Tim), Jiangsu Province, China \newline
                \hspace*{1.8em}Field of Research: Mathematical Analysis \newline
                \hspace*{1.8em}Email: liu\_yuan@aliyun.com
        } \\
        {\small Hengbao Co.,Ltd.}
}
\date {\small Rev. Jan. 2024}

\begin{document}
\maketitle

\begin{abstract}
This paper is to build a primitive framework for a new possible extended system of
real mathematical analysis - the Isomorphic Mathematical Analysis System (IMAS). It
is based on some new concepts: e.g. isomorphic frame, dual-variable-isomorphic
function, and isomorphic coordinate system. More concepts are introduced to
constitute the whole IMAS framework. This IMAS attempts to provide current MA with new implements, which include isomorphic frame, thereby effectively putting the MA to work as well in some unevenly distributed coordinate spaces, realized by the isomorphic coordinate system. With these implements, many existing MA concepts are extended and incorporated into the new IMAS.
\vspace{0.1cm}

\\[0.5em]
\end{abstract}

\section{Introduction}

In classical real Mathematical Analysis (MA) system, the number marked at a point on a number-axis actually is the distance from the axis 0 point to it. But in this paper we are going to introduce a special variant of MA system which bases on coordinate axes that correlates the numbers on the axis with their positions by monotone bijections, refined as the concept of ``isomorphic frame'', whereby various concepts featured in the new system have certain forms of ``isomorphism'' with respect to their counterparts in the classical MA. This idea is analogous with the concept of ``isomorphism'' in the realm of Algebra, therefore the item ``isomorphic'' is borrowed in this work. \par

Corresponding alterations are made to definitions of concepts related to number, variable and function etc. With these changes, the original MA system sees an opportunity to upgrade to a higher ``extended MA'' with a great deal of new MA concepts and math problems. It also integrates lots of existing ones into a higher and nicer hierarchy, including the fact that the classical MA is the simplest case of the new MA system with identity bijections are considered. But in a broader sense, such extension of MA is also an integrated and complementary part of current MA system. The proposed new MA is referred to as Isomorphic Mathematical Analysis System (IMAS) in this paper, which framework must be built from scratch hereinafter.\par

The reference article \cite{LIUY2}, on the ``isomorphic means'', is dedicated for a specific topic belonging to this IMAS and some parts of it are incorporated into this work. A simple example system of IMAS is attached as the last section(Section \ref{sec:Example}) for better understanding of it. Typical extended concepts generated by logarithmic mappings are summarized in Table \ref{tab:GeometricGenSysSummary} for references.


\section{The Basics of IMAS}

There are 5 basic concepts to be introduced. Namely ``Isomorphic frame'', ``Isomorphic number and Isomorphic variable'', ``Dual-variable-isomorphic function'', ``Isomorphic number-axis'', and ``Dual-isomorphic rectangular coordinate system (pair)''.

\subsection{Isomorphic frame}
\subsubsection{Definition}
\sllem{\label{lem:OrderSetbeBJ}Let $X_1,...,X_n,U_1,...,U_n\subseteq\mathbb{R}$, ~$X=X_1\times ...\times X_n$, $~U=U_1\times ...\times U_n$ and $g_i\colon X_i\to U_i(i=1,...,n)$ be $n$ bijections. The ordered set $\{g_1,...,g_n\}$ that can map $\forall x=(x_1,...,x_n)\in X$ to $u=(g_1(x_1),...,g_n(x_n))\in U$ is a bijection, if it is treated as a function with $X$ being its domain and $U$ being the range(the image).}

\slprf{Firstly $\forall x=(x_1,...,x_n)\in X, ~y=(y_1,...,y_n)\in X$ satisfying $(g_1(x_1),...,g_n(x_n))=(g_1(y_1),...,g_n(y_n))$, then due to the $n$ bijections, $x_1=y_1,~...,~x_n=y_n \Rightarrow x=y$(injective). ~Secondly $\forall u=(u_1,...,u_n)\in U$, there always be an $x=(g_1^{-1}(u_1),...,g_n^{-1}(u_n))\in X$ such that $\{g_1,...,g_n\}(x)=u$(surjective).}

\sldef{The ordered set $\{g_1,...,g_n\}$ given by Lemma \ref{lem:OrderSetbeBJ} as a bijection as well as a collection, denoted by $\mathscr{I}\{g_1,...,g_n\}$, ~is called an $n$-dimensional isomorphic frame. It is further expressed as the following notational forms:
\begin{equation} 
\begin{split}
   \mathscr{I}\{g_1,...,g_n\}:X\to U &= [X~\sharp~U]_{g_1,...,g_n}\\
                                     &= [X_1\times ...\times X_n~\sharp~ U_1\times...\times U_n]_{g_1,...,g_n}\\
                                     &= [X_1, ...,X_n~\sharp~ U_1,...,U_n]_{g_1,...,g_n}.
\end{split}
\end{equation}
$X$ is called the base frame of the isomorphic frame(or ``the base'' for short), and $~U$ is called the image frame of the isomorphic frame(or ``the image''). The bijection $g_i(i=1,...,n)$ is called a (the $i$th) dimensional mapping(DM) of $\mathscr{I}\{g_1,...,g_n\}$.
}
\slnot{We also write $\mathscr{I}^{-1}\{g_1,...,g_n\}$ for the inversion of $\mathscr{I}\{g_1,...,g_n\}$.}
\slthm{$\mathscr{I}\{g_1^{-1},...,g_n^{-1}\}=(\mathscr{I}\{g_1,...,g_n\})^{-1}(=\mathscr{I}^{-1}\{g_1,...,g_n\})$.}

(Proof omitted.)

\subsubsection{Embedding in isomorphic frame}
\slnot{Let $\mathscr{I}\{g_1,...,g_n\} = [X_1\times ...\times X_n~\sharp~U_1\times...\times U_n]_{g_1,...,g_n}$ and $D\subseteq X_1\times ...\times X_n$, $E\subseteq U_1\times ...\times U_n$. \\
\indent (i) $D$ is said to be embedded in the base of $\mathscr{I}\{g_1,...,g_n\}$, for which way we use the sign ``$\vee$'' to write $D\vee\mathscr{I}\{g_1,...,g_n\}$, or $D\vee\big(\mathscr{I}\{g_1,...,g_n\}=[X_1\times ...\times X_n~\sharp~U_1\times...\times U_n]_{g_1,...,g_n}\big)$, etc.\\
 \indent (ii) $E$ is said to be embedded in the image of $\mathscr{I}\{g_1,...,g_n\}$, and for this we write $E\vee\mathscr{I}\{g_1^{-1},...,g_n^{-1}\}$.
}

\slthm{$\mathscr{I}\{g_1,...,g_n\}(D)\vee\mathscr{I}\{g_1^{-1},...,g_n^{-1}\}$ if $D\vee\mathscr{I}\{g_1,...,g_n\}$.
}
This is because the image of $D$ is a subset of the image of the isomorphic frame.

\subsubsection{Bonding on isomorphic frame}\label{subsubsec:bondonIsoFrame}
\slnot{Let $\mathscr{I}\{g_1,...,g_n\} = [X_1\times ...\times X_n~\sharp~U_1\times...\times U_n]_{g_1,...,g_n}$.\\
\indent (i) If function $f:D\to M$ of $(n-1)$ variables is such that $D\subseteq X_1\times ...\times X_{n-1}$ ($D\vee\mathscr{I}\{g_1,...,g_{n-1}\}$) and the range $M\subseteq X_n$ ($M\vee\mathscr{I}\{g_n\}$), then $f$ is said to be bonded on the base of $\mathscr{I}\{g_1,...,g_n\}$. For this we use the sign ``$\wedge$'' to write $(f:D\to M) \wedge\mathscr{I}\{g_1,...,g_n\}$, or the alike. \\
\indent (ii) If $D\subseteq U_1\times ...\times U_{n-1}$ and $M\subseteq U_n$, then $f$ is said to be bonded on the image of $\mathscr{I}\{g_1,...,g_n\}$, and for this we write $(f:D\to M) \wedge\mathscr{I}\{g_1^{-1},...,g_n^{-1}\}$, $f \wedge\mathscr{I}^{-1}\{g_1,...,g_n\}$ or the alike.}

\slrem{In the extreme case $n=1$, the above notation applies to a single variable, e.g. $x\in M\subseteq X$, or $x\in M\subseteq U$, which are treated as functions of 0 variables bonded on (the base or the image of) the isomorphic frame.}

\slnot{Let $\mathscr{I}\{g\}= [X~\sharp~U]_g$ be 1 dimensional, and $k$-tuple $\underline{x}\in X$, then $\underline{x}$ is said to be either embedded in or bonded on the base of $\mathscr{I}\{g\}$, written as e.g. $\underline{x}\vee[X~\sharp~U]_g$, or $\underline{x}\wedge[X~\sharp~U]_g$.
}

\subsubsection{About embedding and bonding}
Embedding and bonding are 2 basic styles the objects of MA ``attach or tie'' to isomorphic frames. A set embedded in an isomorphic frame is a part of the latter while a function bonded on an isomorphic frame is not that way strictly. Most definitions in the paper will be based on these 2 concepts.

From now on, all the isomorphic frames are with $\underline{strictly ~monotone~ bijections}$, i.e. strictly monotone (invertible) real functions, as their dimensional mappings unless otherwise specified. These isomorphic frames are denoted as $\mathscr{I}_m\{g_1,...,g_n\}$.

\subsection{Isomorphic number and isomorphic variable}
\sldef{\label{defin:IsoNumVar} With 1 dimensional $\mathscr{I}_m\{g\} = [X~\sharp~U]_g$, $\forall x\in X$, $u=g(x)\in U$ is called the isomorphic number of $x$ generated by mapping $g$(or by $\mathscr{I}_m\{g\}$); In terms of variables, $u$ is called the isomorphic variable of $x$ generated by mapping $g$. $u$ is specially denoted as $u=\varphi(x:g)$ or $u=\varphi(x:\mathscr{I}_m\{g\})$.}

\slrem{Any real number $x$ is an isomorphic number of itself generated by identity mapping ($y=x$).}
An example of isomorphic number in Physics is Conductance G with regards to Resistance R, since $G=1/R$.

With Definition \ref{defin:IsoNumVar} $u$ is already called the ``isomorphic number of $x$'' without an operation defined, because with such $u$ it is ready and easy to introduce 2 genuine operations with ``isomorphism''. For examples:
\begin{itemize}\setlength{\itemsep}{-0.1em}
\item[1).] If there is an operation of arithmetic ``+'' in $U$, that $\forall u_1,u_2\in U$, $\exists u_1+u_2\in U$, then we can define a binary operation denoted by $[\ .+.\ ]_g$ on $X$, such that $\forall x_1,x_2\in X$, $\exists x_3\in X$ which holds $x_3=[x_1+x_2]_g=g^{-1}(g(x_1)+g(x_2))$.
\item[2).] Similarly, if there's an operation on $U$ that computes the mean of $u_1$, $u_2$, there is also a mapped operation on $X$, which computes a ``special mean'': $\bar x$ of $x_1,x_2$. 
\end{itemize}

More generalized, $u$ and $x$ are ``2 conjugate variables embedded in $[X~\sharp~U]_g$'' or ``2 imaging functions of 0 variables bonded on $[X~\sharp~U]_g$''. 

\subsection{Dual-variable-isomorphic function}
With a function of 1 variable bonded on the base of a 2-D isomorphic frame, on the image of the latter bonded is the so called ``dual-variable-isomorphic function''.

\subsubsection{Definition}

\sldef{\label{defin:IsoFun} Let $(f\colon D\to M)\wedge\big(\mathscr{I}_m\{g,h\} = [X,Y~\sharp~U,V]_{g,h}\big)$ and $E=g(D),~N=h(M)$,
\begin{itemize}
\item Function  $(h\circ f\circ g^{-1})\colon E\to N$ is defined to be the dual-variable-isomorphic(DVI) function of $f$ generated by mapping $g,h$(or generated by $\mathscr{I}_m\{g,h\}$). It is denoted by $\varphi(f:g,h)$, or $\varphi(f:\mathscr{I}_m\{g,h\})$,
\begin{equation} 
    \varphi(f:g,h)\colon=(h\circ f\circ g^{-1})\colon E\to N.
\end{equation}
\item $g$ and $h$ are called the independent variable's generator mapping(function) or dimensional mapping(IVDM), and dependent variable's generator mapping(function) or dimensional mapping(PVDM) of $\varphi$ respectively.
\item $E\subseteq U$ is called the isomorphic domain of $D$ generated by mapping $g$; $N\subseteq V$ is called the isomorphic range of $M$ generated by mapping $h$.
\end{itemize}
}

\begin{figure}[htb]
\begin{center}
\includegraphics[viewport=0 0 706 559,scale=0.4]{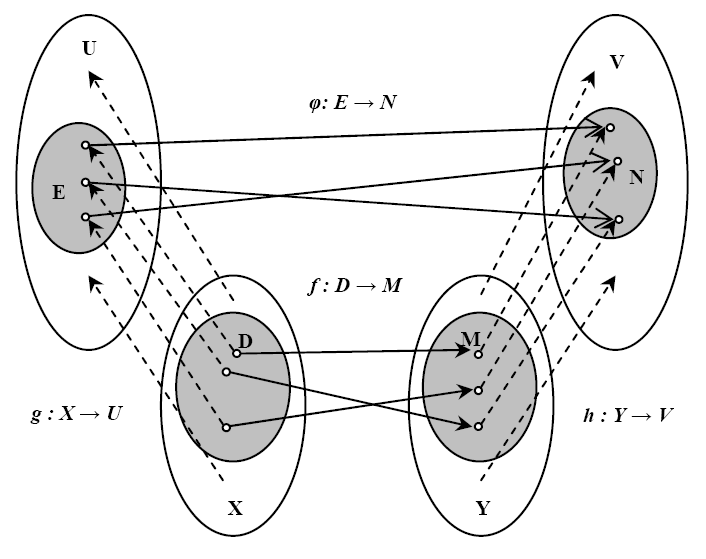}
\end{center}
\begin{illustration}{
\begin{center}
\label{illus:IsoFun}Function $f$ and its DVI function $\varphi\colon E\to N$ \\bonded on a 2-D isomorphic frame
\end{center}
}\end{illustration}
\end{figure}

\slthm{$\varphi(f:g,h) \wedge \big(\mathscr{I}_m\{g^{-1},h^{-1}\} = [U,V~\sharp~X,Y]_{g^{-1},h^{-1}}\big)$.}
That is, $\varphi(f:g,h)$ is bonded on the image of $\mathscr{I}_m\{g,h\}$. Proof omitted.

Illustration \ref{illus:IsoFun} is the visual impression of the definition and the theorem, where $f\wedge\mathscr{I}_m\{g,h\}$ and $\varphi\wedge\mathscr{I}_m\{g^{-1},h^{-1}\}$. A hidden mapping between $f$ and $\varphi$ is formed by $\mathscr{I}_m\{g,h\}$ under which these 2 functions are deemed as either (i) 2 similar elements in different spaces, (ii) 2 sets of relations of same cardinal number, or (iii) 2 unary operations of same structure but each spanning across a pair of number sets.

In traditional terms, the range $M$ here is the image of $f$, and $Y$ is the codomain. The following expression sometimes may also be used for DVI functions, with which the range $N\subseteq V$ is implied, where $V$ is also the codomain.
\begin{equation} 
    \varphi(f:g,h)\colon=(h\circ f\circ g^{-1})\colon E\to V.
\end{equation}

\subsubsection{Sub-classing of DVI function}\label{subsubsec:IsoFunSpecialCase}
\sldef{In view of Definition \ref{defin:IsoFun}, the following special cases can be considered:
\begin{itemize}\setlength{\itemsep}{-0.1em}
  \item[1).] Let $g$ be identity mapping, $\varphi\colon=(h\circ f)\colon D\to N$ is called the dependent-variable-isomorphic(PVI) function of $f$ generated by $h$.
  \item[2).] Let $h$ be identity mapping, $\varphi\colon=( f\circ g^{-1})\colon E\to M$ is called the independent-variable-isomorphic(IVI) function of $f$ generated by $g$.
  \item[3).] Let $Y=X$, $h=g$, $\varphi\colon=(g\circ f\circ g^{-1})\colon E\to N$ is called the same-mapping dual-variable-isomorphic function of $f$ generated by $g$.
  \item[4).] In general case, $h \neq g$, $\varphi\colon=(h\circ f\circ g^{-1})\colon E\to N$ is called the (general) dual-variable-isomorphic function of $f$ generated by $g,h$.
  \item[5).] Let $f(x)=x$, $\varphi\colon=(h\circ g^{-1})\colon E\to N$ is called the dual-variable-isomorphic function of identity generated by $g,h$.
  \item[6).] Let $g$, $h$ both be identity mappings, $\varphi\colon=(h\circ f\circ g^{-1})\colon E\to N$  is equivalent to $f\colon D\to M$. Both domains are $D$, and both ranges are $M$.  i.e. $f$ is a special dual-variable-isomorphic function of itself generated by identities.
  \item[7).] For monotone function $f\colon D\to M(range~M)$, $f^{-1}$ is the dual-variable-isomorphic function of $f$ generated by $f$, $f^{-1}$. 
  \\Above special cases also could be treated as 7 sub-classes of DVI function.
\end{itemize}
}

\subsubsection{Anti-dual-variable-isomorphic function}

\sldef{Let $\varphi\colon=(h\circ f\circ g^{-1})\colon E\to N$  be the DVI function of $f\colon D\to M$. Then $f$ is called the anti-dual-variable-isomorphic function of $\varphi$.}

\slrem{With respect to its DVI function $\varphi$ generated by mapping $g$, $h$, \,$f\colon D\to M$ can be represented by $f=(h^{-1}\circ \varphi \circ g) :D\to M$. Obviously $f$ is the DVI function of $\varphi$ generated by $g^{-1}$, $h^{-1}$. }

This can be observed in Illustration \ref{illus:IsoFun}.

The following are 3 useful special cases of DVI functions.

\subsubsection{V-scaleshift: Vertical scale and shift of a function}
\sldef{For a real $f\colon D\to M$ and 2 constants $k\ne0$, $C$, define the function
\begin{equation}
    V_{ss}\big(f:k,C\big)\colon=~~v=kf(x)+C ~~\big(v\in k(M)+C\big)
\end{equation}
a V-scaleshift of $f$ with scale $k$ and shift $C$. Define the set
\begin{equation}
    \mathbb{V}f = \{ V_{ss}\big(f:k,C\big):k,C\in \mathbb{R},k\ne0\}
\end{equation}
the V-scaleshift space of $f$.
}
\sllem{1).$f\in \mathbb{V}f$; ~2).$\mathbb{V}g=\mathbb{V}f$, if $g\in \mathbb{V}f$.}
\slrem{A V-scaleshift of $y=f(x)$ is a special case of dependent-variable-isomorphic function of $f$, where the PVDM is $v=ky+C$.}
\sllem{Suppose $f\colon D\to M$ be (strictly) convex or (strictly) concave on $D$, then $V_{ss}\big(f:k,C\big)$ and $f$ are of the same (strict) convexity  if $k>0$ or of the opposite (strict) convexity if $k<0$.}

Proof is omitted.

\slnot{In this paper, $\mathbb{V}x$($\mathbb{V}y$) is denoting the V-scaleshift space of identity mapping $g(x)=x$($h(y)=y$).}

\subsubsection{H-scaleshift: Horizontal scale and shift of a function}
\sldef{For a real $f\colon D\to M$ and 2 constants $k\ne0$, $C$, define the function
\begin{equation}
    H_{ss}\big(f:k,C\big)\colon=~~y=f\big(\frac1{k}(u-C)\big) ~~\big(u\in k(D)+C\big)
\end{equation}
an H-scaleshift of $f$ with scale $k$ and shift $C$. Define the set
\begin{equation}
    \mathbb{H}f = \{ H_{ss}\big(f:k,C\big):k,C\in \mathbb{R},k\ne0\}
\end{equation}
the H-scaleshift space of $f$.
}
\sllem{1).$f\in \mathbb{H}f$; ~2).$\mathbb{H}g=\mathbb{H}f$, if $g\in \mathbb{H}f$.}
\slrem{An H-scaleshift of $y=f(x)$ is a special case of independent-variable-isomorphic function of $f$, where the IVDM is $u=kx+C$.}
\sllem{\label{lem:HssSameConvex}Suppose $f\colon D\to M$ be (strictly) convex or (strictly) concave on $D$, then $H_{ss}\big(f:k,C\big)$ is of the same (strict) convexity on $k(D)+C$ as $f$ on $D$.
}
\slprf{Suppose $f$ is convex, then $\forall u_1,u_2\in k(D)+C$, $\forall \lambda\in [0,1]$ $\exists x_1=(u_1-C)/k, x_2=(u_2-C)/k\in D$, such that $f(\lambda x_1+(1-\lambda)x_2)\leq \lambda f(x_1)+(1-\lambda )f(x_2)$. This $\Rightarrow$ $f(\lambda(u_1-C)/k+(1-\lambda)(u_2-C)/k)=f((\lambda u_1+(1-\lambda)u_2)/k-C)\leq \lambda f((u_1-C)/k)+(1-\lambda )f((u_2-C)/k)$, which means $H_{ss}\big(f:k,C\big)$ is also convex. While $f$ has other (strict) convexities, $H_{ss}\big(f:k,C\big)$ will also copy.}

\sllem{\label{lem:HssVssMirror}Let $g,h$ be both invertible. 1).If $g\in \mathbb{H}h$, then $g^{-1}\in \mathbb{V}(h^{-1})$. 2).If $g\in \mathbb{V}h$, then $g^{-1}\in \mathbb{H}(h^{-1})$.}
Proof is omitted.

\subsubsection{HV-scaleshift: Horizontal and vertical scale and shift of a function}
\sldef{For a real $y=f(x):D\to M$ and constants $p\ne0$, $Q$, $k\ne0$, $L$,  define the function
\begin{equation}
\begin{split}
    HV_{ss}&\big(f:p,Q;k,L\big)\colon=~~v=kf\big(\frac1{p}(u-Q)\big)+L \\
    &\big(u\in p(D)+Q,~v\in k(M)+L\big)
\end{split}
\end{equation}
an HV-scaleshift of $f$ with scale $p,k$ and shift $Q,L$. Define the set
\begin{equation}
    \mathbb{HV}f = \{ HV_{ss}\big(f:p,Q;k,L\big):p,Q,k,L\in \mathbb{R},p\ne0,k\ne0\}
\end{equation}
the HV-scaleshift space of $f$.
}
\sllem{1).$f\in \mathbb{HV}f$; ~2).$\mathbb{HV}g=\mathbb{HV}f$, if $g\in \mathbb{HV}f$.}
\slrem{
An HV-scaleshift of $y=f(x)$ is a special case of DVI function of $f$:
\begin{equation}
    HV_{ss}\big(f:p,Q;k,L\big)=\varphi(f:g,h),
\end{equation}
where the DMs $g,h$ are defined by $~u=g(x)=px+Q, ~v=h(y)=ky+L$.
}

\subsection{Isomorphic number-axis and pair}
\subsubsection{Isomorphic number-axis}

\sldef{With an isomorphic frame $\mathscr{I}_m\{g\}=[X~\sharp~U]_{g}$, i.e. a strictly monotone bijection $g\colon X\to U(X,\,U\subseteq\mathbb{R})$, do the following modification to a real number-axis:
\begin{enumerate}
  \item[1).] Replace the mark of  $\,\forall u\in U$ with $g^{-1}(u)$; replace previous 0 point of the axis with mark ``$\alpha$''(alpha), or mark ``$\alpha_g$'';
  \item[2).] Remove all the marks for the rest of the points in the axis if still any;
  \item[3).] If g is strictly increasing, then keep the arrow of the axis pointing to the right; otherwise invert the arrow to point to the left;
  \item[4).] Mark ``$:g$''(a colon followed by `$g$', or by `$g(x)$', `$g(x)=...$' etc.) right beneath the arrow.
\end{enumerate}
Then the modified number-axis is called an isomorphic number-axis generated by
mapping $g$, or generated by isomorphic frame $\mathscr{I}_m\{g\}$. $g$ is called the generator mapping, or dimensional mapping(DM) of the isomorphic number-axis. The point set $X=g^{-1}(U)$ on the axis is called the domain of the isomorphic number--axis. $\alpha$ is called the origin of the axis.}

\slrem{For above definition, if $0\in U$, then $\alpha_g=g^{-1}(0)\in X$; otherwise $\alpha_g$ only serves as a special mark without representing any real number.}

\subsubsection{Isomorphic number-axis pair}

\slnot{Isomorphic number-axis pair is a number-axis system consisting of an isomorphic number-axis and a real number-axis. It is established by the following steps:
\begin{enumerate}
  \item[1).] Given an isomorphic frame $\mathscr{I}_m\{g\}=[X~\sharp~U]_{g}$; 
  \item[2).] Establish 2 parallel real number axes in a plane such that:
  \begin{enumerate}
    \item the 2 axes are with same direction;
    \item the distance between 2 axes is any $\varepsilon>$0;
    \item an imaginary line connecting 0 points of the 2 axes is perpendicular to both axes;
  \end{enumerate}
  \item[3).] Modify the lower real number-axis to an isomorphic number-axis generated by $g$.
\end{enumerate}
Then the system of the 2 number axes is called an isomorphic number-axis pair generated
by mapping $g$, or generated by isomorphic frame $\mathscr{I}_m\{g\}$. $g$ is called the generator mapping, or dimensional mapping(DM) of the isomorphic number-axis pair. The upper real number-axis is called auxiliary(aux.) axis of the pair. The lower
isomorphic number-axis is called the main axis of the pair.
}

\begin{figure}[htb]
\begin{center}
\includegraphics[viewport=0 0 934 261,scale=0.4]{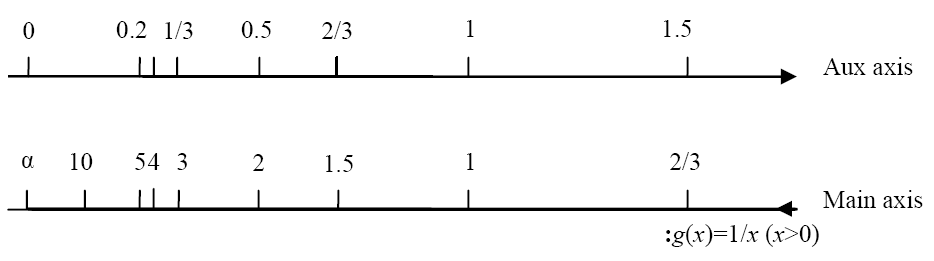}
\end{center}
\begin{illustration}{
\begin{center}
\label{illus:IsoAxPair}An isomorphic number-axis pair
\end{center}
}\end{illustration}
\end{figure}

The Illustration \ref{illus:IsoAxPair} shows an isomorphic number-axis pair generated by mapping $g(x)=1/x (x>0)$. $\alpha$ is not in the domain of main axis.

\subsubsection{Properties of isomorphic number-axis pair}
\paragraph{Correspondence between main axis and aux. axis\\}\label{par:corresMainAux}

According to above, one can conclude easily that:
\begin{enumerate}
  \item[1).] For any $x$ in the domain $X$ of main axis, its orthogonal projection on the aux. axis
is $g(x)$; i.e. on the main axis the actual distance from $\alpha$ to $x$ is the value of $g(x)$,
measured in the unit(also called ``metric'' here) and direction of auxiliary axis.
  \item[2).] For any $u\in U$ on the aux. axis, its orthogonal projection on the main axis $g^{-1}(u)$ is in the domain $X$.
  \item[3).] The $\alpha$ point of main axis always has an orthogonal projection 0 on the aux. axis
and vice versa.
\end{enumerate}

\paragraph{Directed line segment on an isomorphic number-axis\\}
Let $a$ be a point in the domain of an isomorphic number-axis generated by mapping $g$. Based on \ref{par:corresMainAux} case 1), the value of directed line segment $\vec{\alpha a}$ is $\vec{\alpha a}=g(a)-0=g(a)$ if the direction and the unit of aux. axis are taken as standard. As regards any 2 points $a$, $b$ in the domain, the value of directed line segment is
\begin{equation}
  \vec{ab}=\vec{\alpha b} - \vec{\alpha a}=g(b)-g(a).
\end{equation}

\paragraph{\label{sect:FixProp}Fixed proportion division on an isomorphic number-axis\\}
Given 2 points $p_1$, $p_2$ in the domain of an isomorphic number-axis generated by $g$. If an
arbitrary point $p(p\neq p_1, p_2)$ in the domain divides $\vec{p_1p_2}$ into $\vec{p_1p}:\vec{p p_2}=k(k\neq -1,0)$, then
\begin{eqnarray*}
    \frac{g(p)-g(p_1)}{g(p_2)-g(p)}=k
    &\Longrightarrow& g(p)=\frac{g(p_1)+kg(p_2)}{1+k} \\
    &\Longrightarrow& g(p)=\frac{1}{1+k}g(p_1)+\frac{k}{1+k}g(p_2).
\end{eqnarray*}

\noindent
Especially when $p$ is in-between of $p_1$, $p_2$, $k>0$ for $g(p_1)>g(p_2)$ or $g(p_1)<g(p_2)$. Let
\begin{eqnarray*}
    \lambda_1=\frac{\vec{pp_2}}{\vec {p_1p_2}}=\frac{1}{1+k}>0,~
    \lambda_2=\frac{\vec{p_1p}}{\vec {p_1p_2}}=\frac{k}{1+k}>0,
\end{eqnarray*}
then $\lambda_1+\lambda_2=1$, $\lambda_2:\lambda_1=k$, thus
\begin{equation} \label{equ:FixProportion}
    p=g^{-1}\big(\lambda_1g(p_1)+\lambda_2g(p_2)\big)
\end{equation}
with $\lambda_1,\lambda_2\in (0,1)$. If $p$ is the geometrical midpoint of $p_1, p_2$, then
\begin{equation}
    p=g^{-1}\bigl(\frac12g(p_1)+\frac12g(p_2)\bigl).
\end{equation}

\paragraph{Representation of interval on an isomorphic number-axis}

\slthm{\label{thm:IntervalonAxis}Given an isomorphic number-axis pair generated by $\mathscr{I}_m\{g\}=[X~\sharp~U]_{g}$ and an interval $I\subseteq X$, if $g$ is continuous on $I$, then $I$ can be represented by a continuous line segment, line, or ray on the main axis.}


\slprf{As $g$ is continuous on interval $I$, set  $g(I)$ is also an interval, which can be
represented by a continuous line segment, line, or ray on the aux. axis of the axis pair. In
this case, for $\forall a,b(a<b)\in I$ on the main axis, and an arbitrary point $c$ in-between of $a,\,
b$, the orthogonal projection of $c$, $c'$ must be in-between of $g(a), g(b)$. Thus $c'\in g(I)$, then
$c'\in U$, therefore according to \ref{par:corresMainAux} case 2), $c$ is in the domain of main axis, i.e. $c\in X$. As $c$ is in-between of $a, b$ and $g$ is strictly monotone, therefore $a<c<b$. Hence $c\in I$. Finally it concludes with Theorem \ref{thm:IntervalonAxis}.
}

\slrem{If otherwise $g$ is not continuous on $I$, whether $I$ can be represented in such a way is uncertain.}

\subsection{Planar dual-isomorphic rectangular coordinate system pair}
Planar dual-isomorphic rectangular coordinate system pair is basically the 2-dimensional version of isomorphic number-axis pair.

\subsubsection{Establishing of the system pair}

\slnot{\label{not:EstabDualIsoSys}
Planar dual-isomorphic rectangular coordinate system pair is established through the following steps:
\begin{enumerate}
  \item[1).] Given a 2-dimensional $\mathscr{I}_m\{g,h\}=[X,Y~\sharp~U,V]_{g,h}$,  (~i.e. 2 strictly monotone bijection $~g\colon X\to U$, $h\colon Y\to V$ ) and 2 parallel plane $a, b$ which have arbitrary distance of $\varepsilon \geq0$.
  \item[2).] Build a planar rectangular coordinate system $A$ in plane $a$. Build a second planar rectangular coordinate system $B$ in plane $b$ such that $B$ is exactly the orthogonal projection of $A$ on plane $b$.
  \item[3).] Modify $A$'s horizontal axis to an isomorphic number-axis generated by mapping $g$, and modify $A$'s vertical axis to an isomorphic number-axis generated by mapping $h$.
\end{enumerate}
\noindent
Then
\begin{enumerate}
  \item[a).] The modified coordinate system $A$ is called a planar dual-isomorphic rectangular coordinate system generated by mapping $g, h$ or by $\mathscr{I}_m\{g,h\}$, or called for short the dual-isomorphic system, or the isomorphic coordinate system.
  \item[b).] Coordinate system $B$ is called the referential auxiliary coordinate system of $A$, or called for short the referential auxiliary system, or the auxiliary(aux.) system.
  \item[c).] $A$ and $B$ together are called the planar dual-isomorphic rectangular coordinate system pair generated by mapping $g, h$ or by $\mathscr{I}_m\{g,h\}$.
  \item[d).] The horizontal axis of $A$ can be called the $x$ axis of $A$, and vertical axis called $y$ axis of $A$. $g$ is called the horizontal dimensional mapping of $A$, and $h$ is called the vertical dimensional mapping of $A$.
  \item[e).] The intersection point of $x$ axis and $y$ axis is the $\alpha$ point of $x$ or $y$ axis, which is called the origin of $A$, and is represented by $\alpha_{gh}$.
  \item[f).] The horizontal axis and vertical axis of system $B$ can be called the $u$ axis and $v$ axis of $B$ respectively, in order to differentiate $x$ axis and $y$ axis of $A$.\\\\
\noindent
For an arbitrary point $x\in X$ in the domain of $x$ axis of system $A$, and an arbitrary point $y\in Y$ in the domain of $y$ axis, now in plane $a$ draw a line from $x$ perpendicular to $x$ axis and a line from $y$ perpendicular to $y$ axis to cross each other in a unique point $p$ in plane $a$. Here an ordered pair $(x,\,y)$ is used to denote $p$ as $p(x,\,y)$.
  \item[g).] Then point set consisting of all possible $p$, denoted by $\{p(x,\,y)\colon x\in X, y\in Y\}=X\times Y$ is called the domain of dual-isomorphic system $A$.
  \item[h).] For any point $p(x,\,y)$ in the domain, the ordered pair $(x,\,y)$ denoting $p$ is called the coordinate of point $p$, also written as $p=(x,\,y)$.  $x$ is called the horizontal coordinate of $p$, and $y$ the vertical coordinate of $p$. It is stipulated that the coordinate of origin $\alpha_{gh}$ is $(\alpha_{g}, \alpha_{h})$. i.e. $\alpha_{gh} =(\alpha_{g}, \alpha_{h})$.
\end{enumerate}
}

\slrem{\label{rem:coortransform}
According to above notation, it is easy to prove that, for an arbitrary point $p(x, y)$ in the domain of system $A$, its orthogonal projection on the referential aux. system $B$ has the coordinate of $(g(x),h(y))$.
}

\begin{figure}[htb]
\begin{center}
\includegraphics[viewport=0 0 945 640,scale=0.46]{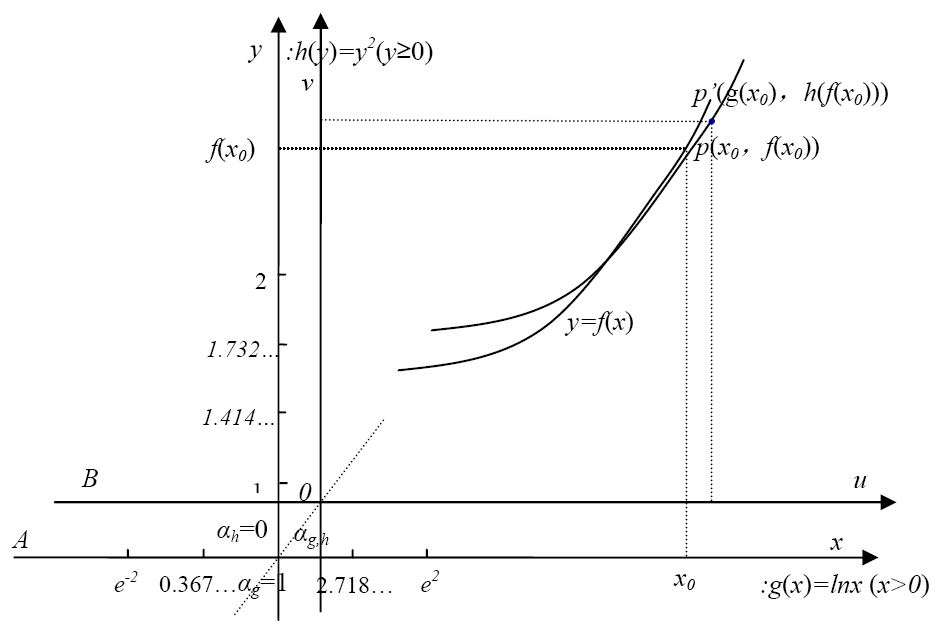}
\end{center}
\begin{illustration}{
\begin{center}
\label{illus:IsoCor}An perspective view of \\ planar dual-isomorphic rectangular coordinate system pair
\end{center}
}\end{illustration}
\end{figure}

A perspective view of the an planar dual-isomorphic rectangular coordinate system pair
generated by mapping  $g(x)=\ln x (x>0)$, $h(y)=y^2(y\geq0)$ is shown in Illustration \ref{illus:IsoCor}. The domain of coordinate system $A$ is the area above $x$ axis joining $x$ axis.

\subsubsection{Special cases of planar dual-isomorphic rectangular coordinate system}

\begin{enumerate}
  \item[(1).] Let  $g, h$  both be identity mapping, the planar dual-isomorphic rectangular
coordinate system is a planar rectangular coordinate system.
  \item[(2).] Let $g$ be identity mapping while $h$ is not, the dual-isomorphic system is called a
vertical-axis-isomorphic coordinate system generated by mapping $h$.
  \item[(3).] Let $h$ be identity mapping while $g$ is not, the dual-isomorphic system is called a
horizontal-axis-isomorphic coordinate system generated by mapping $g$.
  \item[(4).] Let $Y=X, h=g$, the dual-isomorphic system is specifically called a same mapping
dual-isomorphic coordinate system generated by mapping $g$.
\end{enumerate}

\subsubsection{Function's graph on planar dual-isomorphic rectangular coordinate system pair}
Previously in Section \ref{subsubsec:bondonIsoFrame} we have discussed bonding of a function on an isomorphic frame. Here our ready-prepared
isomorphic coordinate system is right for the purpose of visually and geometrically representing of a function in itself in the way that is similar to the way a function is represented in a traditional Cartesian system. That is also important for visually and precisely representing ``how a function is bonded on a 2-D isomorphic frame''.

Let $(f\colon D\to M)\wedge\big(\mathscr{I}_m\{g,h\} = [X,Y~\sharp~U,V]_{g,h}\big)$ and $E=g(D),~N=h(M)$, the graph of $y=f(x)(x\in D)$ is ``plotted'' on a planar dual-isomorphic rectangular coordinate system $A$ generated by the mapping e.g. $g(x)=\ln x(x>0),~ h(y)=y^2(y\geq0)$ as shown in Illustration \ref{illus:IsoCor}.

\slnot{The so-called ``$f$'s graph being plotted on the system $A$'' is just to mark the point set $\{(x,y):y=f(x),~x\in D\}$ in the domain of dual-isomorphic system A, such a geometrical graph can be visually seen on it.}

In addition we consider $f$ has a dual-variable-isomorphic function $\varphi$ generated by $g,h$.

The orthogonal projection of the $f$'s graph in $A$ onto $B$ is considered to be the point set in $B$ constituted by all orthogonal projections of every single point on $f$'s graph in $A$.

According to Remark \ref{rem:coortransform} \,for an arbitrary point $p(x_0,y_0)$ on $f$ its orthogonal projection is $p'(g(x_0),h(f(x_0)))$. Let $u_0=g(x_0)$, then $p'$ can be denoted by $p'(u_0,h[f(g^{-1}(u_0))])$, which in $B$ is right a point on the graph of $\varphi$, denoted by $v=\varphi(u)=h[f(g^{-1}(u))](u\in E=g(D))$, when $u=u_0$. If $x_0$ traverses every value in $D$, then $u_0$ traverses $E$ correspondingly. Therefore the orthogonal projection of graph $f$ onto the referential aux. system is right the graph of $\varphi$ on the referential aux. system. Especially when the distance between system $A, B$, the $\varepsilon=0$, i.e. the two coordinate system is in the same plane, the two graphs are completely overlapped, i.e. two graphs are always congruent.  Therefore the following is true:

\slthm{\label{thm:GraphCorres}Let $(f\colon D\to M)\wedge\big(\mathscr{I}_m\{g,h\} = [X,Y~\sharp~U,V]_{g,h}\big)$ and $E=g(D),~N=h(M)$, then the graph of function $f$ on a planar dual-isomorphic rectangular coordinate system $A$ generated by $\mathscr{I}_m\{g,h\}$ is congruent to the graph of its DVI function $\varphi\colon=(h\circ f\circ g^{-1})\colon E\to N$ generated by $\mathscr{I}_m\{g,h\}$ on the referential auxiliary system $B$.
}

\slrem{ When the distance between coordinate systems $A, B$ is $0$, i.e. plane $a$ and $b$ are overlapped, the two graphs are completely overlapped. For an arbitrary point $(x,f(x))$ on the graph of $f$, the corresponding point on the graph of $\varphi$ is $(u,\varphi(u))=(g(x),h(f(x)))$, where $u=g(x)$.
}

\subsection{Density on isomorphic number-axis and on isomorphic coordinate system}
\subsubsection{Density on an isomorphic number-axis}

\sldef{On the domain $X$ of isomorphic number-axis generated by $\mathscr{I}_m\{g\}=[X~\sharp~U]_g$, define the following function
\begin{equation}
    (x,w)\mapsto \frac{w-x}{\vec{xw}}:=(x,w)\mapsto \frac{w-x}{g(w)-g(x)}
\end{equation}
to be the average density of directed line segment $\vec{xw}$ on the axis, where $x,w\in X, x\neq w$. If there exists the limit
\begin{equation}
    \lim_{w\to x}\frac{w-x}{g(w)-g(x)}~,
\end{equation}
it is called the density at $x$ on the isomorphic number-axis generated by $g$, and it is denoted by $\rho _g(x)$.}

The density of directed line segment $\vec{xw}$ on an isomorphic number-axis reflects the ratio of ``local metric difference of main axis ($w-x$)'' over ``corresponding local metric difference of aux. axis ($g(w)-g(x$))'' with the direction taken into account, whereby it is reasonable to say it actually reflects ``how densely the points are averagely distributed in the line segment there on the main axis''. That implies generally the points are not evenly distributed on an isomorphic number-axis, as opposed we consider points are evenly distributed on its aux. axis, i.e. a real number-axis, for granted.

 As for density at a point $x$ on the isomorphic number-axis, since it is the limit of the former, it reflects how dense the local point $x$ is, comparable in the notion to the case in a real world where we consider how dense the mass of certain material at a specific point of space is, i.e. the density of mass.

\slrem{\label{rem:densityatx}Obviously if and only if $g$ is differentiable at point $x$, and $g'(x)\neq0$, the density at $x$ exists, and it is
\begin{equation}\label{equ:DensiFun}
    \rho_g(x)=\frac{1}{g'(x)}~.
\end{equation}
$\rho_g(x)$ is positive if $g$ is strictly increasing, and negative if $g$ strictly decreasing.
}

\slnot{The reciprocal of $g'$ is called the density function of isomorphic number-axis generated by $g$, denoted by $\rho_g=1/g'$.
}

For the $x$ axis of system $A$ in Illustration \ref{illus:IsoCor} where $g(x)=\ln x(x>0)$, the density at $x$ is $x$, and it is always positive. For $g$ being identity, the density at any point on the axis is always $1$, where the axis is identical to a real number-axis.

\subsubsection{Point of singular-density on isomorphic number-axis}
\slrem{With regards to Remark \ref{rem:densityatx}, if $g$ is differentiable at a point $x$ (in the domain of an isomorphic number-axis generated by $g$) and $g'(x)=0$, the density at $x$ must not exist. In this case, we call such point $x$ is a point of singular-density on the isomorphic number-axis.}

At a point $x$ of singular-density, we have
\begin{equation}
    \lim_{w\to x}\frac{g(w)-g(x)}{w-x}=0.
\end{equation}
It means the limit of the ratio of ``corresponding local metric difference of aux. axis ($g(w)-g(x$))'' over ``local metric difference of main axis ($w-x$)'' is 0, therefore it naturally implies the so-called density at $x$ is being $\infty$. It is reminiscent of a ``gravitational singularity'', i.e. a ``black hole'' in the real world, where its mass density is supposedly being $\infty$.

For example, given an isomorphic number-axis generated by $g(x)=x^3$, then $x=0$ (which is also $\alpha$) is a point of singular-density, where $g'(x)=3x^2=0$.

\subsubsection{Density on a dual-isomorphic rectangular coordinate system}
\slnot{ With arbitrary 2 points $p(x,y), ~q(w,z)(x\neq w, ~y\neq z)$ in the domain of a planar dual-isomorphic coordinate system generated by $[X,Y~\sharp~U,V]_{g,h}$, draw 2 lines from $p$ to be perpendicular to $x$ axis and be perpendicular to $y$ axis respectively, and do the same for $q$. The 4 lines enclose an unique rectangular $r$, use the following to reflect the average density of $r$, denoted by $\rho_{gh}(r)$:
\begin{eqnarray}
    \rho_{gh}(r) &=& \frac{|(w-x)(z-y)|}{|S_r|}=\frac{|(w-x)(z-y)|}{|\vec{xw}\cdot\vec{yz}|}\nonumber\\
    &=& \frac{|(w-x)(z-y)|}{|\bigl(g(w)-g(x)\bigl)\bigl(h(z)-h(y)\bigl)|},
\end{eqnarray}
where $|S_r|$ is the area of rectangular $r$ measured in the metric of referential aux. system.
}

\sldef{With arbitrary 2 points $p(x,y), ~q(w,z)(x\neq w, ~y\neq z)$ in the domain of a planar dual-isomorphic coordinate system generated by $[X,Y~\sharp~U,V]_{g,h}$, if there exists the limit
\begin{equation}
    \lim_{w\to x,z\to y}\frac{|(w-x)(z-y)|}{|\bigl(g(w)-g(x)\bigl)\bigl(h(z)-h(y)\bigl)|}.
\end{equation}
it is called the density at point $p(x,y)$ on the coordinate system, denoted by $\rho_{gh}(p)$.
}

\slrem{\label{rem:densityatxy} Obviously if and only if both $g$ and $h$ are differentiable at $x, y$ respectively, and $g'(x)\neq 0, h'(y)\neq 0$, the density at $p$ exists, and it is
\begin{equation}
    \rho_{gh}(p)=\frac{1}{|g'(x)h'(y)|}=|\rho_{g}(x)\rho_{h}(y)|.
\end{equation}
}

\subsubsection{Point of singular-density on isomorphic coordinate system}
\slrem{With regards to Remark \ref{rem:densityatxy}, if both $g$ and $h$ are differentiable at $x, y$ respectively, whilst $g'(x)=0$ and/or $h'(y)=0$, the density at $p$ does not exist. In these case we call such $p$ is a point of singular-density on the isomorphic coordinate system.
}

\section{Isomorphic arithmetic operations}
If $u$ is the isomorphic number of $x(x\in X, X\subseteq\mathbb{R})$ generated by mapping $g$, and there are
arithmetical operations + or - (addition or subtraction) etc. on $g(X)$, then from the perspective of $X$, there are special operations which hereinafter are referred to as ``isomorphic arithmetic operation''. These special operations on $X$ can be used to describe the phenomenon of composition and decomposition of physical quantity in the real world.

\subsection{Isomorphic addition}\label{subsec:IsoAdd}

\sldef{Given a $\mathscr{I}_m\{g\}=[X~\sharp~U]_g$, and arithmetic addition $+$ on $U$.  For arbitrary $a,b\in X$, define $g^{-1}(g(a)+g(b))$ to be an isomorphic addition operation of $a, b$ on $X$ generated by mapping $g$, denoted by $\bigl[a+b\bigl]_g$.
\begin{equation}
    \bigl[a+b\bigl]_g=g^{-1}\bigl(g(a)+g(b)\bigl).
\end{equation}}

\slrem{ In view of above definition, $<X,\,\bigl[\,+\,\bigl]_g>$ and $<U,\,+ >$ are two algebraic systems. As the bijection $g\colon  X\to U$ satisfies that for arbitrary $a, b\in X$, $g(\bigl[a+b\bigl]_g)=g(a)+g(b)$, therefore $g\colon  X\to U$ is an isomorphism from $<X,\,\bigl[\,+\,\bigl]_g>$ to $<U,\,+ >$. i.e. $<X,\,\bigl[\,+\,\bigl]_g>\,\,\,\cong\,\, <U,\,+ >$.
}

\slnot{ We have the following formula for the isomorphic addition of $n$-tuple($n\geq2$) of addend $\underline a$,
\begin{equation}
    \bigl[a_1+ \ldots + a_n\bigl]_g=g^{-1}\bigl(g(a_1)+\ldots +g(a_n)\bigl).
\end{equation}}

There are quite a lot of instances in the real world, which can be described by isomorphic
addition. Here are 3 examples.

\begin{enumerate}
  \item[(1).] Resistors in parallel. When 2 resistors $R_1, R_2$ are connected in parallel in a circuit, the conductance of 2 resistors $G_1, G_2$ which are the isomorphic numbers of $R_1, R_2$ respectively generated by mapping $G=1/R$ (in physics term, conductance is the reciprocal of resistance), will add up. Therefore the equivalent resistance  $R_3 = \bigl[R_1+R_2\bigl]_{G=1/R} = 1/\big(1/R_1+1/R_2\big)$.
  \item[(2).] Composition of noise. When 2 noises $a_1, a_2$ (their values in dB) are acting on one point (e.g. an ear), their $dB$ values won't simply add up to that of composite noise $a_3$. Here the isomorphic numbers of $dB$ values that add up is the ``ratio of sound power'' in physics term, denoted by $r=10^{0.1a}$ , hence the composite noise $a_3=\bigl[a_1+a_2\bigl]_{r=10^{0.1a}} = 10\lg(10^{0.1a_1}+10^{0.1a_2})$.
  \item[(3).] Multiplication on $\mathbb{R}^+$ is a special case of isomorphic addition, generated by $g(x)=\lg x (x>0)$.
\end{enumerate}

An isomorphic addition generated by identity is the normal addition. By illustrating the addition operation of isomorphic number on an aux. axis of an isomorphic number-axis pair, the corresponding isomorphic addition can be illustrated on the main axis. The elaborations are not covered here. See \ref{subsec:IsoOperExample} for a simple example.

\subsection{Isomorphic subtraction}\label{subsec:IsoSubtra}
Isomorphic subtraction is the inverse operation of isomorphic addition. In practice, it describes the phenomena of physical amounts decomposing, e.g. removing a noise from the composite noise. The following defined is a``closed'' isomorphic subtraction.
\sldef{Given a $\mathscr{I}_m\{g\}=[X~\sharp~U]_g$. $0\in U$ and for arbitrary $u\in U, -u\in U$. There is a closed subtraction ``$-$'' on U. For arbitrary $a, b\in X$, define $g^{-1}(g(a)-g(b))$ to be isomorphic subtraction of $a, b$ on $X$ generated by mapping $g$, denoted by $\bigl[a-b\bigl]_g$:
\begin{equation}
    \bigl[a-b\bigl]_g=g^{-1}\bigl(g(a)-g(b)\bigl).
\end{equation}
}
\slnot{For $n$-tuple($n\geq2$) of subtracter $\underline b$, we have the following formula for the isomorphic subtraction:
\begin{equation}
    \bigl[a-b_1-\ldots -b_n\bigl]_g=g^{-1}\bigl(g(a)-g(b_1)-\ldots -g(b_n)\bigl).
\end{equation}
}

Division on $\mathbb{R}^+$ is a special case of isomorphic subtraction, where $g(x)=\ln x (x>0)$ is the bijection from $\mathbb{R}^+$ onto $\mathbb{R}$.

\subsection{Isomorphic multiplication}\label{subsec:IsoMultip}
\sldef{Given a $\mathscr{I}_m\{g\}=[X~\sharp~U]_g$, $T\subseteq\mathbb{R}$ and a multiplication operation ``$\times$'' $:U\times T\to U$. Then for arbitrary $a\in X, t\in T$, define $g^{-1}\bigl(g(a)\times t\bigl)$ to be the isomorphic multiplication of $a$ and $t$ generated by mapping $g$, denoted by $\bigl[a\times t\bigl]_g$:
\begin{equation}
    \bigl[a\times t\bigl]_g=g^{-1}\bigl(g(a)\times t\bigl).
\end{equation}
}

The power operation is a special case of isomorphic multiplication, generated by $g(x)=\ln x:\mathbb{R}^+\to\mathbb{R}$, whilst $T=\mathbb{R}$.

\subsection{Isomorphic division}
Isomorphic divisions are the inverse operations of isomorphic multiplication.
\subsubsection{Isomorphic division type I}\label{subsubsec:IsoDivT1}
\sldef{Given a $\mathscr{I}_m\{g\}=[X~\sharp~U]_g$, $T\subseteq\mathbb{R}, ~0\notin T$ and a division operation ``$\div$'' : $U\div T\to U$. Then for arbitrary $a\in X$, $t\in T$, define $g^{-1}\bigl(g(a)\div t\bigl)$ to be the isomorphic division type I of $a$ and $t$ generated by mapping $g$, denoted by $\bigl[a\div t\big]_g$:
\begin{equation}
    \bigl[a\div t\bigl]_g=g^{-1}\bigl(g(a)\div t\bigl).
\end{equation}
}

The extraction operation is a special case of isomorphic division type I, generated by $g(x)=\ln x:\mathbb{R}^+\to\mathbb{R}$, whilst $T=\mathbb{R}$(excluding 0).

\subsubsection{Isomorphic division type II}\label{subsubsec:IsoDivT2}
\sldef{Given a $\mathscr{I}_m\{g\}=[X~\sharp~U]_g$, $T\subseteq\mathbb{R}$ and a division operation ``$\div$'' : $U\div U\to T$(the divisor does not equal 0). Then for arbitrary $a, b\in X(g(b)\neq 0)$, define $g(a)\div g(b)$ to be the isomorphic division type II of $a$ and $b$ generated by mapping $g$, denoted by $\bigl[a\div b\bigl]_g^{II}$:
\begin{equation}
    \bigl[a\div b\bigl]_g^{II}=g(a)\div g(b).
\end{equation}
}
The logarithm operation is a special case of isomorphic division type II, generated by $g(x)=\ln x:\mathbb{R}^+\to\mathbb{R}$, whilst $T=\mathbb{R}$.

\section{Dual-isomorphic derivative}
\subsection{Definition of dual-isomorphic derivative}

\sldef{ Given $(f\colon D\to M)\wedge\big(\mathscr{I}_m\{g,h\} = [X,Y~\sharp~U,V]_{g,h}\big)$, if in a certain neighborhood of $x\in D$ composite function $h\circ f$ and $g$ are both derivable with respect to $x$, and $g'(x)\neq 0$, then the ratio $(h\circ f(x))'/g'(x)$ is called the dual-isomorphic derivative of $f$ generated by mapping $g, h$ at $x$. denoted by $\bigl[f'(x)\bigl]_{g,h}$:
\begin{equation}\label{equ:DDerVal1}
    \bigl[f'(x)\bigl]_{g,h}=\frac{\big(h\circ f(x)\big)'}{g'(x)}.
\end{equation}
$f$ is said to be $(g,h)$-dual-isomorphic derivable at point  $x$.
}

\sldef{
If $f$ is $(g,h)$-dual-isomorphic derivable everywhere on interval $(a,b)\subseteq D$, $f$ is said to be $(g,h)$-dual-isomorphic derivable on interval $(a,b)$. Meanwhile on interval $(a,b)$ the function $(h\circ f)'/g'$ is called the dual-isomorphic derivative function of $f$ generated by mapping $g, h$ on interval $(a,b)$, denoted by $\bigl[f'\bigl]_{g,h}$:
\begin{equation}\label{equ:DDerFun1}
    \bigl[f'\bigl]_{g,h}=\frac{(h\circ f)'}{g'}.
\end{equation}
}

\slnot{As with respect to variable $x\in(a,b)$, let $y=f(x)$, let $h'$ denote derivative with respect to variable $y$, then
\begin{eqnarray}
    \label{equ:DDerFun2}
    \bigl[f'\bigl]_{g,h}&=&\frac{(h\circ f)'\cdot \Delta x}{g'\cdot \Delta x}
        =\frac{\mathrm{d}h\circ f}{\mathrm{d}g},\\
    \label{equ:DDerFun3}
    \bigl[f'\bigl]_{g,h}&=&\frac{\mathrm{d}f(x)}{\mathrm{d}x}\cdot \frac{\frac{\mathrm{d}h(y)}{\mathrm{d}y}}{\frac{\mathrm{d}g(x)}{\mathrm{d}x}} = f'\cdot \frac{h'}{g'}.
\end{eqnarray}
In this paper the notation $\frac{\mathrm{d}h\circ f}{\mathrm{d}g}$ and $f'\cdot \frac{h'}{g'}$ are also used to denote the dual-isomorphic derivative function of $f$ with respect to $x$.
}

\slnot{If at $x$ $h\circ f$ and $g$ are both right-derivable (or left-derivable) and $g'(x^+)\neq0$ (or $g'(x^-)\neq0$) then
\begin{enumerate}
  \item[1).] the ratio of right derivative (or left derivative) of $h\circ f$ at $x$ over the right derivative (or left derivative) of $g$ at $x$  is called the dual-isomorphic right-derivative (or left-derivative) of $f$ generated by $g$, $h$ at $x$.
  \item[2).] $f$ is said to be $(g,h)$-dual-isomorphic right-derivable (or left-derivable) at $x$.
  \item[3).] For $[a,b]\subseteq D$, if $f$ is $(g,h)$-dual-isomorphic derivable on $(a,b)$, and at $a$ $f$ is $(g,h)$-dual-isomorphic right-derivable, and at $b$ $f$ is $(g,h)$-dual-isomorphic left-derivable, then $f$ is said to be $(g,h)$-dual-isomorphic derivable on $[a,b]$.
\end{enumerate}
}

\subsection{Geometrical meaning of dual-isomorphic derivative}

\begin{figure}[htb]
\begin{center}
\includegraphics[viewport=0 0 693 616,scale=0.4]{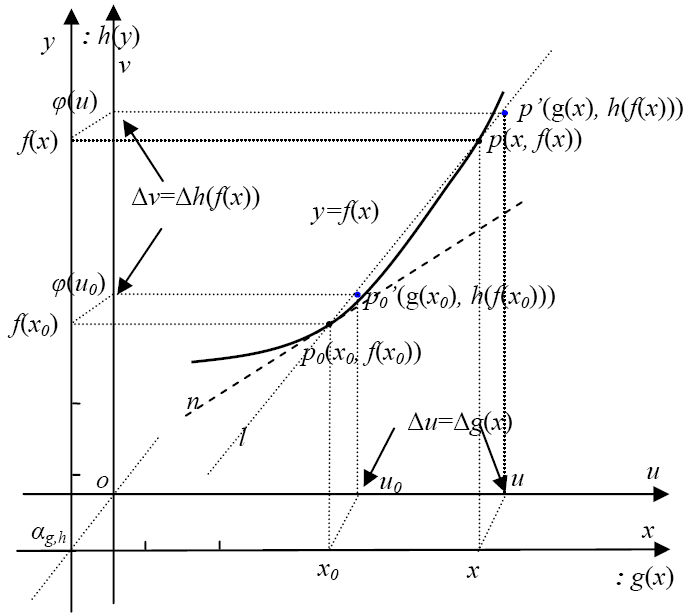}
\end{center}
\begin{illustration}{
\begin{center}
\label{illus:GeoDDer}Geometrical meaning of dual-isomorphic derivative
\end{center}
}\end{illustration}
\end{figure}

In Illustration \ref{illus:GeoDDer} the graph of function $y=f(x)$ is plotted on a planar dual-isomorphic rectangular coordinate system generated by $\mathscr{I}_m\{g,h\}= [X,Y~\sharp~U,V]_{g,h}$, in order to examine the geometrical meaning of dual-isomorphic derivative. Here also consider the dual-variable-isomorphic function of $f$ generated by mapping $g,h$, $v=\varphi (u)=h\bigl[f(g^{-1}(u))\bigl]$, where $u=g(x)$ is the isomorphic number of $x$.

When $f(x)$ is changing from $p_0(x_0,f(x_0))$ to $p(x,f(x))$ in a certain neighborhood of $x_0$, on the referential aux. system the orthogonal projection of $p_0$, i.e. the corresponding point on its dual-variable-isomorphic function: $p_0'(u_0,\varphi (u_0))=p_0'\bigl(g(x_0),h(f(x_0))\bigl)$ (see Remark \ref{rem:coortransform}) is changing accordingly to the orthogonal projection of $p,  p'\bigl(g(x),h(f(x))\bigl)$. Let $\Delta x=x-x_0$, for $v=\varphi(u)$, the increment $\Delta u=g(x)-g(x_0)=\Delta g(x)$, and $\Delta v=h(f(x))-h(f(x_0))=\Delta h(f(x))$ with respect to $\Delta x$. Therefore the ratio $\Delta v/\Delta u=\Delta h(f(x))/\Delta g(x)$ represents the slope of an imaginary line connecting $p_0'$ and $p'$. Employing Theorem \ref{thm:GraphCorres}, these 2 graphs are congruent, the slope of line $p_0'p'$ equals the slope of line $p_0p$(the line $l$). While $u\to u_0$, we see $x\to x_0$ and the slope of line $l$ is changing gradually to the slope of the tangent $n$ of graph $f$ passing point $p_0$ on the dual-isomorphic system. So the slope value is
\begin{eqnarray}
    \lim_{x\to x_0}\frac{\Delta h\big(f(x)\big)}{\Delta g(x)}&=&
    \lim_{x\to x_0}\frac{\Delta h\big(f(x)\big)/\Delta x}{\Delta g(x)/\Delta x} \nonumber \\
    &=&\frac{\mathrm{d}h\big(f(x)\big)/\mathrm{d}x}{\mathrm{d}g(x)/\mathrm{d}x}\bigl|_{x=x_0}=\bigl[f'(x_0)\bigl]_{g,h}.
\end{eqnarray}
which is just the dual-isomorphic derivative of $f$ generated by $g, h$ at $x_0$. On the other hand with $x=g^{-1}(u)$ we have
\begin{equation}
    \frac{\mathrm{d}\varphi(u)}{\mathrm{d}u}
    =\frac{\mathrm{d}h\bigl[f(g^{-1}(u))\bigl]}{\mathrm{d}g^{-1}(u)}
        \cdot\frac{\mathrm{d}g^{-1}(u)}{\mathrm{d}u}
    =\frac{\big(h\circ f(x)\big)'}{g'(x)}.
\end{equation}
Thus :

\slrem{The geometrical meaning of dual-isomorphic derivative of $f$ generated by mapping $g,h$ at $x_0$ can be described as
 \begin{itemize}
   \item the slope of the graph of $f$ at point  $x_0$ on the dual-isomorphic system generated by mapping $g,h$, but measured in the metric of the aux. system.
   \item And the slope of the graph of $\varphi$ (which is the dual-variable-isomorphic function of $f$ generated by mapping $g,h$) on the referential aux. system at point $u_0$ (which is the isomorphic number of $x_0$ generated by mapping g).
 \end{itemize}
}

\subsection{Dual-isomorphic derivative and densities on isomorphic number axes}
For $x\in (a,b)\subseteq D$ and $y=f(x)$, with (\ref{equ:DensiFun}) and (\ref{equ:DDerFun3}) immediately we have:
\begin{equation}
\bigl[f'(x)\bigl]_{g,h}=f'(x)\cdot\frac{h'(y)}{g'(x)}=f'(x)\cdot\frac{\rho_g(x)}{\rho_h(y)}.
\end{equation}

With consideration of foregoing ``geometrical meaning'', this depicts an imaginary picture of how the horizontal and vertical isomorphic axes are ``torturing'' the graph of $f$ with their ``forces of densities''. With identity generator mappings, the dual-isomorphic derivative equals the derivative, since densities on the axes are always 1.

\subsection{A special case of dual-isomorphic derivative - the elasticity of function}\label{subsec:Elasticity}
In Economics, elasticity \cite{TTElas} of a positive differentiable function $f$ at point $x$ is defined to be
\begin{equation}\label{equ:Elasticity}
    Ef(x)=\frac{\frac{\mathrm{d}y}{\mathrm{d}x}}{\frac{y}{x}}=x\frac{f'(x)}{f(x)},
\end{equation}
which evaluates the percentage change of one economic variable with respect to percentage change in another variable.

Because the dual-isomorphic derivative function of a positive function  $f$ generated by mapping $x\mapsto\ln x: \mathbb{R}^+\to\mathbb{R},~y\mapsto\ln y:\mathbb{R}^+\to\mathbb{R}$ is $\mathrm{d}\ln f(x)/\mathrm{d}\ln x=xf'(x)/f(x)$, which is right the elasticity of function $f$. Thus the elasticity of function is a special case of dual-isomorphic derivative.

\slrem{Based on earlier section, the slope of a point on the graph of a positive function $f\colon \mathbb{R}^+\to\mathbb{R}^+$  on the dual-isomorphic system generated by $x\mapsto\ln x: \mathbb{R}^+\to\mathbb{R}, ~y\mapsto\ln y:\mathbb{R}^+\to\mathbb{R}$ is right the elasticity of $f$ at that point.
}

This can be referred to as the geometrical meaning of elasticity of function.

\subsection{Metrical dual-isomorphic derivative}
\subsubsection{Definition}
\sldef{Let $(f\colon D\to M)\wedge\big(\mathscr{I}_m\{g,h\} = [X,Y~\sharp~U,V]_{g,h}\big)$, if at $x\in D$ there is dual-isomorphic derivative $\bigl[f'(x)]_{g,h}\in V$ then
\begin{equation}\label{equ:GenDDer}
    h^{-1}\biggl(\frac{(h\circ f(x))'}{g'(x)}\biggl) ~or~ h^{-1}\biggl(\frac{\mathrm{d}h\circ f(x)}{\mathrm{d}g(x)}\biggl)
\end{equation}
is defined to be the metrical dual-isomorphic derivative of $f$ at point $x$.
}

Metrical dual-isomorphic derivative is a value in $Y$ mapped back by $h^{-1}$ from dual-isomorphic derivative in $V$. As the dual-isomorphic derivative is the slope of the graph of $f$ measured in the metric of the aux. system, after mapped back to $Y$, the metrical dual-isomorphic derivative has an analogous ``metric in the dual-isomorphic system'', which is physically meaningful.

\subsubsection{Special cases of metrical dual-isomorphic derivative}\label{subsubsec:SpecCaseGenDDer}

When identity mappings considered, the metrical dual-isomorphic derivative simplifies to derivative.

When $g(x)=x, ~h(y)=\ln y (y>0)$,
\begin{equation}\label{equ:ExpoDer}
    h^{-1}\biggl(\frac{\mathrm{d}h\circ f(x)}{\mathrm{d}g(x)}\biggl)=\exp \biggl(\frac{f'(x)}{f(x)}\biggl),
\end{equation}
which is the ``exponential derivative'' of function $f$. \cite{GM01}

When $g(x)=\ln x(x>0), ~h(y)=\ln y(y>0)$,
\begin{equation}\label{equ:BigeoDer}
    h^{-1}\biggl(\frac{\mathrm{d}h\circ f(x)}{\mathrm{d}g(x)}\biggl)=\exp \biggl(\frac{xf'(x)}{f(x)}\biggl),
\end{equation}
which is the ``bigeometric derivative'' of function $f$. \cite{GM02}

\section{Isomorphic integral}
\subsection{Isomorphic integral type I}\label{subsec:IsoIntegT1}
\sldef{Given $(f\colon [a, b]\to M)\wedge\big(\mathscr{I}_m\{g,h\} = [X,Y~\sharp~U,V]_{g,h}\big)$ where $g(x)=x$ is identity mapping. If function $h\circ f$ is integrable on [a, b] and the integral $I' \in V$, then
\begin{equation}\label{equ:IsoIntT1}
    I=h^{-1}(I')=h^{-1}\biggl[\int_a^b h(f(x))\mathrm{d}x\biggl]
\end{equation}
is called the isomorphic integral type I of $f$ on $[a, b]$ generated by mapping $h$.
}

Note here $h\circ f$ is considered the dependent-variable-isomorphic function of $f$ generated by $h$ (see \ref{subsubsec:IsoFunSpecialCase} case(1)).

When $h(y)=\ln y(y>0)$, the isomorphic integral type I is
\begin{equation}\label{equ:ProductTypeI}
    I=\exp \bigl(\int_a^b\ln f(x)\mathrm{d}x \bigl),
\end{equation}
Which is the ``geometric integral'' of $f(x)$. \cite{GMKR}

\subsection{Isomorphic integral type II}\label{subsec:IsoIntegT2}
\sldef{Given $(f\colon [a, b]\to M)\wedge\big(\mathscr{I}_m\{g,h\} = [X,Y~\sharp~U,V]_{g,h}\big)$ where $h(y)=y$ is identity mapping and $g$ is continuous and derivable on $[a, b]$. If $f \circ g^{-1}$ is integrable on $[g(a),g(b)]$ and the integral $I'\in U$, then
\begin{equation}\label{equ:IsoIntT2}
    I=g^{-1}(I')=g^{-1}\biggl[\int_{g(a)}^{g(b)}f(g^{-1}(u))\mathrm{d}u\biggl]=g^{-1}\biggl[\int_a^bf(x)g'(x)\mathrm{d}x\biggl]
\end{equation}
is defined to be the isomorphic integral type II of $f$ on [a, b] generated by mapping $g$.
}

\slrem{In view of above definition, in case that $g(a)>g(b)$, $f\circ g^{-1}$ is considered integrable on $[g(a), g(b)]$ if the integral of $f\circ g^{-1}$ on $[g(b), g(a)]$ exists, thus $I'$ takes its negative value.
}

Note here $f\circ g^{-1}$ is considered the independent-variable-isomorphic function of $f$ generated by $g$ (see \ref{subsubsec:IsoFunSpecialCase} case(2)).

\slnot{In this paper, the following also denotes the isomorphic integral type II of $f$
on $[a, b]$
\begin{equation}\label{equ:IsoIntT2variant}
    I=g^{-1}\biggl[\int_a^bf(x)\mathrm{d}g(x)\biggl].
\end{equation}
}

When $g(x)=\ln x(x>0)$, the isomorphic integral type II is
\begin{equation}\label{equ:ElastInt}
    I=\exp\biggl[\int_a^bf(x)\mathrm{d}\ln x\biggl]
    =\exp\biggl[\int_a^b\frac{f(x)}{x}\mathrm{d}x\biggl].
\end{equation}

\slnot{\label{not:ElasticInteg}In this paper, value computed by (\ref{equ:ElastInt}) is called the elastic integral of $f$ on $[a, b]$, as it is the``inversely-related operation'' of the elasticity of function.
}

If we build a function using (\ref{equ:ElastInt}) with its upper integral limit $b$ as a variable $x$, the elasticity of the function with respect to $x$ is $f$.

\slrem{For $f(x)>0$, the elastic integral of the elasticity $Ef(x)$ on $[a, b]$ is the multiplying of $f(x)$ from $a$ to $b$.
}

It is similar with that, the integral of the derivative $f'(x)$ on $[a, b]$ is the increment of $f(x)$ from $a$ to $b$. The proof is easy and omitted.

For instance, if $F(x)=k(x-c)$ is defined on $[a, b](b>a>c,~k>0)$, the multiplying of $F(x)$ from $a$ to $b$ is $F(b)/F(a)=(b-c)/(a-c)$. It's also known the elasticity of $F(x)$ on $[a, b]$ is $EF(x)=xF'(x)/F(x)=x/(x-c)$. Then
\begin{equation}\label{equ:ElastIntExample}
    \exp\biggl[\int_a^b\frac{x}{x-c}\cdot\frac{\mathrm{d}x}{x}\biggl]=\frac{b-c}{a-c}=\frac{F(b)}{F(a)}.
\end{equation}

\section{Isomorphic weighted mean of numbers}\label{sec:IsoWghtMean}
The quasi-arithmetic mean \cite{BULLENPS} or generalized $f$-mean is a generalization of simple means such as the arithmetic mean and the geometric mean, using a function $f$. In this paper it is re-defined as the isomorphic weighted mean with the concept of isomorphic number involved.

\sldef{\label{defin:IsoWeightMean}Given $\mathscr{I}_m\{g\}=[X~\sharp~U]_g$, $n$-tuple($n\geq2$) $\underline{x}\wedge[X~\sharp~U]_g$ and $\underline{u}$ being their respective isomorphic numbers. With positive $\underline{p}$ satisfying $\sum_{i=1}^np_i=1$ if  $\sum_{i=1}^np_iu_i\in U$, then $g^{-1}\bigl(\sum_{i=1}^np_iu_i\bigl)$ is called the isomorphic weighted mean of $n$-tuple (numbers) $\underline x$ generated by $g$(or by $\mathscr{I}_m\{g\}$). Here it is denoted by $\overline{x,p_R}|_g$ (or $\overline{x_{\{i\}},p_R}|_g$, $\overline{x_{i},p_R}|_g$ or simpler $\overline{x,p}|_g$),
\begin{equation}\label{equ:IsoWgtMean}
    \overline{x,p_R}|_g=g^{-1}\bigl(\sum_{i=1}^np_i g(x_i)\bigl).
\end{equation}
$g$ is called the generator mapping, or dimensional mapping of the isomorphic weighted mean.
}

\slrem{The subscript $R$ of $p$ indicates $\underline p$ are the ``relative'' weights (fractions always add up to 1), as comparing to another type of weights known as Frequency Numbers, for which case the definition and formula will be slightly different as in \cite{BULLENPS}. For simplicity, in this paper we use $\overline{x,p}|_g$ which agrees series $\underline p$ are relative (fractional) weights.
}

\slrem{In simple words, isomorphic (weighted) mean is the inverse image of the (weighted) mean of $n$ number of $\varphi(x_i:g)$.}

\subsection{Isomorphic mean of numbers}
\sldef{With Definition \ref{defin:IsoWeightMean},~especially when $p_1=\ldots =p_n=1/n$, the isomorphic weighted mean is called the isomorphic mean of $n$-tuple $\underline x$ generated by $g$, it is denoted by $\overline{x_{\{i\}}}|_g$, $\overline{x_{i}}|_g$ or simpler $\overline{x}|_g$,
\begin{equation}\label{equ:IsoMean}
    \overline{x}|_g=g^{-1}\bigl(\frac 1n\sum_{i=1}^ng(x_i)\bigl).
\end{equation}
}

\subsubsection{ The relationship between isomorphic addition and isomorphic mean}
The relationship  between isomorphic addition and isomorphic mean of $n$-tuple $\underline x$ is connected by isomorphic division type I:
\begin{equation}\label{equ:IsoMean-IsoAdd}
    \overline{x_i}|_g=\biggl[\bigl[x_1+\ldots +x_n\bigl]_g\div n\biggl]_g.
\end{equation}

\subsection{Properties of isomorphic weighted mean}
Among many already known properties, the following are some key properties of the mean.
\subsubsection{Property of mean value}
\slthm{\label{thm:XiofIsoWgtMean} With $n$-tuple $\,\underline{x}\wedge(\mathscr{I}_m\{g\}=[X~\sharp~U]_g)$ that are not all equal, where $g$ is continuous on interval $X$, and with weights $~\underline{p}$, there is an unique $\xi\in (\min\{\underline x\}, \max\{\underline x\})\subseteq X$ such that
\begin{equation}
    \xi=\overline{x,p}|_g=g^{-1}\bigl(\sum_{i=1}^np_i g(x_i)\bigl).
\end{equation}
}

This theorem reflects a basic property of above defined isomorphic weighted mean. It says the mean always exists on the interval provided $g$ is continuous. The proof is omitted. However in the definition of quasi-arithmetic mean in pp.266 of \cite{BULLENPS}, ~$g$ being continuous is prerequisite thus the existence of the mean is ensured. The proof is omitted.

\subsubsection{Property of monotonicity}
\slthm{\label{thm:MonotofIsoWgtMean} Any participating number($x_i$)'s value increasing will result in the isomorphic weighted mean's increasing.}

Thus any one's decreasing results in the mean's decreasing. This property is obvious with (\ref{equ:IsoWgtMean}) as $g$ and $g^{-1}$ are always of the same strict monotonicity.

\subsubsection{Invariant value with vertical scale and shift of DM}
\slthm{$\overline{x,p}|_h=\overline{x,p}|_g$, for  $h\in \mathbb{V}g$.}
\slprf{$\overline{x,p}|_h=g^{-1}\biggl(\Bigl(-C+\sum_{i=1}^np_i \big(kg(x_i)+C\big)\Bigl)/k\biggl)=\overline{x,p}|_g$.}

Obviously the isomorphic weighted mean can be graphically presented on the
corresponding isomorphic number-axis, which topic is not covered here.

\subsection{Comparison problems of isomorphic weighted mean}
Article \cite{LIUY2} is a dedicated work for the topics of ``isomorphic means'' in the scope of this ``isomorphic mathematical analysis system'', which covers ``isomorphic weighted mean of numbers'' in more detail(as well as Section \ref{sec:IsoMeanFunc}:``isomorphic mean of a function'' in this work). The Section 2.2 of \cite{LIUY2} provided already several classes of methods to compare 2 isomorphic weighted means of same tuple of numbers generated respectively by 2 different generator mappings. Please refer to that paper for the details of the topics.

\section{Isomorphic convexity of function}

On a dual-isomorphic system, because the graph of a function is inevitably distorted in some way due to uneven distribution of coordinates on the axes, the ``convexity'' of function on it is worthwhile discussing in both analytic and geometrical aspects.

\subsection{Dual-variable-isomorphic convex function}\label{sec:DVIConvFunc}

\sldef{\label{defin:DVIC} Given $(f\colon D\to M)\wedge\big(\mathscr{I}_m\{g,h\} = [X,Y~\sharp~U,V]_{g,h}\big)$, where $D,X,Y,U,V$ are intervals and both $g,h$ are continuous. If for arbitrary $x_1, x_2\in D$ and arbitrary $\lambda_1, ~\lambda_2\in (0, 1)$ satisfying  $\lambda_1+\lambda_2=1$, the following holds
\begin{equation}\label{equ:inequDVIconvex}
    f(\overline{x_i,\lambda_R}|_g)\leq \overline{f(x_i),\lambda_R}|_h~(i\in\{1,2\}),
\end{equation}
then $f$ is called a dual-variable-isomorphic(DVI) convex function generated by mapping $g, h$ (or by $\mathscr{I}_m\{g,h\}$)on interval $D$. If the inequality inverses, $f$ is called a dual-variable-isomorphic(DVI) concave function generated by mapping $g, h$(or by $\mathscr{I}_m\{g,h\}$) on interval $D$. $g,h$ are called the generator mappings of the DVI convex(concave) function.
}
The property of a function being DVI convex or DVI concave is said for the function to have a ``dual-variable-isomorphic(DVI-) convexity''(it's further referred to as ``Isomorphic convexity of function'' in this paper). From the definition, one can see there is generally one type of ``DVI-convexity'' corresponding to one permutation of arbitrary 2 strictly monotone continuous function $g, h$. When $g, h$ are both identities, the definition simplifies to that of a convex function(and concave function).

In article \cite{LIUY} there is already an equivalent definition but without the concept of isomorphic frame involved. That article(in the 2000's) is a dedicated work for the DVI convex function, which falls into the scope of this IMAS system.

\slrem{Dual-variable-isomorphic convexity of function is a special case of $(M,N)$-convexity of function \cite{GAUMANNG}, with the two isomorphic weighted means of the former generated by $g, h$ being respectively the choices of $M$-mean and $N$-mean of the latter.
}

\subsection{Geometrical meaning of dual-variable-isomorphic convexity}

\begin{figure}[htb]
\begin{center}
\includegraphics[viewport=0 0 877 624,scale=0.4]{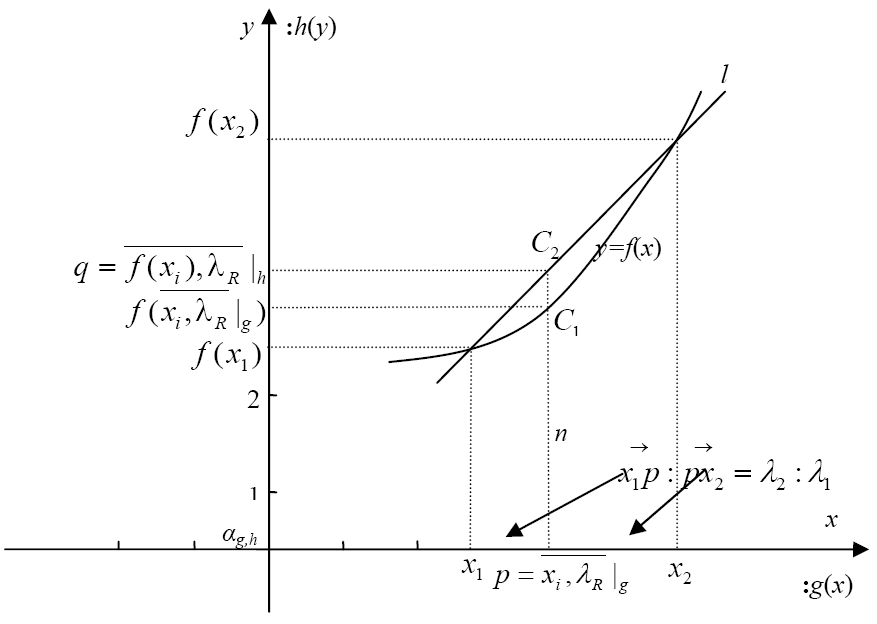}
\end{center}
\begin{center}
$f(\overline{x_i,\lambda_R}|_g)\leq \overline{f(x_i),\lambda_R}|_h~(i\in\{1,2\})$, $h$ strictly increases
\end{center}
\begin{illustration}{
\begin{center}
\label{illus:GeoDVIC}Geometrical meaning of dual-variable-isomorphic convex function
\end{center}
}\end{illustration}
\end{figure}

In Illustration \ref{illus:GeoDVIC} there is a dual-variable-isomorphic convex function $y=f(x)$ generated by mapping $g, h$ having its graph plotted on a dual-isomorphic system generated by $g, h$. The referential aux. system is not displayed here. On the horizontal axis, $\overline{x_i,\lambda_R}|_g(i\in\{1,2\})$ is represented by a point  $p\in (min\{x_1,  x_2\}, ~max\{ x_1,  x_2\})$ on the interval $D$ satisfying $\overrightarrow{x_1p}:\overrightarrow{{p}{x_2}}=\lambda_2:\lambda_1$ (reasoning from section \ref{sect:FixProp}). Similarly $\overline{f(x_i),\lambda_R}|_h(i\in\{1,2\})$ is represented by a point $q$ on the vertical axis. A secant line $l$ transfers the ratio $\lambda_2:\lambda_1$ from horizontal axis to vertical axis with the aid of imaginary parallel lines.

When $h$ is strictly increasing, $f(\overline{x_i,\lambda_R}|_g)\leq \overline{f(x_i),\lambda_R}|_h$ reflects that point $C_1$ is never above $C_2$. Thus any portion of graph of $f$ in between of $x_1, x_2$ will be ``convex'' towards below line $l$ since $p$ is arbitrary. While $h$ is decreasing, the portion of graph will be convex towards the upper of line $l$. In both cases the inequality means the graph of $f$ is convex to the opposite direction of the $y$ axis, to which the coordinates decrease. This is the geometrical meaning of the dual-variable-isomorphic convex function.

\slnot{On Illustration \ref{illus:GeoDVIC}, $f$ is explicitly said to be ``convex to the lower on the dual-isomorphic system''. If in Definition \ref{defin:DVIC} the inequality's sign inverses to ``$\geq$'', then graph of $f$ is convex to the other way and $f$ is said to be ``convex to the upper on the dual-isomorphic system'', in which case $f$ is a dual-variable-isomorphic concave function.
}

\subsection{A special case of dual-variable-isomorphic convex function -- geometrically convex function}\label{subsec:GeometricConvFunc}
The introduction of concept of ``geometrical convexity'' was made by J. Matkowski \cite{MATKOVSKI} in 1992. In monograph \cite{ZHANGXM} by Zhang, there is a simple version of definition:

\sldef{\label{defin:GeoConFun}Let $f(x)\geq0$ be a function defined on $I$, if for $x_1, x_2$ the
following holds
\begin{equation}\label{equ:inequGeoConFun}
    f(\sqrt {x_1x_2})\leq \sqrt{f(x_1)f(x_2)},
\end{equation}
then $f(x)$ is said to be geometrically convex to the lower; if the inequality is inverse $f(x)$ is said to be geometrically convex to the upper. Here $f(x)$ is called geometrically convex function or geometrically concave function respectively.
}

Obviously Definition \ref{defin:GeoConFun} is a special case of Definition \ref{defin:DVIC} when the latter's $g(x)=\ln x(x>0), ~h(y)=\ln y(y>0)$, and $\lambda_1=\lambda_2=1/2$, except for that Definition \ref{defin:DVIC} won't allow $f$ to be defined on point $0$, which case is not critical. The geometrical meaning of the former can also be demonstrated on the Illustration \ref{illus:GeoDVIC}, in which point $p=\sqrt {x_1x_2}$ is the geometrical midpoint of $x_1, x_2$ and $q$ is $\sqrt{f(x_1)f(x_2)}$.

There are also other special cases, e.g. harmonic convex function \cite{WUSH} when $g(x)=1/x(x>0), ~h(y)=1/y(y>0)$, which are not covered in details here.

\subsection{Graphical comparison of mean values - an application of DVI-convexity}\label{sec:GraphCompMean}

Here the IMAS system introduces a class of problems that don't have a counterpart in the classical MA system -- the graphical comparison of mean values in the dual-isomorphic system. Just consider on the Illustration \ref{illus:GeoDVIC}, let $f$ be identity, the inequality turns to be
\begin{equation*}\label{equ:IsoWgtMeanCompare}
    \overline{x_i,\lambda_R}|_g\leq \overline{x_i,\lambda_R}|_h~(i\in\{1,2\}),
\end{equation*}
which right now stands for the geometrical comparison of two types of isomorphic weighted means.

\begin{figure}[htb]
\begin{center}
\includegraphics[viewport=0 0 727 497,scale=0.45]{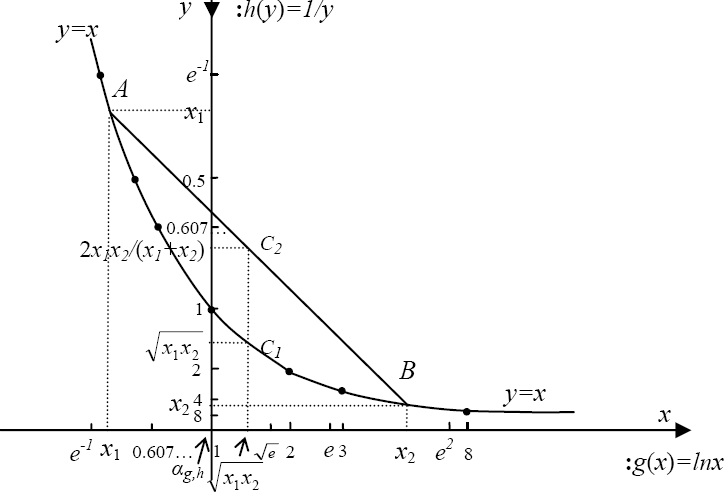}
\end{center}
\begin{center}
$G>H,~i.e. ~\sqrt{x_1x_2}>2x_1x_2/(x_1+x_2)$
\end{center}
\begin{illustration}{
\begin{center}
\label{illus:GraComp}Graphical comparison of mean values in a dual-isomorphic system
\end{center}
}\end{illustration}
\end{figure}

For instance, in Illustration \ref{illus:GraComp} let $g(x)=\ln x(x>0),~\lambda_1=\lambda_2=1/2, ~x_1\ne x_2$, ~$\overline{x_i,\lambda_R}|_g$ is the geometric mean of $x_1, x_2$, denoted by $G=\sqrt{x_1x_2}$. When $h(y)=1/y(y>0),~\overline{x_i,\lambda_R}|_h$ is the harmonic mean value of $x_1, x_2$, denoted by $H=2x_1x_2/(x_1+x_2)$. Next, on the dual-isomorphic system generated by $g, h$, the graph of $y=f(x)=x(x>0)$ is drawn, which is NOT a straight line. Then intuitively the result for $G>H$ can be observed in Illustration \ref{illus:GraComp}(notice the direction of the $y$ axis).

By changing the mapping $g$ and $h$, all types of isomorphic means (including arithmetic mean, square root mean etc.) can be compared by this graphical method. Of course, if $g=h$ the graph of $y=f(x)=x$ is always a straight line, same as is in the Cartesian coordinate system. In this case only the fact $\overline{x_i,\lambda_R}|_g= \overline{x_i,\lambda_R}|_h$ is manifest.

\subsection{Differential criteria of dual-variable-isomorphic convexity}
\subsubsection{The Criteria}
In article \cite{LIUY}, the differential criterion of dual-variable-isomorphic convex function is already elaborated and proved (in Chinese language). As quoted (in translation) below(note: no isomorphic frame involved yet therein).
\sllem{\label{lem:DiffCriDVIC}Let intervals $X, U, Y, V \subseteq\mathbb{R}, ~D \subseteq X$ and set $M\subseteq Y$. There are strictly monotone continuous bijections $g\colon X\to U, ~h\colon Y\to V$ and a function $f\colon D\to M$, and arbitrary $n$-tuple($n\geq2$) $\underline x\in D$, and positive $n$-tuple $\underline p$ satisfying $\sum_{i=1}^np_i=1$. $h\circ f$ and $g$ are derivable on $D$ and $g'\ne 0$,
\begin{enumerate}
  \item[1).] If $(h\circ f)'/g'$ is a constant on $D$, then
    \begin{equation*}
        \hspace{-1.0cm}f(\overline{x_i,p_R}|_g)= \overline{f(x_i),p_R}|_h~(i\in\{1,2,\ldots ,n\});
    \end{equation*}
  \item[2).] If $(h\circ f)'/g'$ is monotone on $D$ and among $(h\circ f)'/g', g, h$ there are odd number (1 or 3) of increasing functions, then
    \begin{equation*}
        \hspace{-1.0cm}f(\overline{x_i,p_R}|_g)\leq \overline{f(x_i),p_R}|_h~(i\in\{1,2,\ldots ,n\}),
    \end{equation*}
i.e. $f$ is a dual-variable-isomorphic convex function; If there are odd numbers (1 or 3) of decreasing functions among $(h\circ  f)'/g', g, h$ then the inequality inverses, i.e. $f$ is a dual-variable-isomorphic concave function. For both situations, if $(h\circ  f)'/g'$ is strictly monotone on $D$, then only if $x_1=\ldots=x_n$ the inequalities are equal;
  \item[3).] If $(h\circ f)'/g'$ and $g$ are both increasing or both decreasing on $D$, then $\varphi$, the dual-variable-isomorphic function of $f$ generated by $g, h$ is a convex function (convex to the lower on the Cartesian coordinate system); if one of $(h\circ  f)'/g'$ and $g$ is increasing and another decreasing, then $\varphi$ is convex to the upper.
\end{enumerate}
}

As a summary of above, the following table is the enumeration of all situations for function $f$'s dual-variable-isomorphic(DVI) convexity, depending on monotonicity of $(h\circ f)'/g', ~g, ~h$, and indicated by the inequality between $f(\overline{x_i,p_R}|_g)$ and $\overline{f(x_i),p_R}|_h$.
\begin{center}
\begin{tabular}{|p{1.8cm}|p{1.3cm}|p{1.1cm}|c|c|c|}\hline
\multicolumn{3}{|c|}{{\rule[-3mm]{0em}{10mm}}Monotonicity of functions} & \multicolumn{2}{|c|}{\tabincell {c}{$f(\overline{x_i,p_R}|_g)$ vs $\overline{f(x_i),p_R}|_h$ \\ ($f$'s DVI-Convexity)} } & $\varphi$'s convexity in R.A.system\\ \cline{1-5}
{\rule[-2mm]{0em}{7mm}} $(h\circ f)'/g'$ & $g$ & $h$ & $\leq$(convex) & $\geq$(concave)  & (Observer's direction)\\ \hline\hline
$\nearrow$ & $\nearrow$ & $\nearrow$ & $\checkmark$ &   & convex to lower \\ \hline
$\nearrow$ & $\nearrow$ & $\searrow$ &   & $\checkmark$ & convex to lower \\ \hline
$\nearrow$ & $\searrow$ & $\nearrow$ &   & $\checkmark$ & convex to upper \\ \hline
$\nearrow$ & $\searrow$ & $\searrow$ & $\checkmark$ &   & convex to upper \\ \hline
$\searrow$ & $\nearrow$ & $\nearrow$ &   & $\checkmark$ & convex to upper \\ \hline
$\searrow$ & $\nearrow$ & $\searrow$ & $\checkmark$ &   & convex to upper \\ \hline
$\searrow$ & $\searrow$ & $\nearrow$ & $\checkmark$ &   & convex to lower \\ \hline
$\searrow$ & $\searrow$ & $\searrow$ &   & $\checkmark$ & convex to lower \\ \hline
\end{tabular}
\footnotesize{~~Remarks:  $\nearrow$: increasing; ~$\searrow$: decreasing; ~R.A.System: Referential aux. system i.e. Cartesian system}
\begin{tables}\label{tab:DVI-conv}
\centering\normalsize{The enumeration of dual-variable-isomorphic convexity}
\end{tables}
\end{center}

The necessity theorem is also quoted from article \cite{LIUY}:
\sllem{\label{lem:DiffCriDVICness}Let intervals $X, U, Y, V \subseteq\mathbb{R}, ~D \subseteq X$ and set $M\subseteq Y$. There are strictly monotone continuous bijections $g\colon X\to U, ~h\colon Y\to V$ and a function $f\colon D\to M$, and arbitrary $n$-tuple($n\geq2$) $\underline x\in D$, and positive $n$-tuple $\underline p$ satisfying $\sum_{i=1}^np_i=1$. $h\circ f$ and $g$ are derivable on $D$ and $g'\ne 0$,
\begin{enumerate}
  \item[1).] If the following holds
    \begin{equation*}
    \hspace{-1.0cm}f(\overline{x_i,p_R}|_g)= \overline{f(x_i),p_R}|_h~(i\in\{1,2,\ldots ,n\}),
    \end{equation*}
   then $(h\circ f)'/g'$ is a constant on $D$.
  \item[2).] Let $f$ be a dual-variable-isomorphic convex function generated by $g, h$ on $D$. If $g, h$ are both increasing or decreasing, $(h\circ f)'/g'$ is increasing on $D$; if one of $g, h$ is increasing and another decreasing, $(h\circ f)'/g'$ is decreasing on $D$.
  \item[3).] Let $f$ be a dual-variable-isomorphic concave function generated by $g, h$ on $D$. If $g, h$ are both increasing or decreasing, $(h\circ f)'/g'$ is decreasing on $D$; if one of $g, h$ is increasing and another decreasing, $(h\circ f)'/g'$ is increasing on $D$.
\end{enumerate}
}

Lemma \ref{lem:DiffCriDVIC} is for judging a function's dual-variable-isomorphic convexity by analytic method, while Lemma \ref{lem:DiffCriDVICness} tells if knowing a function's dual-variable-isomorphic convexity then how its dual-isomorphic derivative will be, depending on $g, h$.

\subsubsection{Application of the differential criteria: the differential criteria of convex function \& geometrically convex function}

When $g, h$ are identities, the Lemma \ref{lem:DiffCriDVIC} \& Lemma \ref{lem:DiffCriDVICness} simplify to differential criteria of a convex function.
When $g(x)=\ln x(x>0), ~h(y)=\ln y(y>0), ~(h\circ f)'/g'$ is the elasticity of function $f$, $xf'(x)/f(x)$ ~($g, h$ are strictly increasing). The Lemmas simplify to differential criteria of a geometrically convex function, that is:
\slcor{When its elasticity $xf'(x)/f(x)$ is increasing, positive function $f(x)$ is a geometrically convex function; on the other hand if derivable function $f$ is a dual-variable isomorphic convex function generated by $g(x)=\ln x(x>0), ~h(y)=\ln y(y>0)$ (hence it's a geometrically convex function), its elasticity $xf'(x)/f(x)$ is increasing.
}

\subsubsection{Application of the differential criteria: a simple method for comparison of means}
When $f$ is identity mapping, Lemma \ref{lem:DiffCriDVIC} turns out to the comparison of 2 types of weighted isomorphic means. Especially when $p_i(i=1,\ldots ,n)=1/n$, it's further about the comparison of isomorphic means. Below Corollary is also quoted from article \cite{LIUY}:

\slcor{\label{cor:IsoMeanComp}Let $g, h$ be strictly monotone continuous and derivable functions on open interval $D$, and $g'\ne0, ~h'\ne0$. $h'/g'$ is monotone. There are arbitrary $n$-tuple($n\geq2$)$~\underline x\in D$, and positive $n$-tuple $\underline p$ adding up to 1,
\begin{enumerate}
  \item[1).] If there are odd numbers (1 or 3) of increasing functions among $h'/g', h, g$, then
      \begin{equation*}
      \hspace{-1.0cm}\overline{x_i,p_R}|_g\leq \overline{x_i,p_R}|_h~(i\in\{1,2,\ldots ,n\});
      \end{equation*}
  \item[2).] If there are odd numbers (1 or 3) of decreasing functions among $h'/g', h, g$, then
      \begin{equation*}
      \hspace{-1.0cm}\overline{x_i,p_R}|_g\geq \overline{x_i,p_R}|_h~(i\in\{1,2,\ldots ,n\}).
      \end{equation*}
\end{enumerate}
For both situations, the inequalities are equal only if $n$-tuple $\underline x$ are all equal.
}

In article \cite{LIUY}, with simple steps Corollary \ref{cor:IsoMeanComp} facilitates the proof of the power mean inequality:  \emph{For positive $n$-tuple($n\geq2$) $\underline x$ and real number $p>q$,
\begin{equation}\label{equ:PowerMean}
    \overline{x_i}|_{x^p}\geq \overline{x_i}|_{x^q},
\end{equation}
the inequality is equal only when $n$-tuple $\underline x$ are all equal.}

Another example by Corollary \ref{cor:IsoMeanComp} in article \cite{LIUY} is: $g(x)=\sinh(x), ~h(x)=\cosh(x)$. On $\mathbb{R}^+$, $g,h$ are increasing, and $h'/g'=(e ^x-e^{-x})/(e^x+e^{-x})=1-2/(e^{2x}+1)$ is increasing, thus
\begin{equation}
    \overline{x_i,p_R}|_{sinh x}\leq \overline{x_i,p_R}|_{cosh x}~(x_i>0).
\end{equation}

\subsubsection{Application of the differential criteria: building new inequality}
The differential criteria can also be used to construct inequalities. The following is a corollary from Lemma \ref{lem:DiffCriDVIC}, as an example, quoted from article \cite{LIUY}:
\slcor{Let $f(x)$ be a derivable function on open interval $I$. If $xf'(x)$ is increasing on $I$, then for arbitrary $n$-tuple($n\geq2$) $\underline x$ in interval $\sqrt[n]{
I}$ , the following holds
\begin{equation}
    nf(\prod_{i=1}^nx_i)\leq \sum_{i=1}^nf(x_i^{\,n}).
\end{equation}
Conversely, if $xf'(x)$ is decreasing on $I$,  the above inequality inverses.
}

\section{Isomorphic convex set}
\subsection{One-dimensional isomorphic convex set}
\sldef{\label{defin:1DIsoConSet}Given $D\vee\big(\mathscr{I}_m\{g\}= [X~\sharp~U]_g\big)$, if set $g(D)$ is an interval, then $D$ is called a one-dimensional isomorphic convex set generated by mapping $g$(or by $\mathscr{I}_m\{g\}$).
}


\slthm{\label{thm:1DIsoConSet}Given $D\vee\big(\mathscr{I}_m\{g\}= [X~\sharp~U]_g\big)$, the necessary and sufficient condition for $g(D)$ being an interval is that: For arbitrary $x_1, x_2\in D$ and $\lambda_1, \lambda_2\in (0,1)$ satisfying $\lambda_1+\lambda_2=1$, it holds that $\overline {x_i, \lambda_R}|_g\in D ~(i\in\{1,2\})$.
}
\slprf{ 
\begin{enumerate}
  \item[1).] Necessity: If $g(D)$ is an interval, because $g$ is strictly monotone, it is easy to prove that $\lambda_1g(x_1)+\lambda_2g(x_2)$ is in between of $g(x_1), ~g(x_2)$, therefore $\lambda_1g(x_1) +\lambda_2g(x_2)\in g(D)$, thus $g^{-1}\bigl(\lambda_1 g(x_1) +\lambda_2 g(x_2)\bigl)\in D$.
  \item[2).] Sufficiency: First, for arbitrary $u_1, u_2 \in g(D)$ there are $x_1=g^{-1}(u_1)\in D, ~x_2=g^{-1}(u_2)\in D$. Then for these $x_1, x_2$ we have $g^{-1}\bigl(\lambda_1g(x_1)+\lambda_2g(x_2)\bigl)\in D$, which means $\lambda_1g(x_1)+\lambda_2g(x_2)\in g(D)$ therefore $\lambda_1u_1+\lambda_2u_2\in g(D)$. Thus $g(D)$ is an interval. Proof done.
\end{enumerate}
}

\slrem{Definition \ref{defin:1DIsoConSet} and Theorem \ref{thm:1DIsoConSet} together are meaning that, for $D$ being a one-dimensional isomorphic convex set, the necessary and sufficient condition is that for arbitrary $x_1, x_2\in D$ and $\lambda_1, \lambda_2\in (0,1)$ satisfying $\lambda_1+\lambda_2=1$, there must be $\overline {x_i, \lambda_R}|_g\in D ~(i\in\{1,2\})$.
}

For an interval $D$ in the domain $X$ of $g$:
\begin{itemize}
  \item Obviously if $g$ is continuous on $D$ then $D$ is a one-dimensional isomorphic convex set generated by $g$ since $g(D)$ is an interval;
  \item Conversely if monotone $g$ is not continuous on $D$, then interval $D$ is not a one-dimensional isomorphic convex set generated by $g$.
\end{itemize}
This is to say:

\slrem{Intervals are not always one-dimensional isomorphic convex sets, but are
also depending on generator function $g$.
}
On the other hand, if $D$ is not an interval while $g(D)$ is, $D$ could still be a one-dimensional isomorphic convex set. That is to say:

\slrem{One-dimensional isomorphic convex sets are not always intervals.
}

\subsection{Two-dimensional isomorphic convex set}

\sldef{\label{defin:2DIsoConSet}Given  $T\vee\big(\mathscr{I}_m\{g_1,g_2\} = [X_1,X_2~\sharp~U_1,U_2]_{g_1,g_2}\big)$ where $U_1, U_2$ are intervals and $S=X_1\times X_2$ is the base frame. If for arbitrary 2 elements $p_1=(x_{11}, x_{21})\in T, ~p_2=(x_{12}, x_{22})\in T\big(x_{11}, x_{12}\in X_1, ~x_{21}, x_{22}\in X_2\big)$ and arbitrary $\lambda_1, \lambda_2\in (0, 1)$ satisfying $\lambda_1+\lambda_2=1$, there always exists a $p=(\overline{x_{1i}, \lambda_R}|_{g_1}, ~\overline{x_{2i}, \lambda_R}|_{g_2})\in T(i\in\{1,2\})$, then $T$ is called a two-dimensional isomorphic convex set generated by mapping  $g_1, g_2$(or by $\mathscr{I}_m\{g_1,g_2\}$) on $S$.
}  

\slrem{In Definition \ref{defin:2DIsoConSet}, $X_1, X_2$ are 2 one-dimensional isomorphic convex sets generated by $g_1, g_2$ respectively. This definition can be extended to ``$n$-dimensional isomorphic convex set'', if $n$ levels of one-dimensional isomorphic convex sets and mappings, e.g. $X_1, X_2, \ldots  X_n, g_1, g_2, \ldots  g_n(n\geq2)$ are involved.
}

\subsubsection{Geometrical meaning of two-dimensional isomorphic convex set}

With the conditions in Definition \ref{defin:2DIsoConSet}, a dual-isomorphic system generated by $g_1, g_2$ can establish. Then $S=X_1\times X_2$ is the domain.

\slrem{Obviously the geometrical meaning of 2-dimensional isomorphic convex set $T$ is that: It satisfies, for arbitrary 2 points  $p_1, p_2$ in it, the line segment connecting them is included by $T$ as the arbitrary point in such a line segment $p\in T$.
}

\begin{figure}[htb]
\begin{center}
\includegraphics[viewport=0 0 908 464,scale=0.45]{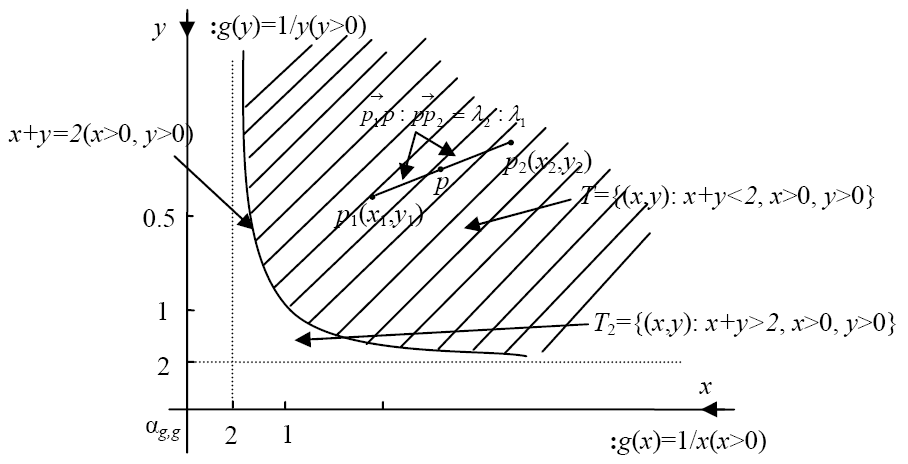}
\end{center}
\begin{center}
$T$ is a 2-dimensional isomorphic convex set.\\
$T_2$ is NOT a 2-dimensional isomorphic convex set,\\ but is a 2-dimensional convex set.
\end{center}
\begin{illustration}{
\begin{center}
\label{illus:2DIsoCnvSet}Two-dimensional isomorphic convex set
\end{center}
}\end{illustration}
\end{figure}

For instance, in Illustration \ref{illus:2DIsoCnvSet} there is a dual-isomorphic system generated by mapping $g(x)=1/x(x>0), ~g(y)=1/y(y>0)$, on which the curve $x+y=2(x>0, ~y>0)$ is convex to the upper and the shadowed area is represented by point set $T=\{(x, y)\colon  x+y<2, ~x>0, ~y>0\}$. We can see for arbitrary 2 points in $T$: $p_1=(x_1, y_1), ~p_2=(x_2, y_2)$ and arbitrary $\lambda_1, \lambda_2\in(0, 1)(\lambda_1+\lambda_2=1)$, the point $(\overline{x_i, \lambda_R}|_g, ~\overline{y_i, \lambda_R}|_g)$ on line segment $\overrightarrow{p_1p_2}$ is always included by set $T$, therefore the whole line segment $\overrightarrow{p_1p_2}$ is included by $T$ . So $T$ is a 2-dimensional isomorphic convex set generated by $g, g$.

To prove $T=\{(x, y)\colon x+y<2, ~x>0, ~y>0\}$ is a 2-dimensional isomorphic convex set generated by $g(x)=1/x(x>0), ~g(y)=1/y(y>0)$, the following equivalent problem needs to prove: There are positive numbers $x_1, y_1, x_2, y_2$ satisfying $x_1+y_1<2, ~x_2+y_2<2$, prove that for arbitrary $\lambda_1, ~\lambda_2\in (0, 1)$ satisfying $\lambda_1+\lambda_2=1$, ~$(x_1x_2)/(\lambda_1x_2+\lambda_2x_1)+(y_1y_2)/(\lambda_1y_2+\lambda_2y_1)<2$ holds. (Proof is omitted.)

Besides, consider $T_2=\{(x, y)\colon x+y>2, ~x>0, ~y>0\}$ which is a 2-dimensional convex set, but in the Illustration \ref{illus:2DIsoCnvSet} obviously it is not a 2-dimensional isomorphic convex set generated by $g, g$.

\subsubsection{Two-dimensional geometrically convex set -– a special case}\label{subsubsec:2DGeoConSet}

\sldef{\label{defin:2DGeoConSet}Let $T\subseteq \mathbb{R}^+\times \mathbb{R}^+$. If for arbitrary $p_1=(x_1, y_1)\in T, ~p_2=(x_2, y_2)\in T$, there always exists $p=(\sqrt{x_1x_2},\sqrt{y_1y_2})\in T$,
then $T$ is called a 2-dimensional geometrically convex set.
}

Obviously, Definition \ref{defin:2DGeoConSet} is a special case of Definition \ref{defin:2DIsoConSet} with $g_1(x)=\ln x(x>0), ~g_2(y)=\ln y(y>0)$ and $\lambda_1=\lambda_2=1/2$, i.e. 2-dimensional geometrically convex set is a special case of two-dimensional isomorphic convex set.

$\mathbb{R}^+\times \mathbb{R}^+$ itself is a 2-dimensional geometrically convex set. Below is the reasoning for another example $T=\{(x,y)\colon x-y>2, ~x>2, ~y>0\}$ :

Draw the graph of $y=x-2(x>2)$ on a dual-isomorphic system generated by $g(x)=\ln x(x>0), ~g(y)=\ln y(y>0)$. The dual-isomorphic derivative of $y$ is $\mathrm{d}\ln(x-2)/\mathrm{d}\ln x=x/(x-2)$ which is decreasing. Meanwhile $g(x), ~g(y)$ are increasing. i.e. among these 3 functions, there is  1 number (odd) of decreasing function. Based on Lemma \ref{lem:DiffCriDVIC}, the curve  $y=x-2(x>2)$ is convex to the upper on the dual-isomorphic system, therefore intuitively it can be asserted that the area below curve $y=x-2(x>2)$, which is represented by $T=\{(x,y)\colon x-y>2, ~x>2, ~y>0\}$, is a two-dimensional geometrically convex set, i.e. for arbitrary 2 elements $(x_1, y_1), ~(x_2, y_2)\in T$,  it holds $\sqrt{x_1x_2}-\sqrt{y_1y_2}>2$. (The analytic proof of the inequality is omitted.)

\section{Isomorphic mean of a function}\label{sec:IsoMeanFunc}

When the mean values of a function bonded on an isomorphic frame is concerned, the dual-variable-isomorphic(DVI) mean of a function can be similarly introduced.

\subsection{The definition of DVI mean of a function}
\sldef{\label{defin:IsoMeanOfFun} Let $f:D\to M$ be bounded on interval $D$ and $f\wedge\big(\mathscr{I}_m\{g,h\}=[X,Y~\sharp~U,V]_{g,h}\big)$. $g$ is continuous on $D$, interval $E=g(D)$, and $h$ is continuous on $[\inf M, \sup M]$. If there exists $M_\varphi \in V$, being the mean value of $\varphi\colon= (h\circ f \circ g^{-1})\colon E\to V$, then $h^{-1}(M_\varphi)\in Y$ is called the (dual-variable-) isomorphic mean(DVI mean) of $f$ on $D$ generated by $g, h$ (or generated by $\mathscr{I}_m\{g,h\}$), denoted by $\isomeanvalue{f}{D}{g,h}$, $\isomeanvalue{f}{\,}{g,h}$ or $M_f|_{g,h}$.\\
\begin{equation}\label{equ:IsoMeanValue0}
\begin{split}
    \isomeanvalue{f}{D}{g,h} &=
    h^{-1}\biggl(\frac{\int_{\scriptscriptstyle E}h\big(f(g^{-1}(u))\big)\mathrm{d}u}
    {\int_{\scriptscriptstyle E}\mathrm{d}u}\biggl). \\
    \biggl(~~&=h^{-1}\bigg(\frac{\int_{\scriptscriptstyle E}h\big(f(g^{-1}(u))\big)\mathrm{d}(-u)}
    {\int_{\scriptscriptstyle E}\mathrm{d}(-u)}\bigg).~~\biggl)
\end{split}
\end{equation}
It's stipulated that $\forall a\in D$, $\isomeanvalue{f}{[a,a]}{g,h}=f(a)$. ~$g,h$ are called the independent variable's generator mapping or dimensional mapping(IVDM), and the dependent variable's generator mapping or dimensional mapping(PVDM) of the isomorphic mean respectively.
}
\slrem{If $M_\varphi$ exists, then there must be $M_\varphi\in h([\inf M, \sup M])\subseteq V$, because $h([\inf M, \sup M])=[\inf N, \sup N]$ ($N=h(M)$) in which $M_\varphi$ must resides.} 

\slrem{In simple words, isomorphic mean of a function $f$ is the inverse image of the mean of $\varphi(f:g,h)$ by the PVDM.}
Intended to be the extension of mean of a function on an interval, this article defines the isomorphic mean of an ``ordinary'' function rather on an interval, than on a general real number set. Another ad-hoc requirement in addition to the ``bonding'' is $[\inf M, \sup M]\subseteq Y$.

\slnot{Taking $[a, b]$ as the form of $D$, the formula (\ref{equ:IsoMeanValue0}) turns into
\begin{equation}\label{equ:IsoMeanValue}
\begin{split}
    \isomeanvalue{f}{[a,b]}{g,h} &= h^{-1}\biggl(\frac{1}{g(b)-g(a)}\int_{g(a)}^{g(b)}h\bigl(f(g^{-1}(u))\bigl)\mathrm{d}u\biggl),\\
    \biggl(~~ &= h^{-1}\bigg(\frac{1}{g(a)-g(b)}\int_{g(b)}^{g(a)}h\bigl(f(g^{-1}(u))\bigl)\mathrm{d}u\bigg) ~~\biggl)
\end{split}
\end{equation}
with which whether $g(b)>g(a)$ is disregarded since the result of the formula is the same with that having $g(a)$ and $g(b)$ exchanged their positions. And if interval $D$ has infinite endpoint(s), the limit form of (\ref{equ:IsoMeanValue}) may be considered, with $a$ and/or $b$ approaching infinity.
}

Whereas if $f$ is not bounded, which falls out of scope of Definition \ref{defin:IsoMeanOfFun}, as long as the value of (\ref{equ:IsoMeanValue0}) exists, it could be considered the generalized isomorphic mean of an unbounded function, in which case $h$ is only to be continuous on $M$.

\subsection{\label{subsec:ExiCompatIsoMean}The existence of \texorpdfstring{$M_f|_{g,h}$}{Mf|g,h} and tolerance of the definition}
With Definition \ref{defin:IsoMeanOfFun}, the existence of $M_f|_{g,h}$ depends on the existence of the mean of $\varphi(f:g,h)$. The following theorem is sufficient but not necessary.
\slthm{\label{thm:IsomeanExist1}The $\isomeanvalue{f}{D}{g,h}$ exists with any applicable $g,~h$, if $f$ is continuous on a close interval $D$.
}
\slprf{Such $f,g,h$ imply a continuous $\varphi(f:g,h)$ on a close interval $E$ with a convex range $h(M)$. This means the numerator integral in the formula (\ref{equ:IsoMeanValue0}) exists and the dominator is non-zero real, such $M_\varphi$ exists in $h(M)\subseteq V$. These lead to an unique $\isomeanvalue{f}{D}{g,h}$ in $M\subseteq Y$.
}

While one expects for $f$ being continuous, the range $M$ being convex(an interval) and $M_f|_{g,h} \in M$, among others the definition however allows for
\begin{itemize} \setlength{\itemsep}{-0.2em}
  \item[(1).] $f$ is not continuous with jump discontinuity, but $M_f|_{g,h}$ exists. e.g.:
   \begin{equation} 
    \begin{split}
        f(x)&=\left\{
        \begin{aligned}
        1, ~~(x&\in [0,1]) \\
        3. ~~(x&\in (1,2]) \\
        \end{aligned} \right.\\
        g(x)&=2x, ~h(y)=y+2.
    \end{split}  
   \end{equation}
   With above, $M_f|_{g,h}=2$.
  \item [(2).] $f$ is not continuous with essential discontinuity, but $M_f|_{g,h}$ exists. e.g.
  If $f$ is a discontinuous but bounded Darboux function, especially being derivative of another continuous function $F$, and being Riemann integrable, then with the same $g,h$ as of above, $M_f|_{g,h}$ exists.
  \item [(3).] $M$ is not convex, and $M_f|_{g,h}\notin M$; but $M_f|_{g,h}\in Y$ as with above case (1).
  \item [(4).] $M$ is not convex, but $M_f|_{g,h}\in M\subseteq Y$.
  \item [(5).] The definition is not contradict with: $f$ is not continuous, $M_f|_{g,h}$ does not exist. e.g. $f(x)=\pi, x\in[0,2]$ and $x$ is irrational; $f(x)=3.14, x\in[0,2]$ and $x$ is rational, with the same $g,h$ as of above.
\end{itemize}
While it does not allow for $M_f|_{g,h}\in Y$ and $M_f|_{g,h}\notin [\inf M, \sup M]$, which will be proved later.

Another sufficient condition for the existence of isomorphic mean is:
\slthm{\label{thm:IsomeanExist2}The $\isomeanvalue{f}{D}{g,h}$ exists if $f$ is Riemann integrable on a close interval $D=[a,b]$ with an applicable $\mathscr{I}_m\{g,h\}$ such $g$ satisfies one of the following:
\begin{itemize} \setlength{\itemsep}{-0.2em}
  \item[1).] $g$ is derivable on $(a,b)$;
  \item[2).] $g$ is a convex or concave function on $[a,b]$;
  \item[3).] $g$ is absolutely continuous on any close sub-interval of $(a,b)$.
\end{itemize}
}
\slprf{Any of the 3 additional conditions adding to the strictly monotone and continuous $g$ ensures $f\circ g^{-1}$ is also integrable on $E$. With $h$ being continuous, $h\circ f\circ g^{-1}$ is Riemann integrable. Meanwhile the dominator is non-zero real, hence $M_\varphi$ exists in $h([\inf M, \sup M])\subseteq V$. This leads to an unique $\isomeanvalue{f}{D}{g,h}$ in $[\inf M, \sup M]\subseteq Y$.
}

Regarding the Riemann integrability of composite functions, the discussions can be founded in \cite{LUJ} and \cite{HUANGQL} which are applicable to e.g. $h\circ f\circ g^{-1}$ here. According to \cite{HUANGQL} any of additional conditions 1) 2) applying to $g$ will lead to condition 3) that will map the set of discontinuity points of $f$ on $[a,b]$ of Lebesgue measure 0 to a counterpart of $f\circ g^{-1}$ on $E$ of Lebesgue measure 0.

With above we claim that isomorphic means of a function are generally available in 2-D isomorphic frames with derivable DMs, for ORDINARY functions on close intervals, which (i) are Riemann integrable; (ii) may be somewhere discontinuous; (iii) have not to be monotone.

\subsection{An equivalent derivation of DVI mean of a function} \label{subsec:IsoMeanFunEquaDev}
Isomorphic mean of a function can also be derived through the limit of the isomorphic weighted mean of numbers expressed in integral form. In Definition \ref{defin:IsoMeanOfFun} assuming $D$ is a close interval $[a,~b]$, then due to monotone and continuous $g$, $E$ is a close interval $[g(a),~g(b)]$(disregarding whether $g(b)>g(a)$). Let $\tau=f\circ g^{-1}$. Do a $n$-tuple partitions of $E$ similar to is done with a definite integral. With each partition $i~(i\in{1,2,...,n}$), there is a corresponding value $\tau(\xi_i)$ for its tagged point $\xi_i$. Then the ratio of $\Delta u_i$($=u_i-u_{i-1}$, where $u_i$ is the end point of each partition towards $g(b)$, $u_{i-1}$ is the other end point of the same partition) over $g(b)-g(a)$, is taken as $\tau(\xi_i)$'s weight $w_i$. We compute the isomorphic weighted mean of $n$-tuple $\tau(\xi_i)$ generated by $h$, denoted by $M_{\tau(\xi_i)}$,
\begin{equation} 
\begin{split}
    M_{\tau(\xi_i)}  &= h^{-1}\biggl( \sum_{i=1}^{n}{\frac{u_i-u_{i-1}}{g(b)-g(a)}h(\tau(\xi_i))}\biggl) \label{equ:PartialSum}\\
        &= h^{-1}\biggl( \frac{1}{g(b)-g(a)}\sum_{i=1}^{n}{h(\tau(\xi_i))\Delta u_i}\biggl).
\end{split}  
\end{equation}

In order for $\tau(\xi_i)$ to enumerate all possible values of $\tau$ on $E$(such $f$ enumerates all on $[a,b]$), let $\|\Delta\|(=\max\{\Delta u_i\})\to0$, i.e. the partitions become infinitely thin, such
\begin{equation} 
\begin{split}
    \lim _{\|\Delta\|\to0}M_{\tau(\xi_i)} &=  h^{-1}\biggl(\frac{1}{g(b)-g(a)}\int_{g(a)}^{g(b)}h\big(\tau(u)\big)\mathrm{d}u\biggl)  \\
    &=  h^{-1}\biggl(\frac{1}{g(b)-g(a)}\int_{g(a)}^{g(b)}h\big(f(g^{-1}(u))\big)\mathrm{d}u\biggl). \label{equ:PartialSumLimit}
\end{split}
\end{equation}
This limit value is deemed as ``function $\tau$'s isomorphic weighted mean generated by mapping $h$''. It is just the same in form as (\ref{equ:IsoMeanValue}), that is function $f$'s dual-variable-isomorphic mean generated by $g,h$.

A special case of above is, the partitions are of all equal size, which is $\frac{1}{n}|g(b)-g(a)|$, such the weights are all $1/n$. Then (\ref{equ:PartialSumLimit}) can be deemed as an integral evolution of isomorphic mean of numbers, which has the same value if exists.

If $E$ is infinite or $D$ is open or half open, the limit forms of (\ref{equ:PartialSumLimit}) will be included cases of (\ref{equ:IsoMeanValue0}).

\subsection{Basic properties of DVI mean of a function}

\subsubsection{Property of intermediate value(IVP)}
\slthm{\label{thm:IVPofIsoMean}For a bounded $f\colon D\to M$($M$ is range) with its isomorphic mean $\isomeanvalue{f}{D}{g,h}$ in the applicable isomorphic frame $\mathscr{I}_m\{g,h\}$ , it holds
\begin{equation} \label{Inequ:IsoMeanFunGenIneq}
    \inf M \leq \isomeanvalue{f}{D}{g,h} \leq \sup M;
\end{equation}
especially if exists $\min\{M\}= \inf M$, $\max\{M\}=\sup M$, then
\begin{equation} \label{Inequ:IsoMeanFunGenIneq2}
   \min\{M\} \leq \isomeanvalue{f}{D}{g,h} \leq \max\{M\}.
\end{equation}
}
\slprf{Let $N=h(M)$. (i) if $M_\varphi\in N \subseteq V$, then due to the bijection $h$, $h^{-1}(M_\varphi)\in M \subseteq[\inf M, ~\sup M]$; (ii) if in general case $M_\varphi \in N$ is not to be considered, with (\ref{equ:PartialSum}), $h$ being continuous on $[\inf M, \sup M]$ by definition, and Theorem \ref{thm:XiofIsoWgtMean}, it holds
\begin{equation}
    \inf M \leq \min\{\tau(\xi_i)\}\leq M_{\tau(\xi_i)} \leq \max\{\tau(\xi_i)\}\leq \sup M
\end{equation}
while $\|\Delta\|\to0$. Thus $\inf M \leq \isomeanvalue{f}{D}{g,h} \leq \sup M$. For all cases, if there exists $\min\{M\}, \max\{M\}$, it holds $\min\{M\} \leq \isomeanvalue{f}{D}{g,h} \leq \max\{M\}$.}
Generally this Intermediate Value Property(IVP) or mean value property of an isomorphic mean does not require the corresponding intermediate value(s) being compulsory within the range of the function.

This IVP holds with ordinary functions. It could be fundamentally accepted as an attribute of ordinary functions on the background of isomorphic frames of continuous DMs, which will make isomorphic means have more coverage of extended unique means of a function than those concepts of IVP only holding with continuous or derivable functions.

The transforming nature of the background as a single bijection also makes the isomorphic mean an extended mean value of more straightforward and genuine origin, as directly from the mean value of an transformed function. As a fact, it naturally refines the development of mean values from simpler arithmetic mean into complicated ones, so that the isomorphic mean is a natural generalization of the simple mean.

In view of these, the isomorphic mean of a function is rather a differentiated concept of IVP than a repeated one with some others. (Also see ``Section 4: Isomorphic mean of a function vs Cauchy mean'' of article \cite{LIUY2}). 

\subsubsection{Property of monotonicity}
\slthm{\label{thm:IsoMeanMonotone}Let $f,m$ be defined on interval $D$, interval $D_2\subseteq D$. If
\begin{equation} \label{equ:MonotoneSubIneq}
m(x)\geq f(x) ~~\forall x\in D_2,
\end{equation}
and $m(x)=f(x)$ $\forall x\notin D_2$, then for any applicable $\mathscr{I}_m\{g,h\}$, it holds
\begin{equation} \label{equ:IsoMeanMonotone}
   M_m|_{g,h} \geq M_f|_{g,h}.
\end{equation}
If (\ref{equ:MonotoneSubIneq}) is strict, then (\ref{equ:IsoMeanMonotone}) is strict.
}
\slprf{While there are several cases according to the ways $D_2$ is located in $D$, we first assume both be close intervals, such $D$ can be expressed as $[a,b]$~ $(a<b)$. And arrange $D_2$ such $D_2=[a,c]$, $(a<c\leq b)$, i.e. $D_2$ locates to the left of $D$.  Let $\tau=f\circ g^{-1}$, $\mu=m\circ g^{-1}$, with the same way in Section \ref{subsec:IsoMeanFunEquaDev}, do a set of $n$-tuple partitions and tagged $\xi_i$ for $D_2=[a,c]$, another set of $k$-tuple partitions and tagged $\eta_j$ for $(c,b]$. Let
\begin{equation} 
\begin{split}
    M_{\tau(\xi,\eta)}  &= h^{-1}\biggl( \sum_{i=1}^{n}\frac{u_i-u_{i-1}}{g(b)-g(a)}h(\tau(\xi_i))+
    \sum_{j=1}^{k}\frac{v_j-v_{j-1}}{g(b)-g(a)}h(\tau(\eta_j))\biggl) \\
    M_{\mu(\xi,\eta)}  &= h^{-1}\biggl( \sum_{i=1}^{n}\frac{u_i-u_{i-1}}{g(b)-g(a)}h(\mu(\xi_i))+
    \sum_{j=1}^{k}\frac{v_j-v_{j-1}}{g(b)-g(a)}h(\mu(\eta_j))\biggl)
\end{split}
\end{equation}
\indent Above each putting in $h^{-1}$ are 2 partial sums. For convenience let's denote them $S_1$, $S_2$ for those inside $M_{\tau(\xi,\eta)}$ resp., and $T_1$, $T_2$ for those inside $M_{\mu(\xi,\eta)}$ resp., such
\begin{equation} \label{equ:IsomeanMonotone}
\begin{split}
    M_{\tau(\xi,\eta)}  = h^{-1}\big(S_1+S_2\big), ~~
    M_{\mu(\xi,\eta)}  = h^{-1}\big(T_1+T_2\big)
\end{split}
\end{equation}
\indent In the case $c=b$, $D_2=D$, $S_2,T_2$ do not exist. As $\forall x\notin D_2(i.e. ~x\in (c,b])$, $m(x)=f(x)$, such $S_2 \equiv T_2$. Meanwhile $\forall x\in D_2$, $m(x)\geq f(x)$, this makes $T_1\geq S_1 \Rightarrow T_1+T_2 \geq S_1+S_2$ with increasing $h$, or $T_1\leq S_1 \Rightarrow T_1+T_2 \leq S_1+S_2$ with decreasing $h$. In both cases $M_{\mu(\xi,\eta)} \geq M_{\tau(\xi,\eta)}$. When the partitions become infinitely thin, the inequality holds, thus it holds $M_m|_{g,h} \geq M_f|_{g,h}$.\\
 \indent For cases where a close $D_2$ is located other ways in a close $D$, the proof is analogous; Then at most we have 3 partial sums each in (\ref{equ:IsomeanMonotone}).\\
 \indent For at least one of $D_2$,$D$ and interval(s) of $D-D_2$ being open/half open intervals including above, the proofs always need to threat the endpoints carefully not to let tagged points be the endpoints that do not belongs, as fortunately finite points exclusion of tagged points does not affected the integral calculations;\\
 \indent For the cases $D_2$ and/or $D$ being infinite intervals, the proofs shall further consider the holding inequality on a selected partial finite interval, and let it hold under a limit of the endpoint(s). \\
 \indent Finally if (\ref{equ:MonotoneSubIneq}) is strict, then the corresponding inequalities within the proof including the result are strict.
}

\slcor{\label{thm:IsoMeanMonotoneCor}Let $f,m$ be defined on interval $D$, there are finite disjoint intervals $D_1,...,D_n\subseteq D$. If
\begin{equation} \label{equ:MonotoneSubIneqCor}
m(x)\geq f(x) ~~\forall x\in D_1\cup ... \cup D_n,
\end{equation}
and $m(x)=f(x)$ $\forall x\notin D_1\cup ... \cup D_n$, then for any applicable $\mathscr{I}_m\{g,h\}$, it holds
\begin{equation} \label{equ:IsoMeanMonotoneCor}
   M_m|_{g,h} \geq M_f|_{g,h}.
\end{equation}
If (\ref{equ:MonotoneSubIneqCor}) is strict, then (\ref{equ:IsoMeanMonotoneCor}) is strict.
}
This can be proved via Theorem \ref{thm:IsoMeanMonotone} and construction of $(n-1)$ bridging functions and passing the same inequality down total $(n+1)$ isomorphic means: $M_m|_{g,h} \geq...\geq M_f|_{g,h}$.

\subsubsection{Property of symmetry with endpoints of interval}
\slprop{(As implied by the definition,) The value of an isomorphic mean of a function on a close interval $[a,b]$ is invariant with $a,b$ exchanging their values in the formulae (\ref{equ:IsoMeanValue}).}

As a result, any function of the pair $(a,b)$ derived by (\ref{equ:IsoMeanValue}), maybe of different forms but of the equivalence of (\ref{equ:IsoMeanValue}), is symmetrical with $a,b$.
\slnot{For an existing $\isomeanvalue{f}{[a,b]}{g,h}$, also denote $\isomeanvalue{f}{[b,a]}{g,h}$ for the same isomorphic mean.
}
Therefore future for such $[a,b]$, $b\ge a$ is no longer a compulsory requirement.

\subsubsection{Invariant value with vertical scale and shift of dimensional mappings}
\paragraph{Invariant value with V-scaleshift of IVDM}
\slthm{\label{thm:IsoMeanInvar1} $M_f|_{m,h}=M_f|_{g,h}$ for $m\in \mathbb{V}g$.}
\slprf{
Let $u=m(x)$, $v=g(x)$, then $m^{-1}(u)=g^{-1}(v)$, ~$u=kv+C$, taking $[a,b]$ as $D$, thus
\begin{equation} \nonumber
\begin{split}
  M_f|_{m,h} &=h^{-1}\biggl(\frac{1}{kg(b)+C-kg(a)-C}\int_{kg(a)+C}^{kg(b)+C}h\big(f(m^{-1}(u))\big)\mathrm{d}u\biggl) \\
  &= h^{-1}\biggl(\frac{1}{kg(b)-kg(a)}\int_{g(a)}^{g(b)}h\big(f(g^{-1}(v))\big)\mathrm{d}(kv+C)\biggl) ~~(v=\frac{u-C}{k})\\
  &= h^{-1}\biggl(\frac{1}{g(b)-g(a)}\int_{g(a)}^{g(b)}h\big(f(g^{-1}(v))\big)\mathrm{d}v\biggl) = M_f|_{g,h}.
\end{split}
\end{equation}
While cases with other forms of $D$ are treated as holding limit forms of above.
}

\paragraph{Invariant value with V-scaleshift of PVDM}
\slthm{\label{thm:IsoMeanInvar2} $M_f|_{g,l}=M_f|_{g,h}$ for $l\in \mathbb{V}h$.}
\slprf{Taking $[a,b]$ as $D$, then
\begin{equation} \nonumber 
\begin{split}
  M_f|_{g,l} &= h^{-1}\biggl(\biggl(\bigg(\frac{1}{g(b)-g(a)}\int_{g(a)}^{g(b)}\big(kh\big(f(g^{-1}(u))\big)+C\big)\mathrm{d}u\bigg)-C\biggl)/k\biggl) \\
  &= h^{-1}\biggl(\biggl(\bigg(\frac{1}{g(b)-g(a)}\int_{g(a)}^{g(b)}kh\big(f(g^{-1}(u))\big)\mathrm{d}u+C\bigg)-C\biggl)/k\biggl) \\
  &= h^{-1}\biggl(\frac{1}{g(b)-g(a)}\int_{g(a)}^{g(b)}h\big(f(g^{-1}(u))\big)\mathrm{d}u\biggl)
  = M_f|_{g,h}.
\end{split}
\end{equation}
While cases with other forms of $D$ are treated as holding limit forms of above.
}

\paragraph{Invariant value with V-scaleshifts of both DMs}
\slcor{\label{cor:IsoMeanInvar12} $M_f|_{m,l}=M_f|_{g,h}$ for $m\in \mathbb{V}g$ and $l\in \mathbb{V}h$.}
It's a result of two-step process by previous 2 theorems.

\slnot{Also denote $M_f|_{\mathbb{V}g,\mathbb{V}h}$ for $M_f|_{g,h}$, $M_f|_{\mathbb{V}g}^{II}$ for $M_f|_g^{II}$, etc., i.e. for any IVDM or PVDM $g$, $\mathbb{V}g\Leftrightarrow g$ in denoting of isomorphic means of a function.}

\subsection{Sub-classing of DVI mean of a function} \label{subsec:IsoMeanClass}
There are 7(seven) typical sub-classes of isomorphic mean of a function, in correspondence with 7 special cases of DVI function:
\sldef{\label{Def:IsoMeanOfFunSub}With Definition \ref{defin:IsoMeanOfFun}, consider the following cases due to special $g,h,f$: 
\begin{enumerate}
  \item[1).] Let $g$ be identity, then $\varphi\colon=(h\circ f)\colon D\to V$. The isomorphic mean of $f$ is called the dependent-variable-isomorphic mean(PVI mean) of $f$ on $D$ generated by mapping $h$, or the isomorphic mean class I of $f$ on $[a, b]$ generated by $h$. It is denoted by $\isomeanvalue{f}{D}{h}$, or $M_f|_h$,
        \begin{equation}
            \isomeanvalue{f}{D}{h}
            = h^{-1}\Biggl(\frac{\int_{\scriptscriptstyle D}h\big(f(x)\big)\mathrm{d}x}
            {\int_{\scriptscriptstyle D}\mathrm{d}x}\Biggl)~(=\isomeanvalue{f}{D}{\mathbb{V}x,\mathbb{V}h}).
        \end{equation}
  If taking $[a, b]$ as $D$,
        \begin{equation}
            \isomeanvalue{f}{[a,b]}{h}
            =h^{-1}\biggl(\frac{1}{b-a}\int_a^bh\big(f(x)\big)\mathrm{d}x\biggl).
        \end{equation}
  \item[2).] Let $h$ be identity, then $\varphi\colon=(f\circ g^{-1})\colon E\to M$. The isomorphic mean of $f$ is called the independent-variable-isomorphic mean(IVI mean) of $f$ on $D$ generated by mapping $g$, or the isomorphic mean class II of $f$ on $D$ generated by $g$. It is denoted by $\isomeanvalueII{f}{D}{g}$, or $M_f|_g^{II}$,
        \begin{equation}
             \isomeanvalueII{f}{D}{g}
            = \frac{\int_{\scriptscriptstyle E}f(g^{-1}(u))\mathrm{d}u}
            {\int_{\scriptscriptstyle E}\mathrm{d}u}~(=\isomeanvalue{f}{D}{\mathbb{V}g,\mathbb{V}y}).
        \end{equation}
     If taking $[a, b]$ as $D$,
        \begin{equation}
            \isomeanvalueII{f}{[a,b]}{g}
            =\frac{1}{g(b)-g(a)}\int_{g(a)}^{g(b)}f\big(g^{-1}(u)\big)\mathrm{d}u.
        \end{equation}
     If $g$ is differentiable,
        \begin{equation}
            \isomeanvalueII{f}{[a,b]}{g}
            =\frac{1}{g(b)-g(a)}\int_a^bf(x)\mathrm{d}g(x).
        \end{equation}
  \item[3).] Let $Y=X, ~h=g$, then $\varphi\colon=(g\circ f\circ g^{-1})\colon E\to V$. The isomorphic mean of $f$ is called the same-mapping (dual-variable-)isomorphic mean(SDVI mean) of $f$ on $D$ generated by mapping $g$, or the isomorphic mean class III of $f$ on $D$ generated by $g$. It is denoted by $\isomeanvalueIII{f}{D}{g}$, $M_f|_{g,g}$, or $M_f|_{g}^{III}$,
      \begin{equation}
            \isomeanvalueIII{f}{D}{g}=
            g^{-1}\Biggl(\frac{\int_{\scriptscriptstyle E}g\big(f(g^{-1}(u))\big)\mathrm{d}u}
            {\int_{\scriptscriptstyle E}\mathrm{d}u}\Biggl)~(=\isomeanvalue{f}{D}{\mathbb{V}g,\mathbb{V}g}).
      \end{equation}
      If taking $[a, b]$ as $D$,
      \begin{equation}
            \isomeanvalueIII{f}{[a,b]}{g}
            =g^{-1}\biggl(\frac{1}{g(b)-g(a)}\int_{g(a)}^{g(b)}
            g\bigl(f(g^{-1}(u))\bigl)\mathrm{d}u\biggl).
      \end{equation}
  \item[4).] If in general case $h\ne g$, then $\varphi\colon=(h\circ f \circ g^{-1})\colon E\to V$. The isomorphic mean $\isomeanvalue{f}{D}{g,h}(=\isomeanvalue{f}{D}{\mathbb{V}g,\mathbb{V}h})$ is called (dual-variable-)isomorphic mean(DVI mean) of $f$ on $D$ generated by mapping $g, h$, or the isomorphic mean class IV of $f$ on $D$ generated by $g, h$, which formula is as (\ref{equ:IsoMeanValue0}), or as (\ref{equ:IsoMeanValue}) if taking $[a, b]$ as $D$.
  \item[5).] Let $f$ be identity, then $\varphi\colon=(h \circ g^{-1})\colon E\to V$. The isomorphic mean simplifies to the mean of one variable $x$. In this paper it is called the (dual-variable-)isomorphic mean of identity on $D$ generated by mapping $g, h$, or the isomorphic mean class V on $D$ generated by $g, h$. It is denoted by $\isomeanvalue{x}{D}{g,h}$, or $M_x|_{g,h}$,
      \begin{equation}
            \isomeanvalue{x}{D}{g,h}=
            h^{-1}\Biggl(\frac{\int_{\scriptscriptstyle E}h\big(f(g^{-1}(u))\big)\mathrm{d}u}
            {\int_{\scriptscriptstyle E}\mathrm{d}u}\Biggl)~(=\isomeanvalue{x}{D}{\mathbb{V}g,\mathbb{V}h}).
      \end{equation}
      If taking $[a, b]$ as $D$,
        \begin{equation} \label{equ:IsoMeanT5}
            \isomeanvalue{x}{[a,b]}{g,h}
            =h^{-1}\Biggl(\frac{1}{g(b)-g(a)}\int_{g(a)}^{g(b)}h\bigl(g^{-1}(u)\bigl)\mathrm{d}u\Biggl).
        \end{equation}
      If $g$ is differentiable,
        \begin{equation} \label{equ:IsoMeanT5dif}
            \isomeanvalue{x}{[a,b]}{g,h}
            =h^{-1}\Biggl(\frac{1}{g(b)-g(a)}\int_a^bh(x)\mathrm{d}g(x)\Biggl).
        \end{equation}
  \item[6).] Let $g, h$ be identities, then $\varphi\colon=f$. The isomorphic mean simplifies to the mean of the function. In this paper it is denoted by $\overline {{\displaystyle f}_{\scriptstyle {D}}}$, or $M_f$,
        \begin{equation}
            \overline {{\displaystyle f}_{\scriptstyle {D}}}=
            \frac{\int_{\scriptscriptstyle D}f(x)\mathrm{d}x}
            {\int_{\scriptscriptstyle D}\mathrm{d}x}~(=\isomeanvalue{f}{D}{\mathbb{V}x,\mathbb{V}y}).
        \end{equation}
      In the case $D$ is a close interval $[a,b]$,
        \begin{equation}
            \overline {{\displaystyle f}_{\scriptstyle {[a,b]}}}=\frac{1}{b-a}\int_a^bf(x)\mathrm{d}x.
        \end{equation}
  \item[7).] For monotone function $f$, its inverse function $f^{-1}$ is the dual-variable-isomorphic function of $f$ generated by mapping $f, f^{-1}$. Correspondingly the dual-variable-isomorphic mean of $f$ generated by $f, f^{-1}$ on a close interval is
      \begin{equation}
          \isomeanvalue{f}{[a,b]}{f,f^{-1}} = f\biggl(\frac{1}{f(b)-f(a)}\int_{f(a)}^{f(b)}f^{-1}(u)\mathrm{d}u\biggl).
      \end{equation}
\end{enumerate}
}
For convenience, while without confusions, we may use ``class I'' or ``Isomorphic mean class I'' for short name of Isomorphic mean class I of a function, and the similar for those of other classes later in this article, e.g. class II, class V.

For Theorem \ref{thm:IsoMeanInvar1}, especially when $g(x)=x$, then $M_f|_{m,h} = M_f|_{g,h} = M_f|_h$.

For Theorem \ref{thm:IsoMeanInvar2}, especially when $h(y)=y$, then $M_f|_{g,l} = M_f|_{g,h} = M_f|_g^{II}$.

For Corollary \ref{cor:IsoMeanInvar12}, especially when $g(x)=x$, $h(y)=y$, then $M_f|_{m,l} = M_f|_{g,h} = M_f$.
\subsection{Isomorphic mean class I of a function}\label{subsec:IsoMeanValFuncT1}
\begin{equation}
    \isomeanvalue{f}{[a,b]}{h}=M_f|_h=h^{-1}\biggl(\frac{1}{b-a}\int_a^bh\big(f(x)\big)\mathrm{d}x\biggl).
\end{equation}

As a special case, isomorphic mean class I also can be an equivalent derivation from isomorphic mean of numbers of $f(\xi_i)$, since in (\ref{equ:PartialSum}), let $g$ be identity mapping, then $\tau = f$. The class I can also be called the quasi-arithmetic mean of a function.

\subsubsection{Relationship btw. class I and isomorphic integral type I}
The relationship between isomorphic mean value class I of a function and the isomorphic integral type I is connected by isomorphic division type I:
\begin{equation*}
    I=h^{-1}\biggl[\int_a^b h\big(f(x)\big)\mathrm{d}x\biggl],
\end{equation*}
\begin{equation}
\isomeanvalue{f}{[a,b]}{h}
    =h^{-1}\Bigl[\frac{1}{b-a}\int_a^bh\big(f(x)\big)\mathrm{d}x\Bigl]
    =\big[I\div (b-a)\big]_h.
\end{equation}

\subsubsection{Some properties of isomorphic mean class I}
\slthm{\label{thm:IsoMeanC1PrdInvariant} $\isomeanvalue{H_{ss}\big(f:k,C\big)}{[ka+C,kb+C]}{h}=\isomeanvalue{f}{[a,b]}{h}$.}
\slprf{
\begin{equation} \nonumber
\begin{split}
    \isomeanvalue{H_{ss}\big(f:k,C\big)}{[ka+C,kb+C]}{h} &=  h^{-1}\biggl(\frac{1}{(kb+C)-(ka+C)}\int_{ka+C}^{kb+C}h\big(f((u-C)/k)\big)\mathrm{d}u\biggl) \\
    &\overset{(x=(u-C)/k)}{=} h^{-1}\biggl(\frac{1}{b-a}\int_a^bh\big(f(x)\big)\mathrm{d}x\biggl).
\end{split}
\end{equation}
}
This means isomorphic mean class I is invariant with the function's horizontal scale and shift.

Isomorphic mean class I is not a homogeneous mean in general cases, as with constant $k$, generally
\begin{equation}
  h^{-1}\biggl(\frac{1}{b-a}\int_a^bh\big(kf(x)\big)\mathrm{d}x\biggl) \neq
      kh^{-1}\biggl(\frac{1}{b-a}\int_a^bh\big(f(x)\big)\mathrm{d}x\biggl),
\end{equation}
i.e. $M_{kf}|_h\neq kM_f|_h$.

There are some special cases or instances of it, as discussed below.

\subsubsection{Arithmetic mean of a function}
When $h\in \mathbb{V}y$,
\begin{equation}
    \isomeanvalue{f}{[a,b]}{h}=\overline {f(x)}
    =\frac{1}{b-a}\int_a^bf(x)\mathrm{d}x.
\end{equation}
It is the arithmetic mean of $f$ on $[a, b]$.

\subsubsection{Geometric mean of a function} \label{subsubsec:GeoMeanFunc}
Let $h(y)=\ln y$, $g$ be identity and $f(x)>0$ in (\ref{equ:PartialSum}), it turns into
\begin{equation} 
    M_{f(\xi_i)} =  \exp\biggl( \frac{1}{b-a}\sum_{i=1}^{n}{\ln f(\xi_i)\Delta x_i}\biggl) =  \sqrt[b-a]{\prod_{i=1}^{n}{f(\xi_i)^{\Delta x_i}}}.
\end{equation}
It is the weighted geometric mean of all $f(\xi_i)$. Especially the partitions being all equal it's the geometric mean of all $f(\xi_i)$,
\begin{equation} 
    M_{f(\xi_i)} = \sqrt[n]{\prod_{i=1}^{n}{f(\xi_i)}}.
\end{equation}
Then its limit form (\ref{equ:PartialSumLimit}) is deemed as positive function ``$f$'s geometric mean''. Therefore when $h(y)=\ln y, ~f(x)>0$, the following is called the geometric mean (value) of $f(x)$ on $(a, b)$
\begin{equation}
    \isomeanvalue{f(x)}{(a,b)}{\ln y}
    =e^{\frac{1}{b-a}\int_{a^{\scriptscriptstyle +}}^{b^{\scriptscriptstyle -}}\ln f(x)\mathrm{d}x}.
\end{equation}
And the following is called the geometric mean (value) of $f(x)$ on $[a, b]$
\begin{equation}
    \isomeanvalue{f(x)}{[a,b]}{\ln y}
    =e^{\frac{1}{b-a}\int_a^b\ln f(x)\mathrm{d}x}.
\end{equation}

The geometric mean of $f(x)$ on $[a, b]$ is related to the ``geometric integral $I$ of $f(x)$ on $[a, b]$'' in the following equations:
\begin{eqnarray}\label{equ:ProductTypeI2}
    I &=& \exp \bigl(\int_a^b\ln f(x)\mathrm{d}x \bigl),\nonumber \\
        &=& \big(\isomeanvalue{f(x)}{[a,b]}{\ln y}\big)^{(b-a)}
\end{eqnarray}

It's obvious that geometric mean of a positive function is a homogeneous mean.

For a point $c\in[a,b]$, it holds that:
\begin{equation}
    {(\isomeanvalue{f}{[a,c]}{\ln y})}^{(c-a)} \cdot {(\isomeanvalue{f}{[c,b]}{\ln y})}^{(b-c)} = {(\isomeanvalue{f}{[a,b]}{\ln y})}^{(b-a)}. \label{eqn:GeoMeanFunMerge}
\end{equation}

Following are some special instances of geometric means on open intervals, in which the PVDM $\ln y$ is not defined on the $\inf M=0$($M$ is range of $f$). These only need an extra step of generalization by shrinking the domains $(a,b)$ of $f$ a little bit to $(a',b')$ so that $\inf M>0$, then taking the limit of the isomorphic mean with $a'\to a$ and/or $b'\to b$.
\begin{itemize}
  \item[(1).] The geometric mean of  $f(x)=x$ on $(0, b)$ : $\isomeanvalue{f(x)}{(0,b)}{\ln y}
        =e^{\frac{1}{b-0}\int_{0^{\scriptscriptstyle +}}^{b^{\scriptscriptstyle -}}\ln x\mathrm{d}x}=\frac be$.
  \item[(2).] The geometric mean of $f(x)=\sin x$ on $(0, \pi)$:
        $\isomeanvalue{\sin x}{(0,\pi)}{\ln y}
        =e^{\frac{1}{\pi}\int_{0^{\scriptscriptstyle +}}^{{\pi}^{\scriptscriptstyle -}}\ln \sin x\mathrm{d}x}.$
    Because the improper integral $\int_{0^{\scriptscriptstyle+}}^{{\pi}^{\scriptscriptstyle -}}\ln \sin x\mathrm{d}x=-\pi\ln2$ \cite{KOSMALA}, such
    \begin{equation}
        \isomeanvalue{\sin x}{(0,\pi)}{\ln y}=e^{-\ln2}=\frac 12.
    \end{equation}
    It's easy to further conclude the geometric mean of $sin x$ on $(0, \frac{\pi}{2})$ is also $\frac 12$.
  \item[(3).] Function $tanx$ is unbounded on $(0, \frac{\pi}{2})$, however we consider the following as the generalized geometric mean for unbounded $tanx$ on $(0, \frac{\pi}{2})$, and we are able to work out the geometric mean value is 1.
    \begin{equation} 
       \isomeanvalue{\tan x}{(0,\pi/2)}{\ln y}
        =e^{\frac{2}{\pi}\int_{0^{\scriptscriptstyle +}}^{{\frac{\pi}{2}}^{\scriptscriptstyle -}}\ln \tan x\mathrm{d}x}=1.
    \end{equation}
  \item[(4).] Given a circle with diameter $d$, radius $r$. To compute the geometric mean of all parallel chords(e.g. in vertical direction): $\bar c$.
    The function of such chords can be written as $c=2 \sqrt{r^2-x^2}, (x\in(-r,r)$), such
    \begin{equation} 
        \bar c = \exp\bigl(\frac{1}{r-(-r)}\int_{(-r)^+}^{r^-}{\ln(2 \sqrt{r^2-x^2})\mathrm{d}x}\bigl) = \frac{4}{e}r = \frac{2}{e}d \approx 0.7358d.
    \end{equation}
\end{itemize}

\begin{figure}[htb]
\begin{center}
\includegraphics[viewport=0 0 639 212,scale=0.7]{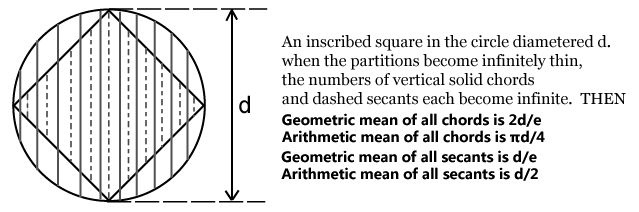}
\end{center}
\begin{illustration}{
\begin{center}
\label{illus:GeoArithRndSqr}Geometrical examples of geometric means of a function
\end{center}
}\end{illustration}
\end{figure}

Summarizing instance (1) and (4), and considering the homogeneity of geometric mean and its property by (\ref{eqn:GeoMeanFunMerge}) with instance (1), we have the following geometrical representation of some geometric means in Illustration \ref{illus:GeoArithRndSqr}.

\subsubsection{Harmonic mean of a function}
Let $h(y)=1/y$, $f(x)>0$ or $f(x)<0$,
\begin{equation}
    \isomeanvalue{f(x)}{[a,b]}{1/y}
    =\frac{b-a}{\int_a^b\frac{\mathrm{d}x}{f(x)}}
\end{equation}
is called the harmonic mean of $f(x)$ on $[a, b]$.

\subsubsection{Mean of power integral of a function}
Let $h(y)=y^p, ~(p\ne0)$, and $f(x)>0$ is continuous on $[a, b], ~b>a>0$,
\begin{equation}
    \isomeanvalue{f(x)}{[a,b]}{y^p}
    =\biggl(\frac{1}{b-a}\int_a^bf^p(x)\mathrm{d}x\biggl)^{\frac 1p},
\end{equation}
It is $f$'s $p$-order mean of power integral. Especially when $f(x)=x(p>-1,~ p\ne0, ~b>a>0)$,
\begin{equation}\label{equ:P-orderMean}
    \isomeanvalue{f(x)}{[a,b]}{y^p}
    =\biggl(\frac{b^{(p+1)}-a^{(p+1)}}{(p+1)(b-a)}\biggl)^{\frac 1p},
\end{equation}

\subsubsection{Function value on the midpoint of an interval}
When it happens that $f=h^{-1}$, then we know $f$ is monotone and $h=f^{-1}$, it follows that
\begin{equation}
    \isomeanvalue{f(x)}{[a,b]}{f^{-1}}
    =f\biggl( \frac{1}{b-a}\int_a^bx\mathrm{d}x \biggl)=f\big(\frac{a+b}{2}\big).
\end{equation}
i.e. the monotone function's value on the midpoint of an interval is a special case of  isomorphic mean class I of the function.

\subsection{Isomorphic mean class II of a function}
 Thanks to the introduction of isomorphic frame and DVI function, isomorphic mean class II of a function of case 2) of Definition \ref{Def:IsoMeanOfFunSub} is the sibling of class I, though these two are quite different in their special forms.
\begin{equation}
    \isomeanvalueII{f}{[a,b]}{g}=M_f|_g^{II}=\frac{1}{g(b)-g(a)}\int_{g(a)}^{g(b)}f\big(g^{-1}(u)\big)\mathrm{d}u.
\end{equation}
If $g$ is differentiable,
\begin{equation}
    \isomeanvalueII{f}{[a,b]}{g}=M_f|_g^{II}=\frac{1}{g(b)-g(a)}\int_a^bf(x)\mathrm{d}g(x).
\end{equation}

Isomorphic mean class II is not invariant with H-scaleshift of the function in general cases, as opposed to that of class I, displayed in Theorem \ref{thm:IsoMeanC1PrdInvariant}.

\subsubsection{Relationship btw. class II and isomorphic integral type II}
The relationship between isomorphic mean value class II of a function and the isomorphic integral type II is connected by isomorphic division type II:
\begin{equation*}
    I=g^{-1}\biggl[\int_{g(a)}^{g(b)}f\big(g^{-1}(u)\big)\mathrm{d}u\biggl],
\end{equation*}
\begin{eqnarray}\label{equ:RelationClassIIwiISOIntegII}
    \isomeanvalueII{f}{[a,b]}{g} &=& \frac{1}{g(b)-g(a)}\int_{g(a)}^{g(b)}f\big(g^{-1}(u)\big)\mathrm{d}u \nonumber\\
    &=& \bigg[I\div\bigl[b-a \bigl]_g\bigg]_g^{II}.
\end{eqnarray}

\subsubsection{Some properties of isomorphic mean class II}
\paragraph{Homogeneity}
\slthm{$\isomeanvalueII{V_{ss}\big(f:k,C\big)}{[a,b]}{g}=k\isomeanvalueII{f}{[a,b]}{g}+C$.}
Proof omitted. This is again opposed to the property of class I.

\paragraph{Conjugation of 2 isomorphic means class II}
\slthm{ Let $g,~f$ be 2 strictly monotone, differentiable functions on interval $[a,b]$, and  $~A=f(a), B=f(b), C=g(a), D=g(b)$. There exists 4 isomorphic means class II: $~E=M_f|_g^{II}, ~F=M_g|_f^{II}$ ($~E\neq A, ~E\neq B, ~F\neq C,~F\neq\ D$),$~G=M_f|_f^{II}=(A+B)/2$, $~H=M_g|_g^{II}=(C+D)/2$. Then the following hold:
\begin{equation}
\begin{split}
  1).~&\big(\overrightarrow{AE} \colon \overrightarrow{EB}\big) \cdot \big( \overrightarrow{CF} \colon \overrightarrow{FD}\big) = 1 \label{thm:IsoMeanC2Conjug}\\
  2). ~&\overrightarrow{GE} \colon \overrightarrow{AB} = \overrightarrow{FH} \colon \overrightarrow{CD}.
\end{split}
\end{equation}
}

\begin{figure}[htb]
\begin{center}
\includegraphics[viewport=0 0 400 170, scale=0.95]{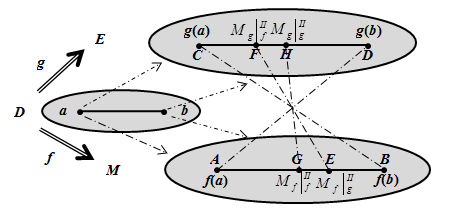}\\
$\big(\overrightarrow{AE} \colon \overrightarrow{EB}\big) \cdot \big( \overrightarrow{CF} \colon \overrightarrow{FD}\big) = 1$
\end{center}
\begin{illustration}{
\begin{center}
\label{illus:IsoMeanC2Conjug}Conjugation of 2 isomorphic means class II
\end{center}
}\end{illustration}
\end{figure}

\slprf{
\begin{equation} 
\begin{split}
        M_f|_g^{II} &=  \frac{1}{g(b)-g(a)}\int_a^bf(x)\mathrm{d}g(x)
         =  \frac{1}{g(b)-g(a)}\biggl( \big(f(x)g(x)\big)|_a^b - \int_a^bg(x)\mathrm{d}f(x)   \biggl) \\
        &=  \frac{1}{D-C}\big( B D -A C - F (B-A) \big)
        =  \frac{1}{D-C}\big( B (D-F) + A (F-C) \big)\\
        \Longrightarrow E &=  B \frac{\overrightarrow{FD}}{\overrightarrow{CD}} + A \frac{\overrightarrow{CF}}{\overrightarrow{CD}}
        \Longrightarrow E-A = B \frac{\overrightarrow{FD}}{\overrightarrow{CD}} + A \big( \frac{\overrightarrow{CF}}{\overrightarrow{CD}}-1 \big)
                            = B \frac{\overrightarrow{FD}}{\overrightarrow{CD}} - A \big( \frac{\overrightarrow{FD}}{\overrightarrow{CD}}\big).
        \nonumber
\end{split}
\end{equation}
It follows that $\overrightarrow{AE} \colon \overrightarrow{AB} = \overrightarrow{FD} \colon \overrightarrow{CD}$; Similarly $\overrightarrow{EB} \colon \overrightarrow{AB} = \overrightarrow{CF} \colon \overrightarrow{CD}$. Then it concludes with the result 1). The result 2) is a natural geometrical corollary of 1).}

Here the result 1) is said to be a kind of ``Conjugation of these 2 related isomorphic means class II($E$ and $F$)'' with exchanged original function and generator function. Illustration \ref{illus:IsoMeanC2Conjug} is the visual impression of the property when $f,g$ are of the same monotonicity.

\subsubsection{Simple special cases of isomorphic mean class II}
\begin{enumerate}
  \item[(1).] Let $g\in \mathbb{V}x$,
    \begin{equation}
        \isomeanvalueII{f}{[a,b]}{\mathbb{V}x}=\overline {f(x)}
        =\frac{1}{b-a}\int_a^bf(x)\mathrm{d}x.
    \end{equation}
It is the arithmetic mean of $f(x)$ on $[a, b]$.
  \item[(2).] When $f(x)$ is monotone and continuous on $[a, b]$, and $g(x)=f(x)$,
    \begin{equation}
        \isomeanvalueII{f}{[a,b]}{g} = \isomeanvalueII{g}{[a,b]}{f}
        =\frac{1}{f(b)-f(a)}\int_{f(a)}^{f(b)}u\mathrm{d}u=\frac12\bigl(f(a)+f(b) \bigl),
    \end{equation}
    i.e. the arithmetic mean of the monotone function's values at the 2 endpoints of an interval is a special case of isomorphic mean class II. 
  \end{enumerate}

\subsubsection{A special case: With a classical mean value theorem}\label{subsubsec:SpecCaseIsoMeanFuncT2}
\slthm{If $f\colon [a, ~b] \to \mathbb{R}$ is continuous and $g$ is an integrable function that does not change sign on $[a, b]$, then there exists $\xi$ in (a, b) such that
\begin{equation}
\begin{split}
  f(\xi)= \isomeanvalueII{f}{[a,b]}{G} ~~\biggl(=\int_a^b f(x) g(x)\mathrm{d}x /\int_a^bg(x)\mathrm{d}x\biggl) ~
         ~~(where~~G\in\int g).
\end{split}
\end{equation}
}
\slprf{With the conditions and the general case of well-known ``First mean value theorem for definite integrals'', there is a mean value $f(\xi)$ such
\begin{equation}
 \int_a^b f(x) g(x)\mathrm{d}x =f(\xi)\int_a^bg(x)\mathrm{d}x~~(\xi\in(a,b)).
\end{equation}
Let $G(x)=\int_a^x g(t)\mathrm{d}t+C$, then $G(x)$ is monotone and continuous, $G'(x)=g(x)$, thus
\begin{equation} 
\begin{split}
  f(\xi) =  \int_a^b f(x) g(x)\mathrm{d}x \div \int_a^b g(x)\mathrm{d}x = \frac{1}{G(b)-G(a)} \int_a^b f(x) \mathrm{d}G(x).
\end{split}
\end{equation}
Hence $f(\xi)$ is the isomorphic mean class II generated by the antiderivative of $g$.}
\slrem{While $f$ has not to be monotone, such $\xi$ has not to be unique, but the isomorphic mean $f(\xi)$ is unique.}
Especially when $g(x)=k,~G(x)=kx+C$, $G\in\mathbb{V}x$,
\begin{equation} 
\begin{split}
  f(\xi) = \isomeanvalueII{f}{[a,b]}{\mathbb{V}x} = \int_a^b kf(x)\mathrm{d}x \div \int_a^b k~\mathrm{d}x = \frac{1}{b-a} \int_a^b f(x) \mathrm{d}x
         =  \overline {{\displaystyle f}_{\scriptstyle {[a,b]}}}.
\end{split}
\end{equation}

\subsubsection{A special case: Elastic mean of a function}\label{subsubsec:SpecCaseIsoMeanFuncT2ElasMeanFun}
\paragraph{Introduction\\}
\sldef{The isomorphic mean class II of $f$ on $[a,b]~(b>a>0)$ generated by a logarithmic function e.g. $g(x)=\ln x(x>0)$,
    \begin{equation}\label{equ:logMeanFun}
        \isomeanvalueII{f}{[a,b]}{\mathbb{V}\ln x}
        =\frac{1}{\ln b-\ln a}\int_a^b\frac{f(x)}{x}~\mathrm{d}x,
    \end{equation}
is defined as the elastic mean of $f$ on $[a,b]$ in this paper.
}
\slrem{Any other logarithmic function is in $\mathbb{V}\ln x$.}
With $f(x)=x$,
\begin{equation} \label{equ:logMeanX}
        \isomeanvalueII{f}{[a,b]}{\ln x}
        =\frac{b-a}{\ln b-\ln a}.
\end{equation}

\paragraph{Why elastic mean of a function\\}
In Economics term, if $f$ is the elasticity \cite{TTElas} of another function $F$ such $f=xF'/F= (\mathrm{d}F/F)/(\mathrm{d}x/x)$, then it reflects the local relative change of $F$ against that of $x$, e.g. percentage change of demand (as an ``elastic'' response) against percentage change of price. Above $f(x)=x$ is the elasticity of $F(x)=ce^x (c > 0)$.

Let $K=F(b)\div F(a)$, and $k=b\div a$, it is easy to get
\begin{equation}\label{equ:ElastAgainstF}
        \isomeanvalueII{f}{[a,b]}{\ln x}
        =\ln \big(F(b)\div F(a)\big) \div \ln(b \div a)=\log_{~k} K. 
\end{equation}
On the other hand, let $M=F(x+\mathrm{d}x)\div F(x)$, $m=(x+\mathrm{d}x)\div x$, then
\begin{equation}
\begin{split}
        f(x) &=(\mathrm{d}F/F)/(\mathrm{d}x/x)=\mathrm{d}\ln F/\mathrm{d}\ln x \approx \Delta \ln F/\Delta \ln x \\
            &=(\ln F(x+\mathrm{d}x)-\ln F(x))\div(\ln(x+\mathrm{d}x)-\ln x)=\log_m M,
\end{split}
\end{equation}
i.e. $f$ is a ``micro-logarithm'' of 2 ``micro-multiplications'' for every local $x$, while $\isomeanvalueII{f}{[a,b]}{\ln x}$ is the logarithm of 2 overall multiplications over $[a,b]$.

Therefore the elastic mean of $f$ is actually the ``average of the elasticity $f$''. It is very similar to arithmetic mean of a function $f$ computed via the Newton-Leibniz formula:
\begin{equation}
        M_f = {\big(F(b)-F(a)\big)}\div {(b-a)},
\end{equation}
where $F$ is now the anti-derivative of $f$.

\paragraph{Relationship with isomorphic integral type II\\}
Furthermore, let $I$ be the related isomorphic integral type II of $f$ generated by $g(x)=\ln x$:
\begin{equation} 
    I   = \exp\biggl[\int_a^bf(x)\mathrm{d}\ln x\biggl].\nonumber \\
\end{equation}
With (\ref{equ:RelationClassIIwiISOIntegII}), we have
\begin{eqnarray}\label{equ:RelationClassIIwiISOIntegII2}
    \isomeanvalueII{f}{[a,b]}{\ln x} &=& \bigg[I\div\bigl[b-a \bigl]_{\ln x}\bigg]_{\ln x}^{II}\nonumber \\
    &=& \ln \big(I\big) \div \ln(b \div a)~.\nonumber
\end{eqnarray}
With (\ref{equ:ElastAgainstF}), it follows that
\begin{equation}
        I = F(b)\div F(a).
\end{equation}
From there we conclude that the so-called elastic integral (Notation \ref{not:ElasticInteg}) of $f$ on $[a,b]$ can be used to compute the overall multiplication of $F$(as is such related) from $a$ to $b$.

\subsubsection{Instances of isomorphic mean class II}
\paragraph{Elastic mean of \texorpdfstring{$\tan x$\\}{tan x}}
A special instance of elastic mean is about unbounded $tanx$ on $(0, \frac{\pi}{2})$, as following limit form:
\begin{equation}
\begin{split}
   \isomeanvalueII{\tan x}{(0,\pi/2)}{\ln x} &= \lim_{x\to 0,y\to\frac{\pi}{2}} \frac{1}{\ln \frac{\pi}{2}-\ln x}\int_0^y \frac{\tan(t)}{t}\mathrm{d}t.
\end{split}
\end{equation}
The numerator part is an improper definite integral, and the denominator part is approaching $+\infty$. We transform it with $~x=\rho\cos\theta,~y=\frac{\pi}{2}+\rho\sin\theta ~~(\rho>0, ~-\frac{\pi}{2}<\theta<0)$, and apply L~'Hopital's ~rule for twice in the following:
\begin{equation}
\begin{split}
   \isomeanvalueII{\tan x}{(0,\pi/2)}{\ln x} &= \lim_{\rho \to 0} \frac{1}{\ln \frac{\pi}{2}-\ln (\rho\cos\theta)}\int_0^{\frac{\pi}{2}+\rho\sin\theta} \frac{\tan(t)}{t}\mathrm{d}t \\
    &= \lim_{\rho \to 0}\frac{1}{-\cos\theta\frac{1}{\rho\cos\theta}}\cdot\frac{\tan(\frac{\pi}{2}+\rho\sin\theta)}{\frac{\pi}{2}+\rho\sin\theta}\cdot \sin\theta
    = \frac{2}{\pi}.
\end{split}
\end{equation}

\paragraph{Elastic mean of power function\texorpdfstring{\\}{}}
Let $f(x)=x^p, p\ne0, b>a>0, x\in [a, b], g(x)=\ln x(x>0)$,
\begin{equation} 
\begin{split}
    \isomeanvalueII{f}{[a,b]}{\ln x} =  \frac{1}{\ln b-\ln a}\int_{\ln a}^{\ln b}f(e^u)\mathrm{d}u 
    =  \frac{f(b)-f(a)}{\ln f(b)-\ln f(a)},
\end{split}
\end{equation}
which is the logarithmic mean of $f(a)$ and $f(b)$.

\paragraph{Power function's isomorphic mean class II generated by \texorpdfstring{$1/x(x>0)$\\}{1/x(x>0)}}
Let $f(x)=x^p, p\ne1, b>a>0, x\in [a, b], g(x)=1/x(x>0)$,
\begin{equation} 
\begin{split}
    \isomeanvalueII{f}{[a,b]}{1/x} = \frac{1}{(1/b)-(1/a)}\int_{1/a}^{1/b}(1/u)^p\mathrm{d}u 
    = \frac{ab(b^{p-1}-a^{p-1})}{(p-1)(b-a)}.
\end{split}
\end{equation}
When $p=2,~\isomeanvalueII{f}{[a,b]}{1/x}=ab$; when $p=3,~\isomeanvalueII{f}{[a,b]}{1/x}=\frac12ab(a+b)$.

While in the case $p=1, ~\isomeanvalueII{f}{[a,b]}{1/x}=\frac{1}{(1/b)-(1/a)}\int_{1/a}^{1/b}(1/u)\mathrm{d}u=\frac{ab(\ln b-\ln a)}{(b-a)}$, It is the product of $a, b$ divided by the logarithmic mean of $a, b$, it is also the limit of above $\frac{ab(b^{p-1}-a^{p-1})}{(p-1)(b-a)}$ when $p\to 1$.

\subsection{Isomorphic mean class III  \& IV of a function }
Isomorphic mean class III \& class IV are just the general forms of the Definition, being the combined form of class I \& class II. Their properties are mainly covered by previous sections.

\subsubsection{A special case of class III}
\slthm{\label{thm:IsoMeanC3ofIden} An isomorphic mean class III of $f(x)=x$ on [a,b] generated by $g$ equals the isomorphic mean of $a,~b$ generated by $g$ (the generalized $g$-mean).
}
\slprf{
\begin{equation}
    \isomeanvalueIII{x}{[a,b]}{g} =g^{-1}\biggl(\frac{1}{g(b)-g(a)}\int_{g(a)}^{g(b)}g\bigl(g^{-1}(u)\bigl)\mathrm{d}u\biggl)
    =g^{-1}\big(\frac{g(a)+g(b)}{2}\big).
\end{equation}}

\subsection{Isomorphic mean class V} 
Isomorphic mean class V in (\ref{equ:IsoMeanT5}) can be deemed as a special mean of a single variable, e.g. in the form of
\begin{equation} \nonumber
    \isomeanvalue{x}{[a,b]}{g,h}=h^{-1}\Biggl(\frac{1}{g(b)-g(a)}\int_a^bh(x)\mathrm{d}g(x)\Biggl).
\end{equation}
While traditionally without isomorphic frame in mind, it is not so worthwhile to discuss, since on $[a,b]$ $y=x$ always has a mean value $\frac{1}{2}(a+b)$. (Recalling Section \ref{sec:GraphCompMean}, the DVI-convexity of $y=x$ is somehow meaningful too.) The class V is generally not a mean value of an ordinary function. It is further sub-classifiable with $g\in \mathbb{V}x$ or $h\in \mathbb{V}y$, corresponding to class I or class II of $f(x)=x$, i.e. of a single variable, which are the closest concepts to so-called ``class 0'': the isomorphic mean of numbers.

\subsubsection{Composite class V}
\sldef{With a strictly monotone and continuous $f$, and applicable $g,h$,
\begin{equation}
M_x|_{g,H}=f^{-1}(M_f|_{g,h})
\end{equation}
is called a composite isomorphic mean class V generated by $g,h,f$, where $H:=h\circ f$.}
Above $H$ is strictly monotone and continuous. Composite class V is a special case of class V, but from the view of a normal class V of $t=M_x|_{g,H}$, the PVDM $H$ can be decomposed to $h\circ f$ whereby class V can be related to mean values of more monotone functions, i.e.
$f(t)=M_f|_{g,h}.$


\subsubsection{Generation of bivariate means by class V}\label{sec:GenBivarMean}
For a close interval $[a,b]$ and applicable $g,h$(and a monotone $f$ as with composite class V), the $M_x|_{g,h}$ (or $M_x|_{g,H}$) is clearly a sort of mean value of $a,b$. Moreover with the ``property of symmetry with endpoints of interval'', such bivariate mean is symmetric. Below are 2 examples among possible others.

\paragraph{Bivariate means regarding trigonometric functions\texorpdfstring{\\}{}}
By choosing $g(x)=\sin x, ~h(y)=\cos y$, $[a,b]\subseteq[0,\pi/2]$,
\begin{equation}
\begin{split}
    \isomeanvalue{x}{[a,b]}{\sin x, \cos y} &= \arccos \biggl( \frac{1}{\sin b- \sin a} \int_a^b \cos x (\sin x)'\mathrm{d}x  \biggl) \\
        &= \arccos \biggl( \frac{b-a+\sin b\cos b-\sin a\cos a}{2(\sin b- \sin a)} \biggl).
\end{split}
\end{equation}
While $g,h$ are exchanged,
\begin{equation}
\begin{split}
    \isomeanvalue{x}{[a,b]}{\cos x, \sin y} &= \arcsin \biggl( \frac{1}{\cos b- \cos a} \int_a^b \sin x (\cos x)'\mathrm{d}x  \biggl) \\
        &= \arcsin \biggl( \frac{a-b+\sin b\cos b-\sin a\cos a}{2(\cos b- \cos a)} \biggl).
\end{split}
\end{equation}

\paragraph{A class of quasi-Stolarsky means\texorpdfstring{\\}{}}
Another example of bivariate mean class generated by the class V, in the case $g(x)=x^p(x>0, ~p\ne0), ~h(y)=y^q(y>0, q\ne0), a>0,b>0,a\ne b$. It's denoted and formulated by:
\begin{equation} \label{equ:quasiStolarsky}
    Q_{p,q}(a,b)= \biggl( \frac{p(b^{p+q}-a^{p+q})}{(p+q)(b^p-a^p)}\biggl)^{1/q}.
\end{equation}

It is very similar to the derivation of the Stolarsky means from the ``Cauchy's extended mean value theorem''(\cite{LEASCHOExtM2}, pp207) by a pair of power functions. This class has also different special cases(including limits of $Q_{p,q}(a,b)$ when $p\to0$ and/or $q\to0$, which are actually equivalent results as with replaced $g(x)=\ln x$ and/or $h(y)=\ln y$):
\begin{equation} \label{equ:quasiStolarsky2}
    \begin{split}
        Q_{p,q}(a,b)=\left\{
        \begin{aligned}
        &\big(\frac{a^p+b^p}{2}\big)^{1/p} &(p=q,pq\ne0), \\
        &\sqrt{ab} &(p=q=0), \\
        &\biggl(\frac{b^p-a^p}{p(\ln b-\ln a)}\biggl)^{1/p} &(p+q=0,pq\ne0), \\
        &\biggl(\frac{b^q-a^q}{q(\ln b-\ln a)}\biggl)^{1/q} &(p=0,q\ne0), \\
        &\exp\biggl(\frac{b^p\ln b-a^p\ln a}{b^p-a^p}-\frac{1}{p}\biggl) &(q=0,p\ne0), \\
        &\frac{2(a^2+ab+b^2)}{3(a+b)} &(p=2,q=1), \\
        &\sqrt[3]{a\cdot\frac{a+b}2\cdot b} &(p=-1,q=3). \\
        \end{aligned} \right.\\
    \end{split}  
   \end{equation}

Details of derivations and proofs are omitted. It features the class of power mean as a quite symmetric children form. Though conversion of this form to Stolarsky mean is easy by substituting power $s=p+q$, but from the perspective of isomorphic mean, it's the very balanced form with respect to $p$ and $q$.

\subsection{Isomorphic mean class VI \& VII of a function}
The case 6) of Definition \ref{Def:IsoMeanOfFunSub} is just the mean of $f$, as ``class VI'' of isomorphic mean. It is not further discussed here. For instance of case 7) as ``class VII'', if $f(x)=x^a (x\in [0, c], a>0)$ then $  \isomeanvalue{f}{[0,c]}{f,f^{-1}} = (\frac{a}{a+1})^ac^a$, which coefficient $(\frac{a}{a+1})^a$ is approaching $\frac{1}{e}$ when $a$ is approaching $+\infty$.

\subsection{The geometrical meaning of isomorphic mean of a function}
For this topic, refer to Section \ref{subsec:GeoMeanDVIMeanCauchyMeanThrm}, where we shall discuss it along with that of the ``Cauchy's mean value theorem''.

\subsection{\label{subsec:RelationIsoMeanCauchy}Isomorphic mean of a function vs Cauchy mean value}

For more detailed discussion of this topic, refer to article \cite{LIUY2} ``Section 4: Isomorphic mean of a function vs Cauchy mean value''. The following sub-sections are incorporated from that section, also serving for preparation of discussion of the geometrical meanings of the both, which are not covered in \cite{LIUY2}.

\subsubsection{About Cauchy mean value}
\slthm{\label{thm:CauchyMeanVal}Cauchy's mean-value theorem states that (\cite{LOSONLCAUCHYCOMP}, pp12): If $f$, $g$ are continuous real functions on $[x_1,x_2]$ which are differentiable in $(x_1,x_2)$, and $g'(u)\ne 0$ for $u\in(x_1,x_2)$, then there is a point  $t\in(x_1,x_2)$ such that
\begin{equation}
\frac{f'(t)}{g'(t)}=\frac{f(x_2)-f(x_1)}{g(x_2)-g(x_1)}.
\end{equation}}

Note such $t$ does not have to be unique in general cases(e.g. an definition in \cite{HILLE} says there is at least one such $t\in(x_1,x_2)$ are satisfactory). Hence none of $f(t)$,$g(t)$,$f'(t)$,$g'(t)$ has to be unique, i.e. there is not uniquely defined mean value of function(s) either.

For such $t$ to be unique, in \cite{LOSONLCAUCHYCOMP} there is further restriction, whereby there comes the definition of ``Cauchy mean value of two numbers'':
\sldef{\label{Def:CauchyMeanof2Number}Assuming now (with Theorem \ref{thm:CauchyMeanVal}) that $f'/g'$ is invertible we get
\begin{equation}
t=\big(\frac{f'}{g'}\big)^{-1}\bigg(\frac{f(x_2)-f(x_1)}{g(x_2)-g(x_1)}\bigg).
\end{equation}
This number is called the Cauchy mean value of the numbers $x_l, x_2$ and will be denoted by $t= D_{fg}(x_1, x_2)$.}
In this paper, we further say the mean value of such 2 numbers ``is generated by $f,g$''. (And in \cite{LOSONLCAUCHYCOMP}, further covered is a generalized form of above to Cauchy mean of $n$ numbers: $D_{fg}(x_1, x_2,...,x_n)$). But even with Definition \ref{Def:CauchyMeanof2Number}, neither $f(t)$ nor $g(t)$ can be deemed as well-defined mean value of a function, as they are symmetrically, temporarily depending on each other.

\subsubsection{Conversions btw. these two genres}
The Cauchy mean value class and certain classes of isomorphic means can be converted to each other, with some criteria. A reasonable point is that since Cauchy mean value deals with 2 functions while isomorphic mean does with 3, the corresponding isomorphic mean is related to or of class V, which is with simplest $f(x)=x$. Refer to article \cite{LIUY2} Section 4.2 \& 4.3.

\subsection{\label{sec:DifferenceIsoCauchy}Differentiations between isomorphic mean and Cauchy's mean value theorem}
Since the IVP property of isomorphic means applies to ordinary functions, while IVP of the Cauchy's mean value theorem applies to derivable functions, with the discussion of conversions between these 2 genres of mean values in \cite{LIUY2}, we claim that (i). Whenever there is a Cauchy mean value, there is a class V(with derivable DMs). (ii). Whenever there is an isomorphic mean unless it's a class V, it's not ensured there is a convertible Cauchy mean value by Theorem 4.5 in article \cite{LIUY2}. Provided therein are some examples of isomorphic means that can not convert to Cauchy mean values.

\subsubsection{About their coverage and categorical intersection}
For more details, please refer to the same topic discussed in Section 4.4.1 of article \cite{LIUY2}. Therein we claimed that \textbf{isomorphic means of a function cover broader range of unique mean values of an ORDINARY function.}(than Cauchy's mean values). Only the class V can be matched by the Cauchy mean value as a special derivation of the Cauchy's mean value theorem in term of uniqueness, whereas more other classes of isomorphic means conform to the ubiquitous theorem.

\subsubsection{About generator functions and identifications of unique means}
For details of the discussion, please refer to Section 4.4.2 of article \cite{LIUY2}. Therein we claimed that (comparing to Cauchy's mean values)\textbf{isomorphic means have better identifications as being well-defined unique extended means of a specific function}(to some extent owing to the well-balanced generator functions). It is with such good identifications, that isomorphic means also have had diversified classifications possibly.

Generally speaking: Isomorphic mean seems to be a concept of mean value with better origin and perspective, better identification and classification, more coverage and more natural generalization.

\subsection{About the geometrical meaning of Cauchy's mean value theorem}
In the ``conclusion \& vision'' part of article \cite{LIUY2}, we stated that ``one also can see that the geometrical meaning of Cauchy mean value in the coordinate system is different than that of isomorphic means''. We will attempt to discuss it in Section \ref{subsec:GeoMeanDVIMeanCauchyMeanThrm} along with that of DVI mean of a function, where we can make use of the graph of Illustration \ref{illus:ExampleOfIMAS}.

\subsection{The comparison problems of isomorphic means of a function} 
Isomorphic means of a function, with their rich sub-classes and abundant special cases, have the important property of IVP, among others. Their comparison inequality problems deserve a specialized long discourse. Please refer to article \cite{LIUY2} Section 5: ``The comparison problems of isomorphic means of a function'' for the details.

\subsection{Isomorphic mean of a function of n-1 variables}
The concept of isomorphic frame has been introduced on the $n$-dimensional basis. Based on this it is possible to extend the concept of isomorphic mean of a function for functions of total $n\geq2$ varibles(($n-1$) independent variables plus 1 dependent variable). The following is a premature definition for the concept.
\sldef{Let intervals $X_1,...,X_n,U_1,...,U_n\subseteq\mathbb{R}$, and $g_i\colon X_i\to U_i(i=1,...,n)$ be $n$ continuous and monotone bijections. Function $f:D\to M$ of ($n-1$) variables is bounded,  $g_n$ is continuous on $[\inf M, \sup M]$ and $f\wedge\mathscr{I}_m\{g_1,...,g_n\}$. $D$ and $E=\mathscr{I}_m\{g_1,...,g_{n-1}\}(D)$ each are connected sets and are measurable of $(n-1)$ dimensional hypervolume. ~If there exists $M_\varphi \in U_n$, being the mean value of $\varphi\colon=g_n\circ f(g_1^{-1},...,{g_{n-1}}^{-1})\colon E\to U_n$ on $E$, ~then $g_n^{-1}(M_\varphi)\in X_n$ ~is called the (all-variable-) isomorphic mean of $f$ on $D$ generated by $\mathscr{I}_m\{g_1,...,g_n\}$, denoted by $\isomeanvalue{f}{D}{g_1,...,g_n}$,
\begin{equation}\label{equ:IsoMeanValueN}
    \isomeanvalue{f}{D}{g_1,...,g_n} =
    g_n^{-1}\biggl(\frac{\int_{\scriptscriptstyle E}g_n\circ f(g_1^{-1},...,{g_{n-1}}^{-1})\mathrm{d}^{(n-1)}u}
    {\int_{\scriptscriptstyle E}\mathrm{d}^{(n-1)}u}\biggl). \\
\end{equation}
}
Upon seeing this definition, one can immediately ask if such defined mean value still has the intermediate value property(IVP) and if yes how to prove it? For studies of such problems and perfection of this MA concept we defer them to our future efforts.

\section{A simple example of IMAS in summary}\label{sec:Example}
Below we have an example of IMAS to examine foregoing topics briefly so as to have a big picture of the system. Let the 2 generator mappings be $g(x)=e^x, ~h(y)=y^p(y>0, ~p\ne0)$ for each topics of IMAS.
\begin{figure}[htb]
\begin{center}
\includegraphics[viewport=0 0 400 280,scale=0.9]{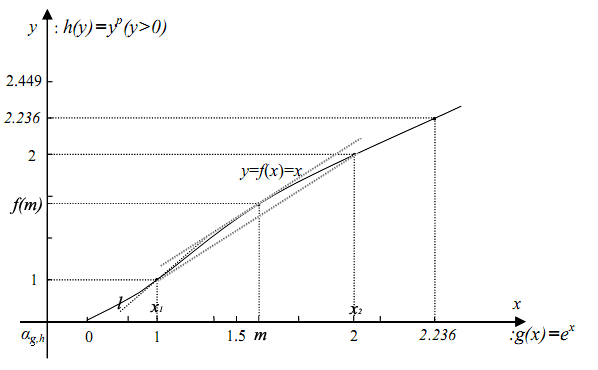}
\end{center}
\begin{center}
    $g(x)=e^x, ~h(y)=y^p ~(e.g. ~p=2)$
\end{center}
\begin{illustration}{
\begin{center}
\label{illus:ExampleOfIMAS}An example of IMAS on a dual-isomorphic system
\end{center}
}\end{illustration}
\end{figure}

\subsection{About isomorphic number}

$g\colon =x\mapsto e^x: \mathbb{R}\to\mathbb{R}^+$, for any $x\in\mathbb{R}$, $u=e^x\in\mathbb{R}^+$
 is the isomorphic number generated by $g$.

$h\colon =y\mapsto y^p: \mathbb{R}^+\to\mathbb{R}^+$, for any $y\in\mathbb{R}^+$, $v=y^p\in\mathbb{R}^+$ is the isomorphic number generated by $h$.

\subsection{About isomorphic number-axis and dual-isomorphic system}
In Illustration \ref{illus:ExampleOfIMAS} there is a dual-isomorphic system generated by mapping $g, h$, through steps in Notation \ref{not:EstabDualIsoSys}. The horizontal axis ($x$ axis) is an isomorphic number-axis generated by $g$, and the vertical axis ($y$ axis) generated by $h$. The referential aux. system is not shown.
On the $x$ axis, the domain of the axis: $\mathbb{R}$ is to the right side of $\alpha _{g,h}$, not including it, $\mathbb{R}^-$ is ``squeezed'' in between of 0 and $\alpha _{g,h}$. On the $y$ axis, the domain $\mathbb{R}^+$ is above $\alpha _{g,h}$, not including it.

\subsection{About fixed proportion division on isomorphic number-axis}
On the $x$ axis, there is the midpoint $m$ of point ``1'' and point ``2''. According to (\ref{equ:FixProportion}), its value is the isomorphic mean of 1 and 2 generated by $g$:
\begin{equation}
    m=\ln\bigl(\frac12(e^1+e^2)\bigl)~(\approx1.6201).
\end{equation}
Its distance to $\alpha _{g,h}$ is $e^m\approx5.0536$ measured in the unit of the aux. axis, which is not shown here.

\subsection{About function graph and dual-variable-isomorphic function}
Regarding a function $f\colon D\to M(D\subseteq\mathbb{R}, ~M\subseteq\mathbb{R}^+)$, let set $E=g(D)=e^D, ~N=h(M)=M^p$, $f$'s dual-variable-isomorphic function generated by $g, h$ is $\varphi\colon =(h\circ f\circ g^{-1})\colon E\to N$, or in traditional notation $v=(f(lnu))^p
~(u\in E, ~v\in N)$.

According to Theorem \ref{thm:GraphCorres}, $f$'s graph on the dual-isomorphic system is congruent to $\varphi$'s graph on the referential aux. system. On Illustration \ref{illus:ExampleOfIMAS}, the graph of function $y=x(x>0)$ is plotted(which is not a straight line). Its curve is congruent to the graph of function $v=(f(lnu))^p(u>1)$ on the referential aux. system, which is not shown here.

\subsection{About isomorphic arithmetic operations}\label{subsec:IsoOperExample}
Here we take a simple example for graphical representation of isomorphic addition on the $y$ axis, with $p=2$ on Illustration \ref{illus:ExampleOfIMAS}:
\begin{equation}
    \bigl[1+2\bigl]_h = \sqrt{1^2+2^2}=\sqrt5 ~(\approx 2.236).
\end{equation}
As we can see on the $y$ axis, the distance from $\alpha _{g,h}$ to 1 plus the distance from $\alpha _{g,h}$ to 2 equals that of from $\alpha _{g,h}$ to $\sqrt5$, near 2.236, or we can think the corresponding directed line segments are manoeuvred on the $y$ axis to represent the operation of isomorphic addition.

\subsection{About dual-isomorphic derivative and densities on the axes}
Regarding above function  $f\colon D\to M~(D\subseteq\mathbb{R}, ~M\subseteq\mathbb{R}^+)$, let set  $E=g(D)=e^D, ~N=h(M)=M^p$  and $y=f(x)(x\in D)$, $f$'s dual-isomorphic derivative at $x$ is:
\begin{equation}
    \bigl[f'(x)\bigl]_{g,h}
    =f'(x)\cdot\frac{h'(y)}{g'(x)}
    =f'(x)\cdot\frac{\rho_g(x)}{\rho_h(y)}
    =f'(x)\cdot\frac{(y^p)'}{(e^x)'}.
\end{equation}
When it applies to Illustration \ref{illus:ExampleOfIMAS}, $y=f(x)=x~(x>0), ~p=2,$
\begin{equation}\label{equ:XDensityonExample}
    \bigl[f'(x)\bigl]_{g,h}
    =1\cdot\frac{(y^p)'}{(e^x)'}
    =1\cdot\frac{px^{(p-1)}}{e^x}
    =1\cdot\frac{e^{-x}}{2^{-1}x^{-1}}.
\end{equation}
When $x=1$, the dual-isomorphic derivative is $2e^{-1}$, it is the slope of the tangent line $l$ of the function curve passing point $(1, 1)$ on the graph of $f$ on the illustration.

On the $x$ axis, the density at any $x$ is $\rho_g(x)=1/g'(x)=e^{-x}$, which is decreasing with respect to change of $x$. When $x$ is approaching $-\infty$(towards $\alpha _{g,h}$), $\rho_g(x)$ is approaching $+\infty$, as intuitively we can see more points are ``squeezed'' near $\alpha _{g,h}$. On the $y$ axis, the density at any $y$ is $\rho_h(y)=2^{-1}y^{-1}$, which is also decreasing with respect to $y$.

On the rightmost part of (\ref{equ:XDensityonExample}), the dual-isomorphic derivative of $y=f(x)=x~(x>0)$ is considered to be affected by following factors:
\begin{enumerate}
  \item[1).] It is based on $f$'s derivative at $x$ which is a constant 1;
  \item[2).] It is directly proportional to the density at $x$ on $x$ axis, which is $e^{-x}$;
  \item[3).] It is inversely proportional to the density at $y$ on $y$ axis, which is $2^{-1}x^{-1}$, corresponding to $x$ on $x$ axis.
\end{enumerate}

We can imagine that the graph of $y=f(x)=x(x>0)$, which is a straight line in Cartesian coordinate system, is in some way deformed by the ``strength'' of $g$ and $h$ in the dual-isomorphic system. If $g$ and $h$ are identities, the densities on both axes are constant 1, thus the graph is a straight line. Similar situation could be considered for an ordinary function $f$.

\subsection{About isomorphic integral}
The isomorphic integral type I of $f$ on interval $[a, b] \subseteq D$ generated by mapping $h$ is:
\begin{equation}
    I_1=h^{-1}\biggl[\int_a^b h(f(x))\mathrm{d}x\biggl]=\biggl(\int_a^b f^p(x)\mathrm{d}x\biggl)^{\frac1p}.
\end{equation}

The isomorphic integral type II of $f$ on interval $[a, b] \subseteq D$ generated by mapping $g$ is:
\begin{equation}
    I_2=g^{-1}\biggl[\int_a^bf(x)g'(x)\mathrm{d}x\biggl]=\ln\biggl(\int_a^bf(x)e^x\mathrm{d}x\biggl).
\end{equation}

\subsection{About geometrical meanings of DVI mean and Cauchy's mean value theorem}\label{subsec:GeoMeanDVIMeanCauchyMeanThrm}

\subsubsection{Geometrical meaning of DVI mean of a function}
In simple words, the DVI mean of a function $f$ generated by $g,~h$ is the inverse image of the mean of $\varphi(f:g,h)$ by mapping $h$. With the implement of dual-isomorphic system, that remark is self-explanatory for the geometrical meaning of DVI mean of a function. Below is a detailed explanation of the meaning with a special instance.

In Illustration \ref{illus:ExampleOfIMAS} a function $f$'s graph is plotted on the dual-isomorphic system. Now assume it's an ordinary function. Its graph is congruent to its DVI-function $\varphi$'s graph on the referential aux system(system B, which is not drawn there). Let $x_1=1$, $x_2=2$ be 2 points on the $x$ axis. Along with the point $m$ on the $x$ axis, there are isomorphic numbers $g(x_1)=e$, $g(m)=e^m$ and $g(x_2)=e^2$ on the aux. axis of system B. For easier understanding, we assume that right at $e^m$ on the aux. axis, there we just have $M_\varphi=\varphi(e^m)=\big(f(m)\big)^p$ being the mean value of $\varphi$ on $[e,~e^2]$ and it falls into $V=\{v:v>0\}$. (With Theorem \ref{thm:IsomeanExist1} it is possible, esp. when $f$ is continuous on [1,2]). Then with the above-mentioned congruency, the point $(e^m, M_\varphi)$ just has its congruent point$(g^{-1}(e^m), ~h^{-1}(M_\varphi))$ on the graph of $f$ which is ($m, f(m)$). Finally on the $y$-axis of the dual-isomorphic system, we can find this DVI mean value $f(m)$ of $f$ on $[1,2]$.

However in case that $f$ is not continuous, e.g. with some jump discontinuity such that both range $M$ of $f$ and set $h(M)$ are not intervals, it's still possible that an eligible $M_\varphi$ exists and falls out of $h(M)$ but in $V$, that will still end up with a DVI mean of $f$ being found on the $y$-axis of the dual-isomorphic system.

\subsubsection{Geometrical meaning of Cauchy's mean value theorem}
Still in Illustration \ref{illus:ExampleOfIMAS} a function $f(x)$'s graph is plotted on the dual-isomorphic system. Let for example $x_1=1$, $x_2=2$ be 2 points on the $x$-axis. Here $h$, $g$ are continuous real functions on $[x_1,x_2]$ which are differentiable in $(x_1,x_2)$, and $g'(x)\ne 0$ for $x\in(x_1,x_2)$. According to Theorem \ref{thm:CauchyMeanVal} there is a point  $m\in(x_1,x_2)$ such that
\begin{equation}
\frac{h'(m)}{g'(m)}=\frac{h(x_2)-h(x_1)}{g(x_2)-g(x_1)}~.
\end{equation}
Here especially with $f(x)=x$ plotted, the geometrical meaning for the Cauchy's mean value Theorem \ref{thm:CauchyMeanVal} is easy to explain:

For any $m\in(x_1,x_2)$, there is a point $(m,m)$ on the graph of $f(x)=x$ and at that point $\frac{h'(m)}{g'(m)}$ is just the dual-isomorphic derivative value $\frac{(h\circ f(x))'}{g'(x)}\bigl|_{x=m}$. This value's geometrical meaning is interpreted as the slope of tangent line of the graph at the point $(m,m)$. On the other hand, $\frac{h(x_2)-h(x_1)}{g(x_2)-g(x_1)}$ is interpreted as the slope of a secant line that passes both $(x_1,x_1), (x_2,x_2)$ on the graph of $f$.

So the equation in that Theorem literally says that under these conditions, there is always a point $m\in(x_1,x_2)$ such that the tangent line there is parallel to the secant line passing $(x_1,x_1), (x_2,x_2)$.

This is an explanation of the geometrical meaning of the Cauchy's mean value theorem in a dual-isomorphic system.

\subsection{About convexity on the dual-isomorphic system and isomorphic means}
According to Lemma \ref{lem:DiffCriDVIC}, the monotonicity of $g, h$ and $(h\circ f)'/g'$ can be used  to determine the convexity of $f$ on the dual-isomorphic system. In Illustration \ref{illus:ExampleOfIMAS}, $y=f(x)=x (x>0)$,
\begin{eqnarray}
    \frac{(h\circ f)'}{g'} &=& \frac{(x^p)'}{(e^x)'}, \nonumber \\
    \biggl(\frac{(h\circ f)'}{g'}\biggl)' &=& \biggl(\frac{px^{p-1}}{e^x}\biggl)' \nonumber \\
    &=& \frac{p(p-1)x^{p-2}e^x-px^{p-1}e^x}{e^{2x}} \nonumber \\
    &=& p(p-1-x)\frac{x^{p-2}}{e^x}.
\end{eqnarray}

Now consider the sign of above derivative is determined by (same as) sign of $p(p-1-x)$:
\begin{enumerate}
  \item[1).] If $p>0$ and $x>p-1$, then $p(p-1-x)<0$, $(h\circ  f)'/g'$  is decreasing. As $g$ is increasing, $h$ is increasing, in this case $f$ is convex to the upper in the dual-isomorphic system. For example in Illustration \ref{illus:ExampleOfIMAS} $p=2$, when $x>1$ the curve is convex to the upper. Meanwhile according to Corollary \ref{cor:IsoMeanComp} we get the inequality between two types of isomorphic means
        \begin{equation*}
            \hspace{-1.0cm}\overline{x_i,p_R}|_g\geq \overline{x_i,p_R}|_h~(i\in\{1,2,\ldots ,n\}, ~n\geq2),
        \end{equation*}
and this can be observed in the illustration by using the method in section \ref{sec:GraphCompMean}.
  \item[2).] If $p>0$ and $0<x<p-1$, then $p(p-1-x)>0$, $(h\circ f)'/g'$  is increasing. As $g$ is increasing, $h$ is increasing, $f$ is convex to the lower, and we've got
        \begin{equation*}
            \hspace{-1.0cm}\overline{x_i,p_R}|_g\leq \overline{x_i,p_R}|_h~(i\in\{1,2,\ldots ,n\}, ~n\geq2).
        \end{equation*}
These can be seen in the illustration for $p=2, ~0<x<1$.
  \item[3).] If $p<0$, then $p(p-1-x)>0$, and $(h\circ  f)'/g'$ is increasing. Since $g$ is increasing, $h$ is decreasing, in this case $f$ is convex to the upper, and we've got
        \begin{equation*}
            \hspace{-1.0cm}\overline{x_i,p_R}|_g\geq \overline{x_i,p_R}|_h~(i\in\{1,2,\ldots ,n\}, ~n\geq2).
        \end{equation*}
\end{enumerate}

\subsection{About 2-dimensional isomorphic convex set}
In Illustration \ref{illus:ExampleOfIMAS}, the chamfered area segregated by $y=x(x>1),~x=1$ and the $x$ axis is a 2-dimensional isomorphic convex set generated by $g, h$. It's literally the set $T=\{(x, y)\colon  y<x, ~x>1, ~y>0\}$, and it satisfies: For arbitrary $p_1(x_1, y_1), ~p_2(x_2, y_2)\in T$, and arbitrary $\lambda_1 ,\lambda_2\in (0,1)$ satisfying $\lambda_1 +\lambda_2=1$, the point $p=\big(\ln (\lambda_1e^{x_1}+\lambda_2e^{x_2}),~\sqrt{\lambda_1y_1^2+\lambda_2y_2^2}\big)\in T$, i.e.
\begin{equation*}
    \sqrt{\lambda_1y_1^2+\lambda_2y_2^2} < \ln (\lambda_1e^{x_1}+\lambda_2e^{x_2}).
\end{equation*}

\subsection{About isomorphic mean of a function}
The isomorphic mean class I of $f$ on interval $[a, b] \subseteq D$ generated by $h$ is:
\begin{eqnarray}
    \isomeanvalue{f(x)}{[a,b]}{h}
    &=& h^{-1}\Bigl[\frac{1}{b-a}\int_a^bh\big(f(x)\big)\mathrm{d}x\Bigl] \nonumber\\
    &=& \biggl[\frac{1}{b-a}\int_a^bf^p(x)\mathrm{d}x\biggl]^{\frac 1p}.
\end{eqnarray}
Taking $y=f(x)=x, ~p=2, ~[a, b]= [1, 2]$, ~and according to (\ref{equ:P-orderMean}),
\begin{eqnarray}
    \isomeanvalue{f(x)}{[a,b]}{y^p}
    &=& \biggl[\frac{b^{(p+1)}-a^{(p+1)}}{(p+1)(b-a)}\biggl]^{\frac 1p} \nonumber\\
    &=& \biggl[ \frac{b^3-a^3}{3(b-a)}\biggl]^{\frac12} \\
    &=& \sqrt{\frac73}~(\approx 1.5275). \nonumber
\end{eqnarray}

The isomorphic mean class II of $f$ on interval $[a, b] \subseteq D$ generated by $g$ is:
\begin{eqnarray}
    \isomeanvalueII{f(x)}{[a,b]}{g}
    &=& \frac{1}{g(b)-g(a)}\int_{g(a)}^{g(b)}f\big(g^{-1}(u)\big)\mathrm{d}u \nonumber\\
    &=& \frac{1}{e^b-e^a}\int_{e^a}^{e^b}f(\ln u)\mathrm{d}u.
\end{eqnarray}
Taking $f(x)=x, ~[a, b]=[1, 2]$,
\begin{eqnarray}
    \isomeanvalueII{f(x)}{[a,b]}{e^x} &=& \frac{1}{e^2-e^1}\int_{e^1}^{e^2}\ln u\mathrm{d}u \nonumber\\
    &=& \frac{1}{e^2-e}(u\ln u-u)\bigl|_e^{e^2}\\
    &=& \frac{e^2}{e^2-e}~(\approx 1.5820). \nonumber
\end{eqnarray}

The corresponding isomorphic mean class III and class IV can also be formulated and computed, which are not discussed here.

\section{About the logarithmic-mapping-generated systems}
In previous discussion of various IMAS topics, the logarithmic generator mapping e.g.~$x\mapsto \ln x:\mathbb{R}^+\to\mathbb{R}$ (It can also be log. function with other base e.g. $\lg x$) are frequently used for instances , which yield some MA concepts as special cases of the IMAS topics. In fact the IMAS instance systems generated by logarithmic mapping(s) are most typical ones. Depending on topics, the IMAS systems may be one-dimensional, i.e. in the isomorphic number-axis generated by logarithmic mapping, or two-dimensional, i.e. in the dual-isomorphic system which has 1 or 2 logarithmic mapping(s) as its generator-mapping(s).  The instances are summarized in the following table:

\begin{center}
\begin{tabular}{l|l|l}\hline
Topics                              &       As special case of IMAS topics &   Sections \\ \hline \hline
Multiplication operation            &       Isomorphic addition         &   $\ref{subsec:IsoAdd}$ \\ \hline
Division operation                  &       Isomorphic subtraction      &   $\ref{subsec:IsoSubtra}$ \\ \hline
Power operation                     &       Isomorphic multiplication   &   $\ref{subsec:IsoMultip}$ \\ \hline
Extraction operation                &       Isomorphic division type I  &   $\ref{subsubsec:IsoDivT1}$ \\ \hline
Logarithmic operation               &       Isomorphic division type II &   $\ref{subsubsec:IsoDivT2}$ \\ \hline
Geometric mean                      &       Isomorphic mean of numbers  &   $\ref{sec:IsoWghtMean}$ \\ \hline
Geometric mean of a function        &       Isomorphic mean class I of a function & $\ref{subsubsec:GeoMeanFunc}$ \\ \hline
Elastic mean of a function*         &       Isomorphic mean class II of a function& $\ref{subsubsec:SpecCaseIsoMeanFuncT2ElasMeanFun}$ \\ \hline
Exponential derivative              &       Metrical dual-isomorphic derivative & $\ref{subsubsec:SpecCaseGenDDer}$ \\ \hline
Bigeometric derivative              &       Metrical dual-isomorphic derivative & $\ref{subsubsec:SpecCaseGenDDer}$ \\ \hline
Geometric integral                  &       Isomorphic integral type I  &   $\ref{subsec:IsoIntegT1}$ \\ \hline
Elastic integral*                   &       Isomorphic integral type II &   $\ref{subsec:IsoIntegT2}$ \\ \hline
Elasticity of function              &       Dual-isomorphic derivative  &   $\ref{subsec:Elasticity}$ \\ \hline
2-D geometrically convex set*       &       2-D isomorphic convex set   &   $\ref{subsubsec:2DGeoConSet}$ \\ \hline
Geometrically convex function       &       Dual-variable-isomorphic convex function   &   $\ref{subsec:GeometricConvFunc}$ \\ \hline
\end{tabular}
\footnotesize{\newline * Defined a name in this paper~~~~~~~~~~~~~~~~~~~~~~~~~~~~~~~~~~~~~~~~~~~~~~~~~~~~~~~~~~~~~~~~~~~~~~~~~~~~~~~~~~~~~~~~~~~~~}
\begin{tables}\label{tab:GeometricGenSysSummary}
\centering\normalsize{The instances in logarithmic-mapping-generated systems}
\end{tables}
\end{center}


\section{Conclusion \& Vision}
Thus far, based on the 5 basic concepts introduced in the beginning, a series of new mathematical analysis concepts have been introduced to form a primitive framework of a possibly upgraded MA system –- the Isomorphic Mathematical Analysis System. In a sense, the IMAS is a new perspective of existing MA system, in a way distorted but with the order maintained by ``isomorphism'', which should also have practical meanings when applying to the real world. The IMAS can be deemed as a generalized MA system and will be instantiated to classical MA when identity generator mappings are considered.

Furthermore as the factor of generator mappings or isomorphic frame is involved, the new IMAS could bring about new classes of problems which do not have counterparts in normal MA system, e.g. graphical comparison of quasi-arithmetic means, mean value and convexity of function $f(x)=x$, densities on a number-axis or a coordinate system, point of singular density etc.

Finally, in addition to each foregoing discussed topics, the IMAS are further exploitable –- there could be problems and topics such as \begin{itemize}
  \item The curvature of function on the dual-isomorphic system;
  \item The isomorphic coordinate system of $n$-dimension and the IMAS on them;
  \item More detailed coverage of ``Isomorphic calculus'' topics which could develop directly from isomorphic operations and have their physical meanings;
  \item Relating the unevenly distributed coordinate space with possible applications in the realm of Physics $\ldots$ etc
\end{itemize}
In conclusion, it is possible to introduce more topics by applying the general idea of the IMAS to more MA concepts, to make this IMAS richer and better established.




\end{document}